\pdfoutput=1
\documentclass{elsarticle}
\usepackage[margin=1in]{geometry}
\usepackage{algorithm}
\usepackage{algpseudocode}
\usepackage{url}

\providecommand{\headers}[2]{}

\makeatletter
\newcommand{\LoadPackageUnlessAcm}[2][]{%
    \@ifclassloaded{acmart}{%
        \PackageInfo{preamble}{Skipping package `#2' because acmart class is used}%
    }{%
        \if\relax\detokenize{#1}\relax
            \usepackage{#2}%
        \else
            \usepackage[#1]{#2}%
        \fi
    }%
}
\makeatother

\LoadPackageUnlessAcm{caption}
\LoadPackageUnlessAcm{subcaption}
\makeatletter
\@ifclassloaded{elsarticle}{}{\LoadPackageUnlessAcm{authblk}}
\makeatother
\LoadPackageUnlessAcm{hyperref}

\makeatletter
\@ifclassloaded{acmart}{}{
    \captionsetup{
        font=small,
        labelfont=bf,
        labelsep=period,
        justification=justified,
        singlelinecheck=false,
        skip=8pt
    }
}
\makeatother

\usepackage{etoolbox}
\AtBeginEnvironment{tabular}{\small}

\usepackage{amsthm}

\usepackage{amstext}
\usepackage{stmaryrd}
\usepackage{algpseudocode}

\usepackage{units}
\usepackage{cancel}
\usepackage{graphicx}

\usepackage{xcolor}
\usepackage{array}
\usepackage{verbatim}
\usepackage{booktabs}
\usepackage{makecell}
\usepackage{calc}
\usepackage{textcomp}
\usepackage{multirow}
\usepackage{tikz,tikzscale}

\usepackage{pgfplots}
\usepackage{pgfplotstable}
\usepackage{pgfkeys}
\pgfplotsset{compat=1.18}

\definecolor{deepgreen}{RGB}{0,100,0}

\newif\ifrunscripts
\runscriptstrue
\definecolor{degQthree}{RGB}{140,81,10}
\definecolor{degQseven}{RGB}{230,159,0}
\pgfplotsset{
    p2/.style={teal, mark=*},
    p3/.style={degQthree, mark=triangle*},
    p4/.style={red, mark=pentagon*},
    p5/.style={violet, mark=otimes*},
    p7/.style={degQseven, mark=square*},
}

\newcommand{\plotwidth}{0.49\columnwidth}
\newcommand{\plotheight}{0.4\columnwidth}

\pgfplotscreateplotcyclelist{paper markers}{
    {teal,    solid, mark=*},
    {orange,   solid, mark=triangle*},
    {blue,     solid, mark=square*},
    {red,      solid, mark=pentagon*},
    {violet,     solid, mark=otimes*},
    {brown,    solid, mark=halfcircle*},
    {black,   solid, mark=diamond*},
}

\colorlet{BarColorOne}{orange!60!white}
\colorlet{BarColorTwo}{teal!80!white}
\colorlet{BarColorThree}{red!80!black}
\colorlet{BarColorFour}{blue!60!white}
\colorlet{BarColorFive}{brown!70!black}

\pgfplotsset{legend swatch/.style={area legend, draw=none}}

\pgfplotsset{
    paperplot/.style={
            width=\plotwidth,
            height=\plotheight,
            cycle list name=paper markers,
            every axis/.append style={font=\small},
            grid=major,
            thick
        }
}

\pgfplotsset{
    longplot/.style={
            width=0.8\columnwidth,
            height=\plotheight,
            cycle list name=paper markers,
            every axis/.append style={font=\small},
            grid=major,
            thick
        }
}

\usepackage[all]{xy}
\usepackage{mathtools}
\usepackage{enumitem}
\usepackage{placeins}
\usetikzlibrary{calc}
\usetikzlibrary{shadings}
\usepgfplotslibrary{fillbetween}

\usepackage[mode=multiuser,status=draft]{fixme}
\fxsetup{innerlayout=inline}
\fxsetup{targetlayout=colorcb}
\fxusetheme{color}
\FXRegisterAuthor{mw}{envmw}{MW}

\definecolor{apply}{rgb}{0.3,.7,0.}
\definecolor{applyface}{rgb}{.9,.95,0.}
\definecolor{invert}{rgb} {0.5,0.,1.}
\colorlet{invertface}{invert!60!red}

\newif\ifprintversion
\printversiontrue

\usepackage[colorinlistoftodos,prependcaption,textsize=small]{todonotes}
\ifprintversion
    \setuptodonotes{
        backgroundcolor=white,
        linecolor=black,
        bordercolor=black,
        textcolor=black
    }
\fi

\def\mesh{\mathbb M}
\def\cell{K}
\def\R{\mathbb R}
\def\Q{\mathbb Q}

\def\vertex{X}

\usepackage{colortbl}

\usepackage{booktabs}
\usepackage[font=footnotesize,position=b]{subcaption}
\usepackage{adjustbox}

\makeatletter
\def\pgfplots@getautoplotspec into#1{%
        \begingroup
        \let#1=\pgfutil@empty
        \pgfkeysgetvalue{/pgfplots/cycle multi list/@dim}\pgfplots@cycle@dim
        \let\pgfplots@listindex=\pgfplots@numplots
        \pgfkeysgetvalue{/pgfplots/cycle list set}\pgfplots@listindex@set
        \ifx\pgfplots@listindex@set\pgfutil@empty
        \else
            \c@pgf@counta=\pgfplots@listindex
            \c@pgf@countb=\pgfplots@listindex@set
            \advance\c@pgf@countb by -\c@pgf@counta
            \globaldefs=1\relax
            \edef\setshift{%
                \noexpand\pgfkeys{
                    /pgfplots/cycle list shift=\the\c@pgf@countb,
                    /pgfplots/cycle list set=
                }
            }%
            \setshift%
            \globaldefs=0\relax
        \fi
        \pgfkeysgetvalue{/pgfplots/cycle list shift}\pgfplots@listindex@shift
        \ifx\pgfplots@listindex@shift\pgfutil@empty
        \else
            \c@pgf@counta=\pgfplots@listindex\relax
            \advance\c@pgf@counta by\pgfplots@listindex@shift\relax
            \ifnum\c@pgf@counta<0
                \c@pgf@counta=-\c@pgf@counta
            \fi
            \edef\pgfplots@listindex{\the\c@pgf@counta}%
        \fi
        \ifnum\pgfplots@cycle@dim>0
            \c@pgf@counta=\pgfplots@cycle@dim\relax
            \c@pgf@countb=\pgfplots@listindex\relax
            \advance\c@pgf@counta by-1
            \pgfplotsloop{%
                \ifnum\c@pgf@counta<0
                    \pgfplotsloopcontinuefalse
                \else
                    \pgfplotsloopcontinuetrue
                \fi
            }{%
                \pgfkeysgetvalue{/pgfplots/cycle multi list/@N\the\c@pgf@counta}\pgfplots@cycle@N
                \pgfplotsmathmodint{\c@pgf@countb}{\pgfplots@cycle@N}%
                \divide\c@pgf@countb by \pgfplots@cycle@N\relax
                \expandafter\pgfplots@getautoplotspec@
                \csname pgfp@cyclist@/pgfplots/cycle multi list/@list\the\c@pgf@counta @\endcsname
                {\pgfplots@cycle@N}%
                {\pgfmathresult}%
                \t@pgfplots@toka=\expandafter{#1,}%
                \t@pgfplots@tokb=\expandafter{\pgfplotsretval}%
                \edef#1{\the\t@pgfplots@toka\the\t@pgfplots@tokb}%
                \advance\c@pgf@counta by-1
            }%
        \else
            \pgfplotslistsize\autoplotspeclist\to\c@pgf@countd

            \pgfplots@getautoplotspec@{\autoplotspeclist}{\c@pgf@countd}{\pgfplots@listindex}%
            \let#1=\pgfplotsretval
        \fi
        \pgfmath@smuggleone#1%
        \endgroup
    }

\pgfplotsset{
    cycle list set/.initial=
}
\makeatother

\AtBeginEnvironment{table}{\footnotesize}
\AtBeginEnvironment{table*}{\footnotesize}
 \usetikzlibrary{patterns,arrows.meta,fit,backgrounds}

\usepackage{lipsum}
\usepackage{amsfonts}
\usepackage{graphicx}
\usepackage{epstopdf}
\usepackage{amsmath,amssymb,amsthm}
\newtheorem{lemma}{Lemma}
\newtheorem{theorem}{Theorem}
\newtheorem{corollary}{Corollary}
\theoremstyle{definition}
\newtheorem{assumption}{Assumption}
\newtheorem{proposition}{Proposition}
\theoremstyle{remark}
\newtheorem{remark}{Remark}
\theoremstyle{plain}
\ifpdf
    \DeclareGraphicsExtensions{.eps,.pdf,.png,.jpg}
\else
    \DeclareGraphicsExtensions{.eps}
\fi

\headers{Short title }{Michał Wichrowski}

\makeatletter
\@ifclassloaded{elsarticle}{}{
    \title{
        A Multigrid Method for CutFEM and its Convergence
    }
    \author{Michał Wichrowski$^0$}
}
\makeatother

\usepackage{amsopn}

\newcommand{\Keywords}[1]{\begin{keyword}#1\end{keyword}}
\newcommand{\AMS}[1]{}

\title{Shedding the ghost weight: unconditionally elliptic CutFEM and high-order multigrid}

\author{Micha\l{} Wichrowski\corref{cor1}}
\cortext[cor1]{Corresponding author. \texttt{mwichro@mimuw.edu.pl}}

\begin{document}

\begin{abstract}
In cut finite element methods, the ghost penalty stabilizes arbitrarily small cut cells, but its weight $\gamma$
pollutes the discretization and degrades multigrid solvers. Both the error and the solver improve as the weight is
reduced. For the symmetric Nitsche formulation, however, this reduction is capped: coercivity fails below a
cut-dependent floor. Switching to the nonsymmetric Nitsche formulation yields a method that is unconditionally elliptic
for every positive weight, cut configuration, and mesh level. The weight becomes a free design variable that may be
made arbitrarily small. Theoretically, the a priori error constants that depend on the decreasing weight are bounded,
saturating below a geometry-determined threshold on any fixed configuration. Numerically, reducing the weight lowers
the error and largely eliminates the degree-dependence of the iteration counts in a vertex-patch multigrid method with
full-residual patches. The price of this formulation is adjoint inconsistency: the energy error converges at the
optimal rate $h^p$, while the $L^2$ estimate yields $h^{p+1/2}$. In practice, this half-order loss appears only for
deliberately adversarial data; smooth solutions retain optimal convergence.
 \end{abstract}


\maketitle

\newcommand{\mgpaper}{\cite{Wichrowski26CutFEMMG}}

\newcommand{\phd}[1]{\textcolor{red}{#1}}

\section{Introduction}
\label{sec:intro}
High-order discretizations call for solvers whose cost does not grow with the polynomial degree $p$. On fitted
Cartesian meshes the vertex-patch multigrid method has this property~\cite{Pavarino94additive}, but with the cut finite
element method (CutFEM)~\cite{BurmanClausHansboLarsonMassing15} its iteration counts grow with
$p$~\cite{CuiKanschat25CutFEM,Wichrowski26CutFEMMG}. The convergence analysis of~\mgpaper{} traces this growth to the
ghost penalty~\cite{burman2010ghost} that stabilizes the small cut cells, and its bounds improve as the ghost weight
$\gamma$ is reduced, or \emph{shed}. Under the symmetric Nitsche formulation the weight hits a floor: coercivity holds
only above a cut-dependent value of $\gamma$. The solver favors a small weight; stability requires a large one. The
symmetric form cannot give both.

In this paper we remove the floor by changing the sign of one term of the Nitsche form. The resulting nonsymmetric
formulation is \emph{unconditionally elliptic}: elliptic for every positive weight, every cut configuration, and every
mesh level. The ghost weight thus becomes a free design variable that may be made arbitrarily small, though never zero:
a positive weight keeps the local problems of the multigrid smoother away from singular on small cuts. A smaller weight
also benefits the discretization itself, since the ghost penalty is only weakly consistent. We then develop the
discretization theory for every positive weight: well-posedness and a priori error estimates whose constants do not
change under mesh refinement for every degree-dependent weight held fixed, with the remaining $\gamma$-dependence
confined to explicit constants that are bounded by geometry-determined caps on a fixed configuration. The price is
adjoint consistency, which costs half an order in the $L^2$ estimate. Shedding the ghost weight improves both the
discretization and the solver at a moderate price; we study the effect on the solver numerically.

In CutFEM the domain is embedded in a fixed Cartesian background mesh, the boundary cuts freely through the cells, and
the Dirichlet condition is imposed weakly by Nitsche's method~\cite{nitsche1971variationsprinzip,HansboHansbo02}. The
geometric freedom has a known price: cells retaining an arbitrarily small physical fraction destroy discrete coercivity
and conditioning. The established remedy is the ghost penalty, a penalty on jumps of normal derivatives across faces of
the boundary zone, which lets badly cut cells borrow stability from neighbors with substantial physical support and
restores stability independently of the position of the
boundary~\cite{BurmanHansbo12,MassingLarsonLogg14,dePrenterVerhooselvanBrummelenLarsonBadia23}.

Cell aggregation instead merges badly cut cells into well-supported aggregates and removes the instability by
construction~\cite{BadiaVerdugoMartin18}; the two approaches are linked, aggregation acting as a ghost penalty of
infinite weight~\cite{BadiaNeivaVerdugo22}. Discrete extension operators build the finite element space itself from
well-supported degrees of freedom and dispense with stabilization~\cite{BurmanHansboLarson22Ext}. The finite cell
method~\cite{ParvizianDusterRank07} fills the exterior with a soft fictitious material, and the shifted boundary
method~\cite{MainScovazzi18} moves the boundary condition to a surrogate boundary aligned with the mesh. Each of these
trades the cut-cell problem for a different discrete structure; the ghost penalty instead keeps the standard space and
mesh connectivity, which is what a matrix-free multigrid hierarchy on the Cartesian background grid
exploits~\cite{CuiKanschat25CutFEM}, and it exposes its stabilization strength as a single scalar $\gamma$.

The ghost penalty is not free either: the face penalty forces the solution towards a single polynomial across
neighboring cut cells, the method locks, and the error becomes sensitive to the
weight~\cite{BadiaNeivaVerdugo22,wichrowski2027unified}. Locking-free variants change what is
penalized~\cite{BurmanHansboLarson26Locking}. In Stokes flow, the standard ghost penalty on the pressure destroys the
pointwise divergence-free velocity of $H(\mathrm{div})$-conforming elements, and the stabilization has to be
redesigned~\cite{FrachonHansboNilssonZahedi24,FrachonNilssonZahedi24}. In finite-strain elasticity the error grows with
the ghost weight~\cite{wichrowski2027unified}. In~\mgpaper{} the growth of the multigrid iteration counts with $p$
enters through a ghost-energy constant on the boundary strip that grows superexponentially in $p$, and redistributing
the weight among the derivative orders does not remove it from the bounds. The conventional fixed $O(1)$
weight~\cite{burman2010ghost} is thus a compromise, and two refinements apply the stabilization only where the cut
requires it. Frei et al.~\cite{FreiKnokeSteinbachWenskeWick24,FreiKnokeSteinbachWenskeWick25} scale the weight with the
cut fraction of each cell, below the conventional value on well-cut cells and above it on badly cut ones, which reduces
spurious fluctuations of the solution on cut cells. Larson and Zahedi~\cite{LarsonZahedi23Macro} restrict the ghost
penalty to the interior of macroelements. Both keep the stabilization on the badly cut cells, which is where the
symmetric floor sits.

How small, then, may the weight be taken? For the symmetric formulation, not very: coercivity requires the ghost term
to control the Nitsche flux on badly cut cells, which bounds $\gamma$ from below in terms of the cut configuration, and
no schedule can pass below this floor without losing well-posedness. On top of that, in elasticity the flux carries the
material tangent, so the floor, and the Nitsche penalty with it, grow with the stiffness~\cite{wichrowski2027unified}.
At finite strain near a corner singularity, where the deformation gradient is unbounded, so are the required penalties.

The floor, however, belongs to the symmetric formulation, not to the problem: flip the sign of one Nitsche consistency
term, and it is gone. Freund and Stenberg~\cite{FreundStenberg95} introduced the nonsymmetric variant and observed that
it is stable for every positive penalty, at the cost of symmetry.\footnote{For the Poisson problem the original paper
    notes that, with the symmetry lost, the variant ``does not seem to be too
    attractive''~\cite[p.~330]{FreundStenberg95}.} In the discontinuous Galerkin family the same variant appeared as the
method of Oden, Babu\v ska and Baumann~\cite{OdenBabuskaBaumann98} and as the nonsymmetric interior penalty method,
whose a priori analysis by Rivi\`ere, Wheeler and Girault~\cite{RiviereWheelerGirault01} established the pattern the
present paper inherits: unconditional coercivity gives optimal energy rates, and adjoint inconsistency degrades the
$L^2$ rate; the unified framework of~\cite{arnold2002unified} places both variants in one family. The two consistency
terms cancel on the diagonal, so the form equals its symmetric part there, and ellipticity holds regardless of how the
flux behaves.

The sign change has also been used to dispense with the boundary penalty. Burman~\cite{Burman12PenaltyFree} showed that
on fitted meshes the nonsymmetric method remains stable without it, Boiveau and Burman~\cite{BoiveauBurman16} extended
this to elasticity, and Schillinger et al.~\cite{SchillingerHarariHsu16} advocated the penalty-free nonsymmetric method
for immersed discretizations because it spares the user a cut-dependent penalty estimate, reporting robust accuracy
across embedded-domain and coupling problems. In a divergence-free CutFEM for Stokes flow, Frachon, Nilsson and
Zahedi~\cite{FrachonNilssonZahedi24} chose the nonsymmetric form because the penalty parameter is then less critical.
We borrow the trick but aim it at the other parameter: we keep the boundary penalty. In the coercivity identity it
appears as a separate nonnegative term, so every $\tau_0>0$ is admissible, as Freund and Stenberg observed, and keeps
the symmetric part of the operator definite; what we remove is the lower bound on the \emph{ghost} weight. Because the
consistency terms cancel, the ghost weight appears in no lower bound, and well-posedness at every positive weight
follows from an algebraic identity and a connectivity argument, with no trace, inverse, or extension estimate;
ellipticity itself takes one line.
To my knowledge, the cancellation has not been used to free the ghost weight before.

Adjoint consistency is lost with the sign change, and with it half an order in the duality argument: the energy error
converges at the optimal rate $h^p$ at each fixed positive weight, while the $L^2$ estimate gives $h^{p+1/2}$. For the
nonsymmetric DG relative, suboptimal $L^2$ convergence at odd degrees was proved in~\cite{GuzmanRiviere09}. For fitted
nonsymmetric Nitsche methods, Chen, Liu and Zhang~\cite{ChenLiuZhang25} recently clarified the picture with a unified,
regularity-dependent $L^2$ estimate: smooth solutions converge at the optimal order, and the half order is lost only
near the low end of the regularity scale, which their piecewise-linear experiments confirm near the $H^2$ endpoint.
Their analysis tells us where to look on cut meshes. We observe the half-order loss only on a boundary-layer solution
modeled on theirs; there the estimate is sharp, at least for the lowest degrees, while smooth solutions converge at the
optimal order $p+1$.

All weight dependence of the error constants enters through a single flux-control constant with two branches:
cut-uniform at the price $\gamma^{-1/2}$, or weight-uniform at the price of the cut geometry, and on a fixed
configuration the constant saturates once $\gamma$ falls below a geometry-determined threshold. The weight also enters
the error through the explicit interpolation contribution of the ghost penalty, which carries a factor $\gamma^{1/2}$;
the ghost form is consistent when the solution extension has sufficient regularity.

A discretization is half a method; the other is an efficient solver for the systems it produces. In the unfitted
literature, it appears to be the smaller half; see the survey~\cite{dePrenterVerhooselvanBrummelenLarsonBadia23}.
Multigrid methods for unfitted \emph{interface} problems were developed in~\cite{LudescherGrossReusken20,Ludescher20},
and Gross and Reusken~\cite{GrossReusken23} analyze an overlapping domain-decomposition preconditioner in the same
setting. For CutFEM proper, the cell-based additive Schwarz method of~\cite{PrenterVerhooselBrummelen19} achieves
cut-independent but $h$-dependent bounds, lacking a coarse space. The vertex-patch multigrid method
of~\cite{CuiKanschat25CutFEM} builds on the fast tensor-product Schwarz smoothers of~\cite{WitteArndtKanschat21} and
the localized residual computations of~\cite{WichrowskiMunchKronbichlerKanschat25}, and its degree dependence is
explained by the analysis of~\mgpaper{}. All of these take the stabilization as given. This paper delivers a whole
method, and the reduced weight improves both halves.

We apply the vertex-patch multigrid method of~\mgpaper{} to the nonsymmetric system without modification. Its
full-residual patch spaces adapt the construction of~\cite{wichrowski2026SBMmg} to CutFEM: on cut patches they also
update every degree of freedom whose active support lies entirely inside the patch. Each patch solve restricts the same
nonsymmetric operator and the same global residual, and unconditional ellipticity makes every patch problem uniquely
solvable at every weight. The V-cycle preconditions GMRES, and we report iteration counts across mesh levels,
polynomial degrees, and ghost weights, for the symmetric and the nonsymmetric form. The convergence theory
of~\mgpaper{} covers the symmetric formulation above its floor; we do not extend it to the nonsymmetric form. Filling
this gap would combine a stable decomposition at a reduced weight with field-of-values bounds for
GMRES~\cite{EisenstatElmanSchultz83,Starke97fov}; we leave it open.
In the experiments, the two forms need nearly the same number of iterations at equal weights, but only the nonsymmetric
form may go below the floor. At its floor, the symmetric form still needs more iterations as the degree grows; at the
reduced weight, the degree has a much smaller influence on the iteration counts, up to $\Q_7$. The counts change little
under mesh refinement, in two and in three dimensions, and do not react to arbitrarily small cuts.

\begin{figure}[!ht]
    \centering
    \begin{tikzpicture}[
            res/.style={draw=black!70, rounded corners=2pt, fill=white, align=center, font=\footnotesize,
                    text width=3.0cm, inner sep=3pt, minimum height=0.95cm},
            hub/.style={res, line width=1.1pt, draw=blue!60!black, fill=blue!6},
            sec/.style={draw=none, fill=black!5, rounded corners=4pt, inner sep=5pt},
            seclabel/.style={font=\scriptsize\itshape, text=black!60, anchor=south west, inner sep=1pt},
            use/.style={-{Stealth[length=2.2mm]}, black!75, line width=0.6pt},
            interp/.style={use, dashed},
        ]
        \node[res] (tools) at (7.4,0) {Trace and extension\\ \eqref{eq:cuttrace}, \eqref{eq:extension}};
        \node[hub] (twobranch) at (11.1,0) {Two-branch bound \eqref{eq:twobranch}\\
            $\Xi_\gamma=\min\{\gamma^{-1},1+R\}$};

        \node[hub] (coer) at (0,-2.1) {Thm.~\ref{thm:coercivity}\\ $a_\gamma(v,v)=\|v\|_\gamma^2$};
        \node[res] (wellposed) at (3.7,-2.1) {Cor.~\ref{cor:wellposed}\\ well-posedness};
        \node[res] (cont) at (7.4,-2.1) {Thm.~\ref{thm:continuity}\\ $M_{\mathrm{geo}}$, $M_{\mathrm{wt}}(\gamma)$};
        \node[res] (active) at (11.1,-2.1) {Prop.~\ref{prop:activemesh}, Cor.~\ref{cor:nofloor}\\
            active-mesh control};

        \node[res] (ortho) at (0,-4.2) {Lem.~\ref{lem:orthogonality}, \ref{lem:interp}\\
            orthogonality, interpolation};
        \node[res] (energy) at (3.7,-4.2) {Thm.~\ref{thm:energy}\\ energy rate $h^p$};
        \node[res] (Mgamma) at (7.4,-4.2) {flux constant\\ $M_\gamma$ \eqref{eq:Mgamma}};
        \node[res] (ltwo) at (11.1,-4.2) {Thm.~\ref{thm:ltwo}\\ $L^2$ rate $h^{p+1/2}$};

        \node[res] (deriv) at (0,-6.7) {Lem.~\ref{lem:weight-derivative}, Cor.~\ref{cor:onedim}\\
            weight derivatives, monotonicity in 1D};
        \node[res] (transfer) at (3.7,-6.7) {Lem.~\ref{lem:transfer}\\ weight transfer};
        \node[res] (monotone) at (7.4,-6.7) {Cor.~\ref{cor:monotone}\\ monotone error bound};
        \node[res] (path) at (11.1,-6.7) {Lem.~\ref{lem:gamma-path}\\ weight path};

        \begin{scope}[on background layer]
            \node[sec, fit=(tools)(twobranch)] (s2) {};
            \node[sec, fit=(coer)(active)] (s3) {};
            \coordinate (s4low) at (7.4,-5.0);
            \node[sec, fit=(ortho)(ltwo)(s4low)] (s4) {};
            \node[sec, fit=(deriv)(path)] (s44) {};
        \end{scope}
        \node[seclabel] at (s2.north west) {Section~\ref{sec:setting}};
        \node[seclabel] at (s3.north west) {Section~\ref{sec:wellposed}};
        \node[seclabel] at (s4.north west) {Section~\ref{sec:error}};
        \node[seclabel] at (s44.north west) {Section~\ref{sec:weight-path}};

        \draw[use] (tools) -- (cont);
        \draw[use] (tools) -- (active);
        \draw[use] (twobranch) -- (cont);
        \draw[use] (twobranch) -- (active);
        \draw[use] (coer) -- (wellposed);
        \draw[use] (cont) -- (wellposed);

        \draw[use] (cont) -- (Mgamma);
        \draw[use] (coer) -- (energy);
        \draw[use] (ortho) -- (energy);
        \draw[use] (Mgamma) -- (energy);
        \draw[use] (energy.south) .. controls (6,-4.95) and (9,-4.95) .. (ltwo.south);
        \draw[use] (Mgamma) -- (ltwo);

        \draw[use] (energy.south) -- (monotone);
        \draw[use] (energy.south) -- (path);
        \draw[use] (path) -- (monotone);
        \draw[interp] (transfer) -- (monotone);
        \draw[use] (ortho) -- (deriv);

        \draw[use] (twobranch.east) -- ++(0.55,0) |- (path.east);
    \end{tikzpicture}
    \caption{Dependencies between the results. Solid arrows mean that a result is used in the proof of another; the
        dashed arrow means that Lemma~\ref{lem:transfer} is used to interpret Corollary~\ref{cor:monotone}. The two
        framed results are the hubs: the coercivity identity of Theorem~\ref{thm:coercivity} enters every proof from
        Section~\ref{sec:error} on, and only its first uses are drawn; the two-branch bound~\eqref{eq:twobranch} is the source of every weight-dependent constant.}
    \label{fig:proof-map}
\end{figure}
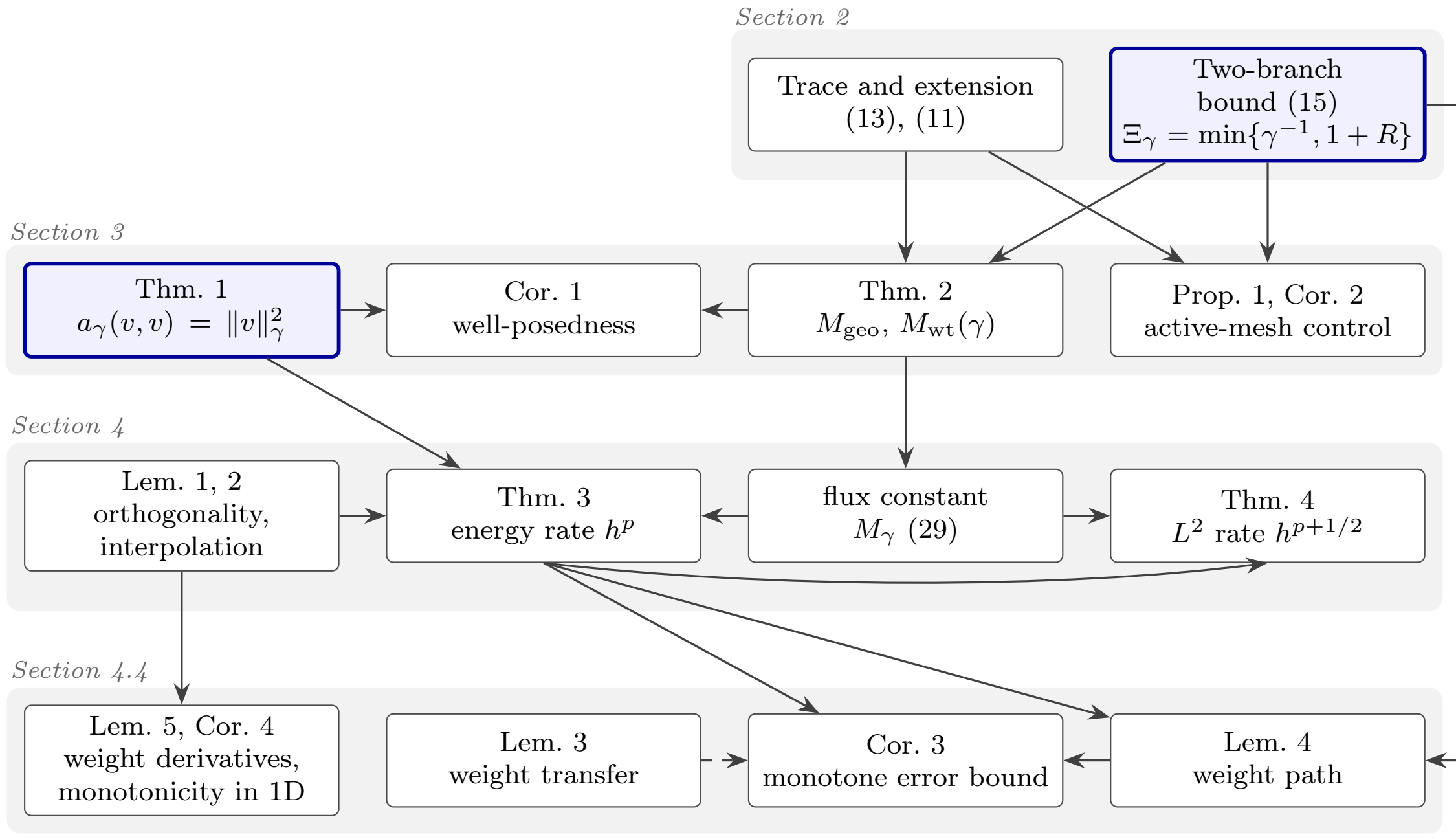

The paper is organized as follows. Section~\ref{sec:setting} fixes the discretization, the geometric assumption, and
the single constant through which the ghost weight enters the estimates. Section~\ref{sec:wellposed} proves
unconditional ellipticity, continuity, and well-posedness for every schedule. Section~\ref{sec:error} develops the a
priori error analysis, including the $L^2$ half-order bound, the exact transfer between ghost weights, and monotonicity
of an error bound along the weight path. Section~\ref{sec:numerics} reports the convergence of the discretization as
the ghost weight decreases, and Section~\ref{sec:multigrid} the multigrid method with full-residual patches and its
iteration counts. Section~\ref{sec:conclusions} concludes. Graph in Figure~\ref{fig:proof-map} provides a roadmap of
the theoretical results.

\section{Problem setting and formulation}
\label{sec:setting}
Let $\Omega\subset\R^d$, $d\in\{2,3\}$, be a bounded domain with boundary $\Gamma=\partial\Omega$. The model problem is
the Poisson equation with Dirichlet data,
\begin{gather}\label{eq:model}
    -\Delta u=f\quad\text{in }\Omega,
    \qquad
    u=g\quad\text{on }\Gamma,
\end{gather}
with data $f\in L^2(\Omega)$ and $g\in H^{1/2}(\Gamma)$, which makes the
problem uniquely solvable in $H^1(\Omega)$.

Let $\mesh_h$ be a Cartesian background mesh covering $\Omega$: every cell $\cell\in\mesh_h$ is an axis-parallel cube
of edge $h$, and every face is axis-parallel. The \emph{active mesh} consists of the cells meeting the physical domain,
\begin{gather}\label{eq:active}
    \mesh_{h,\Omega}:=\{\cell\in\mesh_h:\ |\cell\cap\Omega|>0\},
\end{gather}
and its union is the \emph{active domain} $\Omega_h\supseteq\Omega$. The
active mesh is face-connected: any two active cells can be joined by a
sequence of active cells in which consecutive cells share a
$(d-1)$-dimensional face. This follows from the connectedness of $\Omega$ and
the conformity of the Cartesian background mesh. The
discrete space $V_h^p$ is the continuous $\Q_p$ space of a fixed degree
$p\ge1$ on $\mesh_{h,\Omega}$, where
\begin{gather}\label{eq:Qp}
    \Q_p(\cell):=\operatorname{span}
    \bigl\{x_1^{\alpha_1}\cdots x_d^{\alpha_d}:\ 0\le\alpha_i\le p\bigr\}
\end{gather}
is the tensor-product space of degree at most $p$ \emph{in each variable
    separately}. Since the face normals are coordinate directions, a normal
derivative acts as a pure coordinate derivative and maps $\Q_p$ into itself,
which lets univariate polynomial inequalities be applied in the normal
direction with the tangential variables frozen. The nodal basis is the
Lagrange basis on the tensor-product Gauss--Lobatto nodes of each cell.
Functions of $V_h^p$ are cellwise polynomials on all of $\Omega_h$, so whether
a norm is taken over $\Omega$ or over $\Omega_h$ always matters and is
displayed in the subscript.

\begin{figure}[!ht]
    \centering
    \begin{tikzpicture}[scale=1.]
        \draw[step=1cm, gray, very thin] (-0.2,0) grid (5.2,3.2);

        \node at (-0.35,2.5) {$\Gamma$};
        \draw[line width=2pt, red, dashed] plot [smooth, tension=0.5] coordinates {(-0.2,2.3) (1,2.2) (2,2.5) (3,2.2)
                (4,1.5)
                (5.2,1.5)};

        \node at (2.5,0.7) {$\Omega$};
        \node at (3.2,2.7) {$\mathcal F_G$};

        \foreach \x/\y in {0/2, 1/2, 2/2, 3/1, 3/2, 4/1}
            {
                \fill[blue, opacity=0.3] (\x,\y) rectangle (\x+1,\y+1);
            }
        \foreach \x/\y in {0/0, 0/1, 2/1, 1/0, 1/1, 2/0, 3/0, 4/0}
            {
                \fill[green, opacity=0.3] (\x,\y) rectangle (\x+1,\y+1);
            }

        \foreach \x/\y in {0/2, 1/2, 2/2, 2/1, 3/1, -1/2, 4/1} {
                \draw[blue, line width=1.5pt] (\x+1,\y) -- (\x+1,\y+1);
            }

        \foreach \x/\y in {0/1, 1/1, 2/1, 3/0, 4/0, 3/1} {
                \draw[blue, line width=1.5pt] (\x,\y+1) -- (\x+1,\y+1);
            }

    \end{tikzpicture}
    \caption{The cut setting: the boundary $\Gamma$ (red dashed line) intersects a Cartesian background mesh,
        creating cut cells (blue) and interior cells (green). The ghost penalty~\eqref{eq:ghost} acts on the faces
        $\mathcal F_G$ (thick blue lines) of the boundary zone, which weakly glue a badly cut cell to neighbors with
        substantial physical support.}
    \label{fig:cut-setting}
\end{figure}
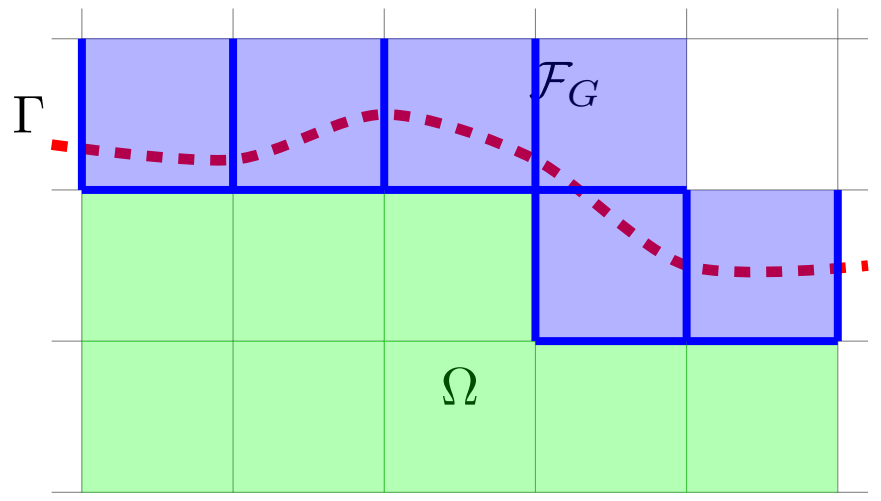

The ghost penalty acts on the faces of the boundary zone,
\begin{gather}\label{eq:ghostfaces}
    \mathcal F_G:=\{F=\partial\cell_1\cap\partial\cell_2:\
    \cell_1,\cell_2\in\mesh_{h,\Omega},\ \Gamma\cap(\cell_1\cup\cell_2)\ne\emptyset\},
\end{gather}
the interior faces of the active mesh with at least one cut neighbor
(Figure~\ref{fig:cut-setting}). Let
\begin{gather}\label{eq:ghost}
    g_h(u,v)
    =
    \sum_{F\in\mathcal F_G}
    \sum_{k=1}^p
    \frac{h_F^{2k-1}}{(k!)^2}
    \bigl(
    \llbracket\partial_{n_F}^ku\rrbracket,
    \llbracket\partial_{n_F}^kv\rrbracket
    \bigr)_F
\end{gather}
be the Burman ghost penalty~\cite{burman2010ghost}, with a single multiplier
$\gamma>0$ on all derivative orders. Here $h_F$ is the diameter of $F$,
$n_F$ a fixed unit normal of $F$, $\partial_{n_F}^k$ the $k$-th normal
derivative and $\llbracket\cdot\rrbracket$ the jump across $F$; on $\Gamma$,
$n$ denotes the unit normal pointing out of $\Omega$ and
$\partial_nv=n\cdot\nabla v$.

\begin{assumption}[geometric assumption]\label{ass:geometry}
    The boundary $\Gamma$ is $C^2$ and has positive reach $\delta_0>0$. If
    $\Gamma$ intersects a cell $\cell$, then $\Gamma\cap\cell$ is a single
    curve (surface) intersecting $\partial\cell$ exactly twice (in a single
    closed curve for $d=3$), each closed face is intersected at most once, and
    $\cell\setminus\Gamma$ has exactly two connected components. Moreover,
    $\Gamma$ contains no vertex of $\mesh_h$, and the mesh resolves the
    geometry,
    \begin{gather}\label{eq:resolution}
        c_R\,h\le\delta_0,
        \qquad
        c_R=c_R(d):=24d+10\sqrt d .
    \end{gather}
\end{assumption}

The assumption is the one of~\mgpaper{}. The $C^2$ regularity and the positive reach make the closest-point projection
onto $\Gamma$ well defined and $C^1$ in a tube of radius $\delta_0$; for a $C^2$ boundary this amounts to a curvature
bound together with a lower bound on the separation of distinct sheets. The cut-cell conditions, which ask the boundary
to cross each cut cell in a single piece, go back to Hansbo and Hansbo~\cite{HansboHansbo02}, and the exclusion of mesh
vertices is a generic-position assumption that can be enforced by an arbitrarily small perturbation of the level set.
The value of $c_R$ carries a safety margin for the geometric constructions it feeds and is not optimized. One
consequence is used repeatedly and stated separately: every active cell meets the physical domain in a set of positive
measure, $|\cell\cap\Omega|>0$, which holds by the definition~\eqref{eq:active} of the active mesh but would fail if an
extra layer of exterior cells were carried for the ghost chains.

For $\varepsilon\in\{-1,1\}$ define the Nitsche form
\begin{gather}\label{eq:aform}
    a_\gamma^\varepsilon(u,v)
    :=
    (\nabla u,\nabla v)_\Omega
    -(\partial_nu,v)_\Gamma
    +\varepsilon\,(u,\partial_nv)_\Gamma
    +\tau_h(u,v)_\Gamma
    +\gamma\, g_h(u,v).
\end{gather}
The choice $\varepsilon=-1$ gives the symmetric form $a_\gamma^{\mathrm{sym}}$, and $\varepsilon=1$ the nonsymmetric
form $a_\gamma$ studied in this paper. We use the standard boundary-penalty scaling
\begin{gather}\label{eq:tau}
    \tau_h=\tau_0\,p(p+1)/h,
    \qquad
    \tau_0>0\ \text{fixed},
\end{gather}
and the discrete problem: find $u_h\in V_h^p$ with
\begin{gather}\label{eq:discrete-problem}
    a_\gamma(u_h,v_h)=\ell(v_h)
    :=(f,v_h)_\Omega+(g,\partial_nv_h)_\Gamma+\tau_h(g,v_h)_\Gamma,
    \qquad v_h\in V_h^p,
\end{gather}
the load functional $\ell$ being the one that renders $a_\gamma$ consistent with the model problem~\eqref{eq:model}
(Lemma~\ref{lem:orthogonality}); it needs only $g\in L^2(\Gamma)$, and for $a_\gamma^{\mathrm{sym}}$ the sign of the
term $(g,\partial_nv_h)_\Gamma$ changes accordingly. For $\varepsilon=1$
the two consistency terms cancel on the diagonal, so the ghost weight is absent from every lower bound.

For $\gamma\ge0$, the symmetric form and \emph{graph norm} associated with~\eqref{eq:aform} are
\begin{gather}\label{eq:graphnorm}
    s_\gamma(v,v)
    :=
    \underbrace{\|\nabla v\|_{0,\Omega}^2+\tau_h\|v\|_{0,\Gamma}^2}_{=:\,n(v,v)}
    +\ \gamma\, g_h(v,v)
    \ =:\ \|v\|_\gamma^2 .
\end{gather}
The form $s_\gamma$ is symmetric. It is the form of a ghost-stabilized
\emph{Neumann} problem with a boundary mass term: no cut-dependent Nitsche
consistency term appears in it.

The analysis uses the following notation. For a measurable set $X\subseteq\R^d$, a face $F$, or a subset of $\Gamma$,
we write $(u,v)_X$ for the $L^2(X)$ inner product with respect to the Lebesgue or the surface measure,
$\|v\|_{0,X}:=(v,v)_X^{1/2}$ for the induced norm, and $|v|_{1,X}:=\|\nabla v\|_{0,X}$,
$\|v\|_{1,X}^2:=\|v\|_{0,X}^2+|v|_{1,X}^2$ for the $H^1$-seminorm and norm. Constants $c$ and $C$ are independent of
the mesh size and of the position of $\Gamma$ relative to the mesh, but may depend on the polynomial degree $p$, on the
shape of $\Gamma$ and on the boundary-penalty parameter $\tau_0$. Dependence on $p$ and on the ghost weight $\gamma$ is
displayed wherever it is tracked.

The extension estimate is used in its weight-resolved form~\cite{burman2010ghost,BurmanHansbo12}: there is $C_E(p)$,
independent of $h$ and the cut configuration, with
\begin{gather}\label{eq:extension}
    \|\nabla v_h\|_{0,\Omega_h}^2
    \le
    C_E(p)\max\{1,\gamma^{-1}\}\,
    \bigl(\|\nabla v_h\|_{0,\Omega}^2+\gamma\, g_h(v_h,v_h)\bigr),
    \qquad
    v_h\in V_h^p,
\end{gather}
for every $0<\gamma\le1$. Second, the converse facewise trace--Markov bound,
also in gradient form: there is $C_U(p)$, independent of $h$ and the cut, with
\begin{gather}\label{eq:converse}
    g_h(v_h,v_h)\ \le\ C_U(p)\,\|\nabla v_h\|_{0,\Omega_h}^2,
    \qquad
    v_h\in V_h^p.
\end{gather}
Only the gradient enters~\eqref{eq:converse}: in the $k$-th ghost term the
scaling $h_F^{2k-1}$ cancels against the $h^{-1}$ of the face trace and the
$h^{-2(k-1)}$ of the inverse inequality, leaving
$\|\nabla v_h\|_{0,\Omega_h}^2$ with no $L^2$ contribution. Third, the
cut-cell trace inequality~\cite{HansboHansbo02,MassingLarsonLogg14}: under
Assumption~\ref{ass:geometry} there is $C_{\mathrm{tr}}(p)$, independent of
$h$ and of the position of the cut, with
\begin{gather}\label{eq:cuttrace}
    h\,\|q\|_{0,\Gamma\cap\cell}^2
    \le
    C_{\mathrm{tr}}(p)^2\,\|q\|_{0,\cell}^2,
    \qquad
    q\in\Q_p(\cell),
\end{gather}
for every cut cell $\cell$. The right side is the norm on the \emph{full}
cell, not on its physical part, so~\eqref{eq:cuttrace} is an inverse trace
estimate from the cell to the embedded surface piece; its uniformity in the
cut position rests on the surface-measure bound
$|\Gamma\cap\cell|\le Ch^{d-1}$, which Assumption~\ref{ass:geometry}
provides through the positive reach and the resolution
condition~\eqref{eq:resolution}. Applied componentwise to $q=\partial_iv_h$ and summed over the
(disjoint) cut cells, it gives the global bound
$h^{1/2}\|\partial_nv_h\|_{0,\Gamma}\le
    C_{\mathrm{tr}}(p)\|\nabla v_h\|_{0,\Omega_h}$. The constants
$C_E(p)$, $C_U(p)$ and $C_{\mathrm{tr}}(p)$ can grow rapidly with the
degree. Their rate is not
needed below because no statement in this paper is claimed to be $p$-uniform.

Every weight-dependent constant in this paper is a multiple of a single quantity. Write
\begin{gather}\label{eq:shedfactor}
    R:=\sup_{v_h\ne0}\frac{g_h(v_h,v_h)}{n(v_h,v_h)}\ \ (<\infty),
    \qquad
    \Xi_\gamma:=\min\bigl\{\gamma^{-1},\,1+R\bigr\},
    \qquad 0<\gamma\le1,
\end{gather}
where $R$ is the \emph{amplification quotient}: over all discrete functions, the largest factor by which the ghost
form~\eqref{eq:ghost} amplifies the Neumann energy $n$ of~\eqref{eq:graphnorm}. The amplification comes from the cut
cells. On a small piece $\Omega\cap\cell$ a polynomial can carry little energy and still have large derivative jumps
on the ghost faces of $\cell$. The
elementary bound behind every later use of~\eqref{eq:shedfactor} is that, for
$0<\gamma\le1$,
\begin{gather}\label{eq:twobranch}
    \|\nabla v_h\|_{0,\Omega}^2+g_h(v_h,v_h)
    \ \le\
    \Xi_\gamma\,\|v_h\|_\gamma^2 ,
\end{gather}
since the left side is at most
$\gamma^{-1}(\|\nabla v_h\|_{0,\Omega}^2+\gamma g_h)$ and also at most
$n(v_h,v_h)+R\,n(v_h,v_h)\le(1+R)\|v_h\|_\gamma^2$. The two branches are the
two ways the unweighted ghost energy can be paid for: out of the weight, or
out of the geometry. The bound $\Xi_\gamma\le\gamma^{-1}$ holds uniformly in
the cut, and $R$ is finite for every fixed configuration but diverges as the
cut degenerates, $R\gtrsim\kappa^{-(2p-1)}$ in the slab model, so no
cut-uniform bound on $\Xi_\gamma$ better than $\gamma^{-1}$ exists. For a
\emph{fixed} mesh and cut, however, $\Xi_\gamma$ saturates at $1+R$ and stops
responding to a further decrease of the weight. This is the same dichotomy that governs
weight transfer in Lemma~\ref{lem:transfer}, and it is the
reason the estimates below cease to degrade once $\gamma$ falls below
$(1+R)^{-1}$.

\section{Stability and well-posedness}
\label{sec:wellposed}
\begin{theorem}[unconditional ellipticity]\label{thm:coercivity}
    Let $\tau_0>0$. Then, for every $\gamma>0$ and every $v_h\in V_h^p$,
    \begin{gather}\label{eq:coercivity-identity}
        a_\gamma(v_h,v_h)=s_\gamma(v_h,v_h)=\|v_h\|_\gamma^2
    \end{gather}
    and $\|\cdot\|_\gamma$ is a norm on $V_h^p$; the Neumann part $n$
    of~\eqref{eq:graphnorm} is itself positive definite. In particular the discrete
    inf-sup constant of $a_\gamma$ in the weight-dependent graph norm is at
    least one, for every $\gamma$, every cut configuration, and every $h$.
\end{theorem}

\begin{proof}
    The two Nitsche consistency terms in~\eqref{eq:aform} are exchanged with
    a sign under $u\leftrightarrow v$, so their sum is skew and vanishes on
    the diagonal; this gives~\eqref{eq:coercivity-identity}. For
    definiteness it suffices to show that $n(v_h,v_h)=0$ forces $v_h=0$,
    since $s_\gamma\ge n$. Suppose therefore $n(v_h,v_h)=0$. Then $\nabla v_h=0$ on
    $\Omega\cap\cell$ for every active cell $\cell$. Since $v_h|_\cell$ is a
    polynomial and $\Omega\cap\cell$ contains an open set, the polynomial
    identity theorem forces $v_h|_\cell$ to be constant on the whole cell;
    continuity and face-connectivity of the active mesh propagate one
    constant globally. Its boundary trace vanishes because
    $\tau_h\|v_h\|_{0,\Gamma}^2=0$ and $\tau_h>0$, so the constant is zero.
    The inf-sup statement follows from
    $\sup_{w}a_\gamma(v,w)/\|w\|_\gamma\ge
        a_\gamma(v,v)/\|v\|_\gamma=\|v\|_\gamma$.
\end{proof}

Theorem~\ref{thm:coercivity} uses no ghost-penalty transfer estimate, no inverse inequality and no trace inequality. It
is an algebraic identity plus a connectivity argument. The one thing the connectivity argument does need is that every
active cell meets $\Omega$ in a set of positive measure, $|\cell\cap\Omega|>0$ for all active $\cell$: the standard
active mesh $\{\cell:\ |\cell\cap\Omega|>0\}$ has this by definition, but a mesh carrying an extra layer of exterior
cells for the ghost chains does not, and such a cell supports a function whose graph norm is carried by the ghost term
alone, so that definiteness would degenerate with the weight. This is assumed throughout. Its content is that \emph{the
    ghost weight never appears in the lower bound}: there is no coercivity floor to respect, and hence no variational
obstruction to an arbitrarily aggressive schedule.

\begin{theorem}[continuity at the same weight, two branches]\label{thm:continuity}
    There is $C$ depending only on the shape regularity and $\tau_0$ such
    that, for $0<\gamma\le1$,
    \begin{gather}\label{eq:continuity}
        |a_\gamma(v_h,w_h)|
        \le
        C_a\,\|v_h\|_\gamma\|w_h\|_\gamma,
        \qquad
        C_a
        \le
        C\bigl(1+\min\{M_{\mathrm{geo}},\,M_{\mathrm{wt}}(\gamma)\}\bigr),
    \end{gather}
    with the two branches
    \begin{align}
        \label{eq:branch-geo}
        M_{\mathrm{geo}}
         & :=C_{\mathrm{Nik}}(p,\mathcal G),
         &                                                      & \text{(geometric branch, independent of $\gamma$),}
        \\
        \label{eq:branch-wt}
        M_{\mathrm{wt}}(\gamma)
         & :=C_{\mathrm{tr}}(p)\,C_E(p)^{1/2}\,\Xi_\gamma^{1/2}
        \le C_{\mathrm{tr}}(p)\,C_E(p)^{1/2}\,\gamma^{-1/2},
         &                                                      & \text{(weight branch, uniform in the cut).}
    \end{align}
    Here
    \begin{gather*}
        C_{\mathrm{Nik}}(p,\mathcal G)
        :=\max_{\cell:\,|\Gamma\cap\cell|>0}\ \sup_{0\ne q\in\Q_{p}(\cell)}
        \frac{h^{1/2}\|q\|_{0,\Gamma\cap\cell}}{\|q\|_{0,\cell\cap\Omega}}
    \end{gather*}
    is the trace--Nikolskii constant of the cut configuration.
\end{theorem}

Both branches bound the same quantity, the normal flux $h^{1/2}\|\partial_nv_h\|_{0,\Gamma}$ in the consistency terms,
and they pay for it in different currencies. The geometric branch~\eqref{eq:branch-geo} uses only the physical parts
$\cell\cap\Omega$ and never the ghost penalty, so it holds at every weight; its size is set by the \emph{shapes} of the
physical parts and degenerates for slivers. The weight branch~\eqref{eq:branch-wt} pays with the ghost energy; it is
bounded by $C(p)\gamma^{-1/2}$ uniformly in the cut position, $h$ and the mesh level, and on a fixed configuration it
saturates at $C_{\mathrm{tr}}C_E^{1/2}(1+R)^{1/2}$ once $\gamma\le(1+R)^{-1}$. At moderate weights the weight branch
gives a cut-uniform constant; as $\gamma\to0$ the geometric branch caps it.
The bound is pessimistic: in our experiments we do not observe a dependence on the cut fraction.

\begin{proof}
    The volume, boundary-penalty, and ghost terms are bounded by
    Cauchy--Schwarz directly in $\|\cdot\|_\gamma$. The first consistency term of~\eqref{eq:aform} is bounded by
    Cauchy--Schwarz and split into two factors,
    \begin{gather*}
        |(\partial_nv_h,w_h)_\Gamma|
        \le
        h^{1/2}\|\partial_nv_h\|_{0,\Gamma}
        \cdot
        h^{-1/2}\|w_h\|_{0,\Gamma};
    \end{gather*}
    the second consistency term is treated identically with the roles of
    $v_h$ and $w_h$ interchanged.
    By~\eqref{eq:tau}, the second factor satisfies
    \begin{gather*}
        h^{-1/2}\|w_h\|_{0,\Gamma}
        \le
        (\tau_0p(p+1))^{-1/2}\|w_h\|_\gamma,
    \end{gather*}
    and the factor $(\tau_0p(p+1))^{-1/2}\le\tau_0^{-1/2}$ is absorbed into the
    constant $C$.
    There are two ways to bound the first factor. \emph{Geometric branch.} Since $|\partial_nv_h|\le|\nabla v_h|$ on
    $\Gamma$ and each component $\partial_iv_h$ belongs to $\Q_p(\cell)$, the definition of $C_{\mathrm{Nik}}$, applied
    componentwise on each cut cell and summed in its squared form over the (disjoint) cut cells, gives
    \begin{align*}
        h^{1/2}\|\partial_nv_h\|_{0,\Gamma}
         & \le
        C_{\mathrm{Nik}}(p,\mathcal G)
        \|\nabla v_h\|_{0,\Omega}
        \\
         & \le
        M_{\mathrm{geo}}\|v_h\|_\gamma.
    \end{align*}
    This bound does not use the ghost term, so it does not depend on the weight.
    \emph{Weight branch.} The cut-cell trace inequality~\eqref{eq:cuttrace} in its summed componentwise form, followed
    by~\eqref{eq:extension} at unit weight and~\eqref{eq:twobranch}, gives
    \begin{align*}
        h^{1/2}\|\partial_nv_h\|_{0,\Gamma}
         & \le
        C_{\mathrm{tr}}(p)\|\nabla v_h\|_{0,\Omega_h}
        \\
         & \le
        C_{\mathrm{tr}}(p)C_E(p)^{1/2}
        \Xi_\gamma^{1/2}\|v_h\|_\gamma
        =M_{\mathrm{wt}}(\gamma)\|v_h\|_\gamma,
    \end{align*}
    and $\Xi_\gamma\le\gamma^{-1}$ gives the bound in~\eqref{eq:branch-wt}.
\end{proof}

\begin{remark}[the measured constant follows the two branches]\label{rem:measured-Ca}
    In the two-cell slab model, one full cell and one cell cut at fraction $\kappa$, with $h=1$ and the boundary
    penalty~\eqref{eq:tau} at $\tau_0=5$, the Nikolskii constant is $C_{\mathrm{Nik}}(p,\mathcal G)=(p+1)\kappa^{-1/2}$
    exactly. The derivative of a discrete function has degree $p-1$, and its trace constant is $p\,\kappa^{-1/2}$. The
    symmetric part of $a_\gamma$ is $s_\gamma$, so $C_a=(1+\|T\|^2)^{1/2}$, where $T$ is the skew pair of consistency terms
    in the graph norm, and with the trace constant of the derivative in place of $C_{\mathrm{Nik}}$ the geometric
    branch~\eqref{eq:branch-geo} takes the explicit form
    \begin{gather*}
        C_a\le\bigl(1+p^2\kappa^{-1}(h\tau_h)^{-1}\bigr)^{1/2},
    \end{gather*}
    a bound free of $\gamma$.

    We computed \footnote{The matrices are assembled in exact rational arithmetic and the singular values are computed in
        floating-point arithmetic with 140 significant decimal digits. Double precision does not suffice: the Gram matrix of
        the cut cell has entries of size $\kappa^{2p}$, down to $2^{-140}$ in the sweep, and the swept weights reach
        $10^{-12}$.} $C_a$ for $p\le7$, $\kappa\in\{2^{-6},\dots,2^{-10}\}$ and $\gamma\in\{1,10^{-2},\dots,10^{-12}\}$;
    Table~\ref{tab:measured-Ca} lists the smallest cut fraction, and the others behave alike with smaller bounds. The
    measured constant never exceeds the bound. At $\gamma=1$ it is close to one for every degree and cut fraction, so the
    weight branch is cut-uniform. As the weight decreases, $C_a$ grows, far more slowly than $\gamma^{-1/2}$, and
    approaches the bound once the weight falls below the saturation threshold $(1+R)^{-1}$ of~\eqref{eq:shedfactor}. From
    $p=3$ on the threshold lies below the swept weights, consistent with Remark~\ref{rem:amplification}, and $C_a$ stays
    below the bound. These values support the minimum in~\eqref{eq:continuity} over the measured range.

    The identification $C_{\mathrm{Nik}}\simeq(p+1)\kappa_{\min}^{-1/2}$ is the slab-model value and does not survive in
    general: for an arbitrary cut cell the constant is governed by the shape of $\cell\cap\Omega$ and not by its volume
    fraction alone, and degenerates for slivers. This is why~\eqref{eq:continuity} is stated with
    $C_{\mathrm{Nik}}(p,\mathcal G)$ and not with a power of $\kappa_{\min}$. Only the conditioning at very small weights
    inherits this dependence; existence and uniqueness rest on Theorem~\ref{thm:coercivity} alone, which involves no
    Nikolskii constant, and the weight branch~\eqref{eq:branch-wt} is bounded by $C(p)\gamma^{-1/2}$, a genuine cut-uniform
    bound.
\end{remark}

\begin{table}[htbp]
    \centering
    \caption{Measured continuity constant $C_a$ in the two-cell slab model at cut fraction $\kappa=2^{-10}$, with the
        $\gamma$-free bound of Remark~\ref{rem:measured-Ca} and the saturation threshold $(1+R)^{-1}$
        of~\eqref{eq:shedfactor}.}
    \label{tab:measured-Ca}
    \small
    \begin{tabular}{lccccccccc}
        \toprule
               &       &                    & \multicolumn{7}{c}{$\gamma$}                                                                           \\
        \cmidrule(l){4-10}
               & bound & $(1+R)^{-1}$       & $1$                          & $10^{-2}$ & $10^{-4}$ & $10^{-6}$ & $10^{-8}$ & $10^{-10}$ & $10^{-12}$ \\
        \midrule
        $\Q_1$ & 10.2  & $9.7\cdot10^{-4}$  & $1.10$                       & $3.19$    & $9.69$    & $10.2$    & $10.2$    & $10.2$     & $10.2$     \\
        $\Q_2$ & 11.7  & $3.1\cdot10^{-10}$ & $1.08$                       & $2.04$    & $5.65$    & $5.93$    & $6.18$    & $10.6$     & $11.7$     \\
        $\Q_3$ & 12.4  & $4.4\cdot10^{-17}$ & $1.08$                       & $1.63$    & $4.06$    & $4.25$    & $4.43$    & $7.54$     & $8.31$     \\
        $\Q_4$ & 12.8  & $4.8\cdot10^{-24}$ & $1.08$                       & $1.43$    & $3.21$    & $3.35$    & $3.49$    & $5.87$     & $6.47$     \\
        $\Q_5$ & 13.1  & $4.6\cdot10^{-31}$ & $1.08$                       & $1.32$    & $2.68$    & $2.80$    & $2.91$    & $4.83$     & $5.31$     \\
        $\Q_6$ & 13.3  & $4.0\cdot10^{-38}$ & $1.08$                       & $1.26$    & $2.33$    & $2.42$    & $2.51$    & $4.12$     & $4.52$     \\
        $\Q_7$ & 13.4  & $3.2\cdot10^{-45}$ & $1.09$                       & $1.21$    & $2.08$    & $2.16$    & $2.23$    & $3.60$     & $3.95$     \\
        \bottomrule
    \end{tabular}
\end{table}

\begin{corollary}[well-posedness for every schedule]\label{cor:wellposed}
    Let $\tau_0>0$ and $\gamma\in(0,1]$. The discrete problem
    $a_\gamma(u_h,v_h)=\ell(v_h)$ has a unique solution, and
    $A_\gamma:V_h^p\to(V_h^p)'$ satisfies, in the graph norm,
    $\|A_\gamma\|\le C_a$, $\|A_\gamma^{-1}\|\le1$,
    $\operatorname{cond}(A_\gamma)\le C_a$.
\end{corollary}

\begin{proof}
    Immediate from Theorems~\ref{thm:coercivity} and~\ref{thm:continuity},
    $A_\gamma$ being square with trivial kernel.
\end{proof}

\begin{proposition}[active-mesh gradient control]\label{prop:activemesh}
    Let $\tau_0>0$ and $0<\gamma\le1$. Then
    \begin{gather}\label{eq:activemesh}
        \|\nabla v_h\|_{0,\Omega_h}^2
        \le
        C_E(p)\,\Xi_\gamma\,\|v_h\|_\gamma^2
        \le
        C_E(p)\,\gamma^{-1}\|v_h\|_\gamma^2,
        \qquad
        v_h\in V_h^p,
    \end{gather}
    with a constant independent of the cut and $h$.
\end{proposition}

\begin{proof}
    Apply~\eqref{eq:extension} at unit weight,
    $\|\nabla v_h\|^2_{0,\Omega_h}\le
        C_E(p)(\|\nabla v_h\|_{0,\Omega}^2+g_h(v_h,v_h))$, and
    then~\eqref{eq:twobranch}.
\end{proof}

\begin{corollary}[no floor, with an extension price]\label{cor:nofloor}
    Every schedule $\gamma_p\in(0,1]$, including
    $\gamma_p=\min\{1,\,Cp^qe^{-cp}\}$, retains cut- and level-uniform
    control of the active-mesh gradient, at the explicit price
    $C_E(p)\Xi_{\gamma_p}\le C_E(p)\gamma_p^{-1}$. On a fixed mesh and cut
    configuration and for fixed $p$, the price saturates as the weight
    decreases: writing $R_p$ for the amplification
    quotient~\eqref{eq:shedfactor} of $V_h^p$, once $\gamma_p\le R_p^{-1}$
    the factor $\Xi_{\gamma_p}$ lies within a factor of two of its limit
    value $1+R_p$.
\end{corollary}

\begin{proof}
    The first claim is~\eqref{eq:activemesh} read along the schedule: for
    every $p$ the constant $C_E(p)\Xi_{\gamma_p}\le C_E(p)\gamma_p^{-1}$ is
    finite and independent of the cut and of $h$, whatever the decay of
    $\gamma_p$; no uniformity in $p$ is claimed. For the saturation claim
    fix the mesh, the cut and the degree, so that $R_p<\infty$. If $R_p=0$
    the ghost form vanishes identically, $\Xi_\gamma=1$ for every
    $0<\gamma\le1$, and there is nothing to prove. If $R_p>0$ and
    $\gamma_p\le R_p^{-1}$, then
    $\gamma_p^{-1}\ge\max\{1,R_p\}\ge(1+R_p)/2$, hence
    \begin{gather*}
        \tfrac12(1+R_p)\le\Xi_{\gamma_p}\le1+R_p ,
    \end{gather*}
    and the price in~\eqref{eq:activemesh} stays within a factor of two of
    $C_E(p)(1+R_p)$.
\end{proof}

\begin{remark}[$\tau_0>0$ is not cosmetic]\label{rem:tau}
    At $\tau_h=0$ a constant finite element function $c$ has
    $a_\gamma(c,c)=0$, so the symmetric part of $A_\gamma$ is singular and
    $0$ lies in the numerical range of $B^{1/2}A_\gamma B^{1/2}$ for
    \emph{every} symmetric positive definite $B$. No congruence repairs
    this, so no choice of preconditioner restores a positive lower bound
    on the field of values required by the corresponding Krylov
    convergence bound. The
    penalty-free case~\cite{Burman12PenaltyFree} is not covered here; a slab
    computation gives an inf-sup constant at least $0.618$ over the tested
    cuts and weights, but this evidence is not a theorem. The
    \emph{magnitude} $\tau_0$ is free, since it never competes with
    $\gamma g_h$ in~\eqref{eq:coercivity-identity}, so it creates no lower
    bound on $\gamma$.
\end{remark}

\section{A priori error analysis}
\label{sec:error}
Well-posedness (Section~\ref{sec:wellposed}) shows that the discrete problem is well posed at every positive weight;
this section says what it converges \emph{to}, and how fast. The result is the classical pair for the nonsymmetric
Nitsche method, transplanted to the cut setting at every positive weight: energy rate $h^{p}$, $L^2$ rate $h^{p+1/2}$,
uniform in the mesh level and, at the price of a constant $\le C(p)\gamma^{-1/2}$ in energy, $\le C(p)\gamma^{-1}$ in
$L^2$, in the cut configuration. The weight enters through the same two-branch constant as Theorem~\ref{thm:continuity}
and nowhere else: cut-uniformity is bought at $\gamma^{-1/2}$ per flux estimate, weight-uniformity (down to arbitrarily
small $\gamma$) at the cut-dependent $C_{\mathrm{Nik}}(p,\mathcal G)$, and $\Xi_\gamma$ of~\eqref{eq:shedfactor}
saturates at $1+R$ on a fixed configuration. No simultaneous uniformity in $\gamma$ and the cut is claimed. The stated
convergence rates refer to refinement at fixed $p$ and fixed $\gamma>0$, including any degree-dependent choice
$\gamma_p$ held fixed during refinement. If the weight also depends on $h$, its contribution through $M_{\gamma(h)}$
must be included in the convergence rate. The fitted-mesh analysis of~\cite{ChenLiuZhang25} gives refined rates under
$W^{p+1,q}$ regularity and numerical evidence of a half-order loss for piecewise-linear elements near the $H^2$
endpoint. These results motivate the comparison in Remark~\ref{rem:import}; they do not prove sharpness for the present
unfitted method. Throughout, exact geometry and exact quadrature are assumed; the constants $C(p)$ are not claimed
$p$-uniform, consistent with the rest of the paper.

\subsection{Consistency and approximation}

Let the solution $u$ of the model problem~\eqref{eq:model} belong to $H^{p+1}(\Omega)$, $p\ge1$, and let $u_h$ solve
the discrete problem~\eqref{eq:discrete-problem}, uniquely solvable for every $\gamma>0$ by
Corollary~\ref{cor:wellposed}. Let $u^e:=Eu\in H^{p+1}(\R^d)$ be a Stein extension~\cite{Stein70}, with
$\|u^e\|_{H^{p+1}(\Omega_h)}\le C_S\|u\|_{H^{p+1}(\Omega)}$, $C_S$ depending on $\Omega$ only, in particular
independent of $h$, of the cut position, and of $\gamma$. Extend the ghost penalty to $V_h^p+H^{p+1}(\Omega_h)$: for
$w\in H^{p+1}(\Omega_h)$ and $k\le p$ the one-sided traces of $\partial_{n_F}^kw$ on an interior face exist and agree,
so $\llbracket\partial_{n_F}^kw^e\rrbracket=0$ and
\begin{gather}\label{eq:ghost-consistency}
    g_h(u^e,v)=0
    \qquad\text{for all }v\in V_h^p+H^{p+1}(\Omega_h).
\end{gather}

\begin{lemma}[exact Galerkin orthogonality]\label{lem:orthogonality}
    For every $\gamma>0$,
    \begin{gather}\label{eq:orthogonality}
        a_\gamma(u^e-u_h,v_h)=0,
        \qquad v_h\in V_h^p.
    \end{gather}
\end{lemma}

\begin{proof}
    Integration by parts with $-\Delta u=f$ and $u=g$ on $\Gamma$ gives, for
    every $v_h\in V_h^p$ (with $\partial_nv_h$ understood elementwise on
    $\Gamma$),
    $(\nabla u,\nabla v_h)_\Omega-(\partial_nu,v_h)_\Gamma
        +(u,\partial_nv_h)_\Gamma+\tau_h(u,v_h)_\Gamma=\ell(v_h)$; the volume and
    boundary terms of $a_\gamma(u^e,v_h)$ see only $u^e|_\Omega=u$, and the
    ghost term vanishes by~\eqref{eq:ghost-consistency}.
    Subtract~\eqref{eq:discrete-problem}.
\end{proof}

The ghost penalty is thus not merely weakly consistent here: on the extension it is \emph{exactly} consistent, and its
entire contribution to the error analysis is the term $\gamma\,g_h(I_hu^e,\cdot)$ of the interpolant below, which
carries a factor $\gamma^{1/2}$ and therefore \emph{shrinks} with the weight. The $k$-th normal derivative has a
single-valued $L^2$ trace on a ghost face whenever $u^e\in H^s(\Omega_h)$ with $s>k+1/2$. Thus exact ghost consistency
holds already for $s>p+1/2$, below the full approximation regularity $p+1$. For lower regularity, some ghost terms need
not be defined on $u^e$.

\begin{remark}[consistency below full regularity]\label{rem:lowreg}
    Let $u\in H^s(\Omega)$ with $3/2<s<p+1$, and take a bounded extension $u^e=Eu\in H^s(\R^d)$. For every integer $1\le
        k\le p$ with $s>k+1/2$, the traces of $\partial^k_{n_F}u^e$ exist and their jumps vanish. If $s>p+1/2$, all ghost terms
    and exact Galerkin orthogonality remain available. For $s\le p+1/2$, the highest-order traces are not guaranteed, so
    the graph norm of $u^e-u_h$ need not be defined.

    An argument valid throughout $3/2<s<p+1$ uses only the physical part of the form: writing
    $a_\gamma=a^{\mathrm{ph}}+\gamma g_h$ (the four physical terms of~\eqref{eq:aform}), integration by parts as in
    Lemma~\ref{lem:orthogonality} gives $a^{\mathrm{ph}}(u^e,v_h)=\ell(v_h)$ for $u\in H^s$, $s>3/2$, hence for every
    $v_h\in V_h^p$
    \begin{gather}\label{eq:orthogonality-lowreg}
        a_\gamma(I_hu^e-u_h,v_h)
        =a^{\mathrm{ph}}(I_hu^e-u^e,v_h)+\gamma\,g_h(I_hu^e,v_h),
    \end{gather}
    both terms well defined. Here $I_h$ is the interpolant of
    Lemma~\ref{lem:interp}, and $\eta:=u^e-I_hu^e$.

    The volume and boundary estimates of~\eqref{eq:interp} hold with $h^p$ replaced by $h^{s-1}$ and
    $\|u\|_{H^{p+1}(\Omega)}$ replaced by $\|u\|_{H^s(\Omega)}$, with constants also depending on $s$. For $3/2<s<2$,
    choose $1/2<t<s-1$ and apply the scaled fractional trace inequality on each full cell,
    \begin{gather*}
        \|w\|_{0,\Gamma\cap\cell}^2
        \le C(t)\bigl(h^{-1}\|w\|_{0,\cell}^2
        +h^{2t-1}|w|_{H^t(\cell)}^2\bigr),
    \end{gather*}
    to each component of $\nabla\eta$. The interpolation estimates in
    $H^1$ and $H^{1+t}$ give
    \begin{gather*}
        h^{1/2}\|\partial_n\eta\|_{0,\Gamma}
        \le C(p,s)h^{s-1}\|u\|_{H^s(\Omega)}.
    \end{gather*}
    For $s\ge2$, the integer-order trace argument of
    Lemma~\ref{lem:interp} applies.

    To bound the ghost term without taking high-order traces of $u^e$, let $B_F$ be a cube of diameter $Ch$ containing both
    cells adjacent to $F$ and their interpolation patches. These cubes can be chosen with bounded overlap. With $m:=\lceil
        s\rceil-1\le p$, choose a polynomial $P_F$ of total degree at most $m$ such that
    \begin{gather*}
        \|u^e-P_F\|_{0,B_F}+h|u^e-P_F|_{1,B_F}
        \le C(p,s)h^s|u^e|_{H^s(B_F)}.
    \end{gather*}
    Polynomial reproduction gives $I_hP_F=P_F$, whose normal-derivative
    jumps vanish. Polynomial inverse trace estimates, followed by the local
    stability bound $\|I_hw\|_{0,\cell}\le C(p)\bigl(\|w\|_{0,\omega_\cell}+h|w|_{1,\omega_\cell}\bigr)$ of the
    Scott--Zhang interpolant (the $H^1$ term enters because the interpolant uses face averages), therefore give, for
    $1\le k\le p$,
    \begin{align*}
        h_F^{k-1/2}
        \|\llbracket\partial_{n_F}^kI_hu^e\rrbracket\|_{0,F}
         & \le C(p)h^{-1}
        \sum_{\cell\supset F}\|I_h(u^e-P_F)\|_{0,\cell}    \\
         & \le C(p)h^{-1}
        \bigl(\|u^e-P_F\|_{0,B_F}+h|u^e-P_F|_{1,B_F}\bigr) \\
         & \le C(p,s)h^{s-1}|u^e|_{H^s(B_F)}.
    \end{align*}
    Squaring and summing over faces and derivative orders, then using
    the extension bound, yields
    \begin{gather}\label{eq:lowreg-ghost}
        g_h(I_hu^e,I_hu^e)^{1/2}
        \le C(p,s)h^{s-1}\|u\|_{H^s(\Omega)}.
    \end{gather}

    Set $\xi:=I_hu^e-u_h$. Testing~\eqref{eq:orthogonality-lowreg} with $\xi$ and estimating the four physical terms as in
    Theorem~\ref{thm:energy} gives
    \begin{gather*}
        \|\xi\|_\gamma
        \le C(p,s,\tau_0)(1+M_\gamma)h^{s-1}\|u\|_{H^s(\Omega)}.
    \end{gather*}
    The ghost residual is bounded by
    $\gamma^{1/2}g_h(I_hu^e,I_hu^e)^{1/2}\|\xi\|_\gamma$.
    Adding the physical interpolation error gives
    \begin{gather}\label{eq:lowreg-energy}
        \bigl(\|\nabla(u-u_h)\|_{0,\Omega}^2
        +\tau_h\|u-u_h\|_{0,\Gamma}^2\bigr)^{1/2}
        \le C(p,s,\tau_0)(1+M_\gamma)h^{s-1}\|u\|_{H^s(\Omega)}.
    \end{gather}
    Only when $s>p+1/2$ do we additionally assert this bound for the full
    graph error $\|u^e-u_h\|_\gamma$.

    Under the dual regularity assumption of Theorem~\ref{thm:ltwo}, its physical duality argument gives
    \begin{gather}\label{eq:lowreg-ltwo}
        \|u-u_h\|_{0,\Omega}
        \le C(p,s,\tau_0)\bigl((1+M_\gamma)h^{s-1/2}
        +(1+M_\gamma)^2h^s\bigr)\|u\|_{H^s(\Omega)}.
    \end{gather}
    Indeed, $a^{\mathrm{ph}}(u-u_h,v_h)=\gamma g_h(u_h,v_h)$ replaces
    full orthogonality, so the ghost contribution to the dual identity is
    $+\gamma g_h(u_h,I_h\Psi^e)$. Since
    $\gamma^{1/2}g_h(u_h,u_h)^{1/2}
        \le\gamma^{1/2}g_h(I_hu^e,I_hu^e)^{1/2}+\|\xi\|_\gamma$,
    this contribution is bounded by
    $C(p,s,\tau_0)(1+M_\gamma)h^s\|u\|_{H^s(\Omega)}\|\Psi\|_{H^2(\Omega)}$.
    The explicit ghost residual in the energy argument has upper bound
    $C(p,s)\gamma^{1/2}h^{s-1}\|u\|_{H^s(\Omega)}$; this contribution
    decreases with the weight, without implying monotonicity of the error.
\end{remark}

\begin{lemma}[interpolation on the active mesh]\label{lem:interp}
    Let $u\in H^{p+1}(\Omega)$ and let $I_h:H^1(\Omega_h)\to V_h^p$ be the
    Scott--Zhang interpolant of the active mesh, built on whole cells, including
    the parts of cut cells outside $\Omega$, so that no cut geometry enters its
    construction, and write $\eta:=u^e-I_hu^e$. Then,
    with $C(p)$ independent of $h$, the cut, and $\gamma$,
    \begin{gather}\label{eq:interp}
        \|\nabla\eta\|_{0,\Omega_h}
        +h^{1/2}\|\partial_n\eta\|_{0,\Gamma}
        +\tau_h^{1/2}\|\eta\|_{0,\Gamma}
        \ \le\
        C(p)\,h^{p}\,\|u\|_{H^{p+1}(\Omega)},
        \\
        \label{eq:interp-ghost}
        g_h(I_hu^e,I_hu^e)^{1/2}
        =g_h(\eta,\eta)^{1/2}
        \ \le\
        C(p)\,h^{p}\,\|u\|_{H^{p+1}(\Omega)}.
    \end{gather}
\end{lemma}

\begin{proof}
    Since $I_h$ is locally stable and polynomial-preserving, and is built on the active mesh alone, its interpolation
    patches $\omega_\cell$ consist of full Cartesian cells, also for cells touching $\partial\Omega_h$. These patches are
    shape regular and have overlap bounded independently of $h$ and of the cut. The elementwise estimates
    $|\eta|_{H^m(\cell)}\le C(p)\,h^{p+1-m}|u^e|_{H^{p+1}(\omega_\cell)}$, $0\le m\le p+1$, hold with cut-independent
    constants. Summation over the active mesh gives the volume estimate $\|\nabla\eta\|_{0,\Omega_h}\le
        C(p)h^p|u^e|_{H^{p+1}(\Omega_h)}$; the Stein bound converts $|u^e|_{H^{p+1}(\Omega_h)}$ to $\|u\|_{H^{p+1}(\Omega)}$.

    For the boundary terms use the trace inequality for $H^1$ functions on the \emph{full} cell
    $\cell\supset\Gamma\cap\cell$, $\|w\|^2_{0,\Gamma\cap\cell}\le C(h^{-1}\|w\|^2_{0,\cell}+h\|\nabla w\|^2_{0,\cell})$,
    whose constant is cut-uniform for $\Gamma$ resolved in the sense of the geometric assumption~\cite{BurmanHansbo12}. It
    is applied cell by cell, since $\nabla I_hu^e$ is discontinuous across interior faces: with $w=\eta$ and with $w$ equal
    to each component of $\nabla\eta$ (legitimate since $p\ge1$ gives $u^e\in H^2$), the elementwise estimates give
    $\|\eta\|^2_{0,\Gamma\cap\cell} \le C(p)\,h^{2p+1}|u^e|^2_{H^{p+1}(\omega_\cell)}$ and
    $\|\partial_n\eta\|^2_{0,\Gamma\cap\cell} \le C(p)\,h^{2p-1}|u^e|^2_{H^{p+1}(\omega_\cell)}$; multiplication by
    $\tau_h=\tau_0p(p+1)/h$ and by $h$, summation over the cut cells with bounded patch overlap, and the Stein bound yield
    $\tau_h\|\eta\|^2_{0,\Gamma}\le C(p)\,h^{2p}\,\|u\|^2_{H^{p+1}(\Omega)}$ and $h\|\partial_n\eta\|^2_{0,\Gamma}\le
        C(p)\,h^{2p}\,\|u\|^2_{H^{p+1}(\Omega)}$.

    For~\eqref{eq:interp-ghost}, the jumps of $\partial^k_{n_F}u^e$ vanish on every ghost face for $k\le p$, the
    observation behind~\eqref{eq:ghost-consistency}, so $\llbracket\partial^k_{n_F}I_hu^e\rrbracket
        =-\llbracket\partial^k_{n_F}\eta\rrbracket$; the sign cancels in the quadratic form and
    $g_h(I_hu^e,I_hu^e)=g_h(\eta,\eta)$. In the $k$-th face term, the facewise trace inequality on either (full) cell
    $\cell$ sharing $F$ and the elementwise estimates give $\|\partial^k_{n_F}\eta\|^2_{0,F} \le
        C(p)\bigl(h^{-1}|\eta|^2_{H^k(\cell)} +h|\eta|^2_{H^{k+1}(\cell)}\bigr) \le
        C(p)\,h^{2p+1-2k}|u^e|^2_{H^{p+1}(\omega_\cell)}$, the case $k=p$ being covered by the estimate with $m=p+1$, since
    $\eta|_\cell\in H^{p+1}(\cell)$. The scaling $h_F^{2k-1}$ absorbs the inverse powers exactly as in~\eqref{eq:converse}
    and leaves $h^{2p}$; summation over the ghost faces with bounded patch overlap and the Stein bound give the
    weak-consistency estimate~\cite{BurmanHansbo12,MassingLarsonLogg14} $g_h(\eta,\eta)\le
        C(p)\,h^{2p}\,\|u\|^2_{H^{p+1}(\Omega)}$.
\end{proof}

\subsection{Energy-norm error estimate}

All weight dependence below rides on one quantity, the flux-control constant of Theorem~\ref{thm:continuity}:
\begin{gather}\label{eq:Mgamma}
    M_\gamma
    :=\min\bigl\{M_{\mathrm{geo}},\,M_{\mathrm{wt}}(\gamma)\bigr\},
    \qquad
    h^{1/2}\|\partial_n\xi\|_{0,\Gamma}\le M_\gamma\|\xi\|_\gamma
    \ \ \forall\xi\in V_h^p,
\end{gather}
with the branches~\eqref{eq:branch-geo} and~\eqref{eq:branch-wt}: the weight branch is cut-uniform at the price
$\gamma^{-1/2}$, the geometric branch is weight-uniform at the price of the cut geometry, and on
a fixed configuration $M_\gamma$ saturates once $\gamma\le(1+R)^{-1}$.

\begin{theorem}[energy estimate at every positive weight]\label{thm:energy}
    Let $\tau_0>0$, $\gamma\in(0,1]$, and $u\in H^{p+1}(\Omega)$. Then
    \begin{gather}\label{eq:energy-rate}
        \|u^e-u_h\|_\gamma
        \ \le\
        C(p,\tau_0)\,(1+M_\gamma)\,h^{p}\,\|u\|_{H^{p+1}(\Omega)},
    \end{gather}
    with $C(p,\tau_0)$ independent of $h$, the mesh level, the cut
    configuration and $\gamma$.

    In particular the physical errors obey
    \begin{gather}\label{eq:physical-rates}
        \|\nabla(u-u_h)\|_{0,\Omega}\ \le\
        C(p,\tau_0)(1+M_\gamma)\,h^{p}\,\|u\|_{H^{p+1}(\Omega)},
        \qquad
        \|u-u_h\|_{0,\Gamma}\ \le\
        C(p,\tau_0)(1+M_\gamma)\,h^{p+1/2}\,\|u\|_{H^{p+1}(\Omega)},
    \end{gather}
    cut-uniformly with $1+M_\gamma\le C(p)\gamma^{-1/2}$, and
    weight-uniformly, down to arbitrarily small $\gamma$, with
    $1+M_\gamma\le C(1+C_{\mathrm{Nik}}(p,\mathcal G))$.
\end{theorem}

\begin{proof}
    Set $\xi:=I_hu^e-u_h\in V_h^p$. The coercivity identity~\eqref{eq:coercivity-identity} and Galerkin
    orthogonality~\eqref{eq:orthogonality} give
    \begin{gather*}
        \|\xi\|_\gamma^2
        =a_\gamma(\xi,\xi)
        =a_\gamma(I_hu^e-u^e,\xi)
        =-a_\gamma(\eta,\xi).
    \end{gather*}
    The volume term satisfies
    \begin{gather*}
        |(\nabla\eta,\nabla\xi)_\Omega|
        \le
        \|\nabla\eta\|_{0,\Omega}\,\|\xi\|_\gamma.
    \end{gather*}
    The boundary part of the graph norm and~\eqref{eq:tau} yield
    \begin{align*}
        |(\partial_n\eta,\xi)_\Gamma|
         & \le
        \|\partial_n\eta\|_{0,\Gamma}\,\|\xi\|_{0,\Gamma}
        \\
         & \le
        (\tau_0p(p+1))^{-1/2}
        h^{1/2}\|\partial_n\eta\|_{0,\Gamma}\,\|\xi\|_\gamma.
    \end{align*}
    The skew term is the only contribution that tests a \emph{discrete} normal derivative. The flux
    estimate~\eqref{eq:Mgamma} gives
    \begin{align*}
        |(\eta,\partial_n\xi)_\Gamma|
         & \le
        \bigl(h^{-1/2}\|\eta\|_{0,\Gamma}\bigr)
        \bigl(h^{1/2}\|\partial_n\xi\|_{0,\Gamma}\bigr)
        \\
         & \le
        (\tau_0p(p+1))^{-1/2}M_\gamma
        \tau_h^{1/2}\|\eta\|_{0,\Gamma}\,\|\xi\|_\gamma.
    \end{align*}
    The penalty term satisfies
    \begin{gather*}
        \tau_h|(\eta,\xi)_\Gamma|
        \le
        \tau_h^{1/2}\|\eta\|_{0,\Gamma}\,\|\xi\|_\gamma.
    \end{gather*}
    Finally, the Cauchy--Schwarz inequality for the ghost penalty implies
    \begin{align*}
        \gamma|g_h(\eta,\xi)|
         & \le
        \gamma^{1/2}g_h(\eta,\eta)^{1/2}
        \bigl(\gamma g_h(\xi,\xi)\bigr)^{1/2}
        \\
         & \le
        \gamma^{1/2}g_h(\eta,\eta)^{1/2}\,\|\xi\|_\gamma.
    \end{align*}
    If $\xi=0$, the desired estimate is immediate. Otherwise, collecting these bounds and dividing by $\|\xi\|_\gamma$
    gives
    \begin{align*}
        \|\xi\|_\gamma
        \le{} &
        \|\nabla\eta\|_{0,\Omega}
        +(\tau_0p(p+1))^{-1/2}
        h^{1/2}\|\partial_n\eta\|_{0,\Gamma}
        \\
              & +\bigl(1+(\tau_0p(p+1))^{-1/2}M_\gamma\bigr)
        \tau_h^{1/2}\|\eta\|_{0,\Gamma}
        +\gamma^{1/2}g_h(\eta,\eta)^{1/2}.
    \end{align*}
    Therefore, the interpolation estimates
    \eqref{eq:interp}--\eqref{eq:interp-ghost} imply
    \begin{gather}\label{eq:xi-bound}
        \|\xi\|_\gamma\ \le\ C(p,\tau_0)\,(1+M_\gamma)\,h^{p}\,
        \|u\|_{H^{p+1}(\Omega)}.
    \end{gather}
    Since $\gamma\le1$, the same interpolation estimates give
    \begin{align*}
        \|\eta\|_\gamma^2
         & =
        \|\nabla\eta\|^2_{0,\Omega}
        +\tau_h\|\eta\|^2_{0,\Gamma}
        +\gamma g_h(\eta,\eta)
        \\
         & \le
        C(p,\tau_0)h^{2p}\|u\|^2_{H^{p+1}(\Omega)}.
    \end{align*}
    The decomposition $u^e-u_h=\eta+\xi$, the triangle inequality, and~\eqref{eq:xi-bound} prove~\eqref{eq:energy-rate}.
    The two components of the graph norm then give~\eqref{eq:physical-rates}; for the boundary error,
    \begin{gather*}
        \|u-u_h\|_{0,\Gamma}
        \le
        \tau_h^{-1/2}\|u^e-u_h\|_\gamma
        =
        (\tau_0p(p+1))^{-1/2}h^{1/2}\|u^e-u_h\|_\gamma.
    \end{gather*}
\end{proof}

Exactly one contribution in the proof, the skew Nitsche term that pairs a continuous approximation error with the
discrete flux $\partial_n\xi$ on $\Gamma$, sees the weight, through the same mechanism and the same constant as
Theorem~\ref{thm:continuity}. The ghost-penalty consistency contribution decreases with the weight. The weight never
spreads: it enters the argument once, through the flux constant; no lower bound proving sharpness of the resulting
error constant is asserted.

\subsection{$L^2$-error estimate and the half order}

\begin{theorem}[$L^2$ rate: the nonsymmetric half order]\label{thm:ltwo}
    Assume in addition that $\Omega$ is convex or $C^{1,1}$, so that the
    Dirichlet dual problem has $H^2$ regularity. Then, with $M_\gamma$ as
    in~\eqref{eq:Mgamma},
    \begin{gather}\label{eq:ltwo-rate}
        \|u-u_h\|_{0,\Omega}
        \ \le\
        C(p,\tau_0)\Bigl((1+M_\gamma)\,h^{p+1/2}+(1+M_\gamma)^2\,h^{p+1}\Bigr)
        \|u\|_{H^{p+1}(\Omega)},
    \end{gather}
    uniformly in the mesh level and the cut configuration; cut-uniformly
    the bracket is
    $\le C(p,\tau_0)(\gamma^{-1/2}h^{p+1/2}+\gamma^{-1}h^{p+1})$,
    and weight-uniformly (down to arbitrarily small $\gamma$) it is
    $\le C(1+C_{\mathrm{Nik}})h^{p+1/2}+C(1+C_{\mathrm{Nik}})^2h^{p+1}$.
\end{theorem}

\begin{proof}
    Write $e:=u-u_h$ on $\Omega$ and $\tilde e:=u^e-u_h$ on $\Omega_h$, and
    let $\Psi\in H^2(\Omega)\cap H^1_0(\Omega)$ solve $-\Delta\Psi=e$ with
    $\|\Psi\|_{H^2(\Omega)}\le C\|e\|_{0,\Omega}$. Let $\Psi^e$ be its Stein
    extension and $\zeta:=\Psi^e-I_h\Psi^e$. Lemma~\ref{lem:interp} is stated
    at full regularity, but its proof uses the elementwise Scott--Zhang
    estimates only up to the regularity that is available; at $H^2$ these
    read $|\zeta|_{H^m(\cell)}\le C(p)\,h^{2-m}\,|\Psi^e|_{H^2(\omega_\cell)}$
    for $0\le m\le2$, and combined with the same cellwise trace inequalities
    on full cells they give
    \begin{gather}\label{eq:dual-interp}
        \|\nabla\zeta\|_{0,\Omega}+h^{-1/2}\|\zeta\|_{0,\Gamma}
        +h^{1/2}\|\partial_n\zeta\|_{0,\Gamma}
        \le C(p)\,h\,\|\Psi\|_{H^2(\Omega)},
        \qquad
        g_h(I_h\Psi^e,I_h\Psi^e)^{1/2}\le C(p)\,h\,\|\Psi\|_{H^2(\Omega)},
    \end{gather}
    the ghost bound needing a separate argument at $H^2$ regularity, since
    for $k\ge2$ the jump
    $\llbracket\partial_{n_F}^kI_h\Psi^e\rrbracket$ cannot be compared with
    $\partial^k\Psi^e$, which has no trace. The $k=1$ term is the usual one:
    $\llbracket\partial_{n_F}I_h\Psi^e\rrbracket
        =\llbracket\partial_{n_F}(I_h\Psi^e-\Psi^e)\rrbracket$ since
    $\nabla\Psi^e\in H^1$ has a single-valued trace, and the facewise trace
    inequality with elementwise interpolation gives
    $h_F^{1/2}\|\llbracket\partial_{n_F}I_h\Psi^e\rrbracket\|_{0,F}\le
        Ch\,|\Psi^e|_{H^2(\omega_F)}$. For $k\ge2$ drop the jump structure
    entirely and bound each one-sided contribution on its cell $\cell$: the
    trace--inverse inequality for polynomials and then the inverse estimate
    down to order two give
    \begin{gather*}
        \|\partial_{n_F}^kI_h\Psi^e\|_{0,F}
        \ \le\ Ch^{-1/2}\|D^kI_h\Psi^e\|_{0,\cell}
        \ \le\ Ch^{-1/2}\,h^{-(k-2)}\|D^2I_h\Psi^e\|_{0,\cell}
        \ \le\ Ch^{-(k-3/2)}\,|\Psi^e|_{H^2(\omega_\cell)},
    \end{gather*}
    the last step by $H^2$-stability of the Scott--Zhang interpolant,
    $\|D^2I_h\Psi^e\|_{0,\cell}\le\|D^2(I_h\Psi^e-\Psi^e)\|_{0,\cell}
        +\|D^2\Psi^e\|_{0,\cell}\le C|\Psi^e|_{H^2(\omega_\cell)}$. The face
    scaling then closes exactly:
    $h_F^{(2k-1)/2}\cdot h^{-(k-3/2)}=h$, uniformly in $k\le p$, and summing
    over the ghost faces with bounded overlap of the patches
    $\omega_\cell$ yields the stated bound.

    Integration by parts yields
    \begin{gather*}
        \|e\|^2_{0,\Omega}
        = (\nabla e,\nabla\Psi)_\Omega-(e,\partial_n\Psi)_\Gamma.
    \end{gather*}
    Testing orthogonality~\eqref{eq:orthogonality} with
    $I_h\Psi^e\in V_h^p$ and expanding~\eqref{eq:aform} gives
    \begin{gather*}
        0=a_\gamma(\tilde e,I_h\Psi^e)
        =(\nabla e,\nabla I_h\Psi^e)_\Omega
        +(\partial_ne,\zeta)_\Gamma
        +(e,\partial_nI_h\Psi^e)_\Gamma
        -\tau_h(e,\zeta)_\Gamma
        +\gamma\,g_h(\tilde e,I_h\Psi^e).
    \end{gather*}
    Here, $I_h\Psi^e|_\Gamma=-\zeta|_\Gamma$ because
    $\Psi|_\Gamma=0$. Substituting this identity and
    $\nabla\Psi=\nabla I_h\Psi^e+\nabla\zeta$ into the integration-by-parts
    formula gives
    \begin{gather*}
        \|e\|^2_{0,\Omega}
        =\underbrace{(\nabla e,\nabla\zeta)_\Omega}_{T_1}
        -\underbrace{(\partial_ne,\zeta)_\Gamma}_{T_2}
        +\underbrace{(e,\partial_n\zeta)_\Gamma}_{T_3'}
        -\underbrace{2(e,\partial_n\Psi)_\Gamma}_{T_3}
        +\underbrace{\tau_h(e,\zeta)_\Gamma}_{T_4}
        -\underbrace{\gamma\,g_h(\tilde e,I_h\Psi^e)}_{T_5},
    \end{gather*}
    where we used
    $\partial_nI_h\Psi^e=\partial_n\Psi^e-\partial_n\zeta$ and
    $\Psi^e|_\Omega=\Psi$.
    The energy estimate and~\eqref{eq:dual-interp} imply
    \begin{align*}
        |T_1|
         & \le \|\nabla e\|_{0,\Omega}\|\nabla\zeta\|_{0,\Omega} \\
         & \le C(1+M_\gamma)h^{p+1}
        \|u\|_{H^{p+1}(\Omega)}\|\Psi\|_{H^2(\Omega)}.
    \end{align*}
    To bound $T_2$, first use~\eqref{eq:interp}, \eqref{eq:Mgamma}, and~\eqref{eq:xi-bound} to obtain
    \begin{align*}
        \|\partial_ne\|_{0,\Gamma}
         & \le \|\partial_n\eta\|_{0,\Gamma}
        +h^{-1/2}M_\gamma\|\xi\|_\gamma      \\
         & \le C(1+M_\gamma)^2h^{p-1/2}
        \|u\|_{H^{p+1}(\Omega)}.
    \end{align*}
    Since~\eqref{eq:dual-interp} gives
    $\|\zeta\|_{0,\Gamma}\le Ch^{3/2}\|\Psi\|_{H^2(\Omega)}$, it follows
    that
    \begin{gather*}
        |T_2|
        \le C(1+M_\gamma)^2h^{p+1}
        \|u\|_{H^{p+1}(\Omega)}\|\Psi\|_{H^2(\Omega)}.
    \end{gather*}
    The boundary estimate~\eqref{eq:physical-rates}, the definition of $\tau_h$, and~\eqref{eq:dual-interp} yield
    \begin{align*}
        |T_3'|
         & \le \|e\|_{0,\Gamma}\|\partial_n\zeta\|_{0,\Gamma} \\
         & \le C(1+M_\gamma)h^{p+1}
        \|u\|_{H^{p+1}(\Omega)}\|\Psi\|_{H^2(\Omega)},        \\
        |T_4|
         & \le \tau_h\|e\|_{0,\Gamma}\|\zeta\|_{0,\Gamma}     \\
         & \le C(p,\tau_0)(1+M_\gamma)h^{p+1}
        \|u\|_{H^{p+1}(\Omega)}\|\Psi\|_{H^2(\Omega)}.
    \end{align*}
    By Cauchy--Schwarz,~\eqref{eq:dual-interp}, and Theorem~\ref{thm:energy},
    \begin{align*}
        |T_5|
         & \le \bigl(\gamma g_h(\tilde e,\tilde e)\bigr)^{1/2}
        \gamma^{1/2}g_h(I_h\Psi^e,I_h\Psi^e)^{1/2}             \\
         & \le \|\tilde e\|_\gamma\,\gamma^{1/2}Ch
        \|\Psi\|_{H^2(\Omega)}                                 \\
         & \le C(1+M_\gamma)h^{p+1}
        \|u\|_{H^{p+1}(\Omega)}\|\Psi\|_{H^2(\Omega)}.
    \end{align*}
    In contrast, the adjoint-inconsistency term satisfies only
    \begin{gather*}
        |T_3|=2\bigl|(e,\partial_n\Psi)_\Gamma\bigr|
        \le 2\|e\|_{0,\Gamma}\|\partial_n\Psi\|_{0,\Gamma}
        \le C(1+M_\gamma)h^{p+1/2}
        \|u\|_{H^{p+1}(\Omega)}\|\Psi\|_{H^2(\Omega)}.
    \end{gather*}
    This follows from~\eqref{eq:physical-rates} and the trace estimate
    $\|\partial_n\Psi\|_{0,\Gamma}\le C\|\Psi\|_{H^2(\Omega)}$.
    No additional factor of $h$ is available because $\Psi$ vanishes on
    $\Gamma$, whereas its normal derivative generally does not.
    Collecting the estimates and using dual regularity gives
    \begin{align*}
        \|e\|_{0,\Omega}^2
         & \le C(p,\tau_0)\Bigl((1+M_\gamma)h^{p+1/2}
        +(1+M_\gamma)^2h^{p+1}\Bigr)
        \|u\|_{H^{p+1}(\Omega)}\|\Psi\|_{H^2(\Omega)} \\
         & \le C(p,\tau_0)\Bigl((1+M_\gamma)h^{p+1/2}
        +(1+M_\gamma)^2h^{p+1}\Bigr)
        \|u\|_{H^{p+1}(\Omega)}\|e\|_{0,\Omega}.
    \end{align*}
    If $e=0$, then~\eqref{eq:ltwo-rate} is trivial; otherwise dividing by
    $\|e\|_{0,\Omega}$ proves~\eqref{eq:ltwo-rate}.
\end{proof}

\begin{remark}[the half order and the fitted-mesh scale]\label{rem:halforder}
    Every term of the proof is $O(h^{p+1})$ except $T_3$, produced by the sign of the Nitsche consistency terms. On fitted
    simplicial meshes over convex polytopes, Chen, Liu and Zhang~\cite{ChenLiuZhang25} prove, for $u\in W^{p+1,q}(\Omega)$
    and boundary penalty $h^{-\alpha}$ (the scaling~\eqref{eq:tau} is $\alpha=1$),
    \begin{gather}\label{eq:clz}
        \|u-u_h\|_{0,\Omega}\le C\,h^{r}\,|u|_{W^{p+1,q}(\Omega)},
        \qquad
        r=\min\bigl\{p+1,\ p+\max\{1,\alpha\}-1/q\bigr\},
    \end{gather}
    $q\in[2,\infty]$. At $\alpha=1$ this gives $h^{p+1/2}$ for $q=2$, the
    rate of Theorem~\ref{thm:ltwo}, and the optimal $h^{p+1}$ for
    $q=\infty$, at every degree.

    The two endpoints rest on different evidence. That smooth data restore the optimal order is a theorem, and it explains
    the optimal rates routinely measured on smooth manufactured solutions. That the half order is genuinely \emph{lost} at
    $q=2$ is not: no lower bound is proved, and the support is numerical, piecewise linears in two and three dimensions on
    data just above the $H^2$ endpoint. On cut meshes, Section~\ref{sec:num-halforder} observes the half order for $\Q_1$
    and $\Q_2$ on such data. Theorem~\ref{thm:ltwo} is accordingly an upper bound that is not known to be improvable, not a
    rate proved to be attained.
\end{remark}

\begin{remark}[scope of the import]\label{rem:import}
    Theorem~\ref{thm:ltwo} is proved here, in the cut setting, at every
    positive weight; \eqref{eq:clz} is quoted for fitted simplicial meshes
    on convex polytopes and says nothing directly about Cartesian cut
    meshes. Transplanting it is plausible, the mechanism being
    boundary-local, a boundary-strip Ritz projection with projected trace
    and a maximum-norm boundary-layer estimate, acting exactly where the
    ghost penalty lives and where Lemma~\ref{lem:interp} operates. The
    expected statement is~\eqref{eq:clz} with additional dependence on the
    cut and the weight. This is left open; nothing later depends on it.
\end{remark}

\begin{remark}[weight dependence of the error constants]\label{rem:weight-once}
    The weight enters the error
    analysis through exactly one mechanism, control of the discrete boundary
    flux $h^{1/2}\|\partial_n\xi\|_{0,\Gamma}$ by the graph norm, that is
    the constant $M_\gamma$ of~\eqref{eq:Mgamma}, which appears in the skew
    term $(\eta,\partial_n\xi)_\Gamma$ of the proof of Theorem~\ref{thm:energy}
    and in $T_2$ of the proof of Theorem~\ref{thm:ltwo}; and it \emph{helps}
    at exactly one: the ghost consistency terms $\gamma g_h(\eta,\xi)$ and
    $T_5$ carry $\gamma^{1/2}$ and shrink with the weight. Both
    branches of $M_\gamma$ transfer to the rates: cut-uniformly the
    constants grow like $\gamma^{-1/2}$ (energy) and
    $\gamma^{-1/2},\gamma^{-1}$ ($L^2$). The factor
    $\Xi_\gamma$ saturates at $1+R$ on a fixed configuration;
    weight-uniformly, down to arbitrarily small $\gamma$, the error constants are
    bounded by powers of $1+C_{\mathrm{Nik}}(p,\mathcal G)$. So the method
    converges at the stated rates for every degree-dependent schedule
    $\gamma_p>0$ held fixed during mesh refinement. This includes a weight
    reduced exponentially with the degree, $\gamma_p=\gamma_0e^{-cp}$: the
    cut-uniform constants then grow like $e^{cp/2}$ (energy) and $e^{cp}$
    ($L^2$), a factor that depends on $p$ alone and leaves the rates in $h$
    unchanged, while the ghost consistency terms shrink like $e^{-cp/2}$. For a mesh-dependent
    schedule, the bounds retain $M_{\gamma(h)}$: for example,
    $\gamma(h)=h^a$, $a>0$, gives the cut-uniform energy bound
    $Ch^{p-a/2}\|u\|_{H^{p+1}(\Omega)}$ and the $L^2$ bound
    $C(h^{p+1/2-a/2}+h^{p+1-a})\|u\|_{H^{p+1}(\Omega)}$.
    These bounds alone do not preserve the classical rates or guarantee
    convergence for arbitrarily rapid decay of $\gamma(h)$.
\end{remark}

\subsection{Monotonicity along the weight path}
\label{sec:weight-path}

The estimates so far bound what reducing the weight can \emph{cost}. They do not say how the graph norm and the
discrete solution themselves respond when the weight is lowered on a fixed configuration, nor in which direction the
residual weight dependence of the error points. All three questions have exact answers. The graph norms at two weights
are equivalent with a constant that is computable and follows the same two branches as~\eqref{eq:shedfactor}; the
solution moves along the weight path with a Lipschitz bound whose slope is again that factor; and below the saturation
threshold the only surviving weight dependence of the energy bound is the consistency term, which carries
$\gamma^{1/2}$ and points the \emph{opposite} way from the worst-case constant, so that lowering the weight lowers the
bound. Beyond the bound, the weight derivative of the discrete solution obeys exact identities
(Lemma~\ref{lem:weight-derivative}) which make the physical error itself monotone in one dimension and identify the
improvement as released regularization.

We begin with the norms, since every statement about the weight path is made in one of them. Lowering the weight
weakens the graph norm, so a quantity that shrinks as the weight decreases may be shrinking only because the norm
measuring it does. The following lemma settles how much the two norms can differ, and records that the bilinear form is
continuous with respect to the \emph{larger} of them with a constant independent of the smaller weight. Throughout,
$\gamma'\in(0,1]$ plays the role of the larger, or reference, weight, and we write
\begin{gather}\label{eq:theta}
    \theta:=\frac{\gamma'}{\gamma},
    \qquad
    \theta_{\mathrm{eff}}:=\frac{1+\gamma' R}{1+\gamma R},
    \qquad 0<\gamma\le\gamma' ,
\end{gather}
with $R$ the amplification quotient of~\eqref{eq:shedfactor}. That
quotient is finite because $n(v,v)>0$ for $v\ne0$ when $\tau_0>0$, by
Theorem~\ref{thm:coercivity}; it is cut-dependent, and $R\to\infty$ as the
cut degenerates.

\begin{lemma}[transfer between weights]\label{lem:transfer}
    Let $\tau_0>0$ and $0<\gamma\le\gamma'\le1$. Then
    \begin{gather}\label{eq:transfer-twosided}
        \theta_{\mathrm{eff}}^{-1}\,s_{\gamma'}(v,v)
        \ \le\ s_\gamma(v,v)\ \le\ s_{\gamma'}(v,v),
        \qquad
        \theta_{\mathrm{eff}}\ \le\ \min\bigl\{\theta,\ 1+\gamma' R\bigr\}.
    \end{gather}
    Equivalently,
    $\theta_{\mathrm{eff}}^{-1/2}\|v\|_{\gamma'}
        \le\|v\|_\gamma\le\|v\|_{\gamma'}$.
    The form $a_\gamma$ is continuous in the graph norm \emph{at the larger
        weight} with a constant that is free of $\gamma$ and of the cut:
    \begin{gather}\label{eq:transfer-continuity}
        |a_\gamma(v_h,w_h)|
        \ \le\
        \bar C_a\,\|v_h\|_{\gamma'}\|w_h\|_{\gamma'},
        \qquad
        \bar C_a
        \le
        C\Bigl(1+\tau_0^{-1/2}p^{-1}C_{\mathrm{tr}}(p)\,
        \bigl(C_E(p)/\gamma'\bigr)^{1/2}\Bigr).
    \end{gather}
\end{lemma}

\begin{proof}
    For~\eqref{eq:transfer-twosided}, fix $v\ne0$ and put
    $t:=g_h(v,v)/n(v,v)\in[0,R]$. Then
    \begin{gather*}
        \frac{s_\gamma(v,v)}{s_{\gamma'}(v,v)}
        =\frac{n(v,v)+\gamma g_h(v,v)}{n(v,v)+\gamma' g_h(v,v)}
        =\frac{1+\gamma t}{1+\gamma' t}=:f(t),
        \qquad
        f'(t)=\frac{\gamma-\gamma'}{(1+\gamma' t)^2}\le0 ,
    \end{gather*}
    so $f$ is nonincreasing and $f(R)\le f(t)\le f(0)=1$ on $[0,R]$, which
    is~\eqref{eq:transfer-twosided} with
    $\theta_{\mathrm{eff}}=1/f(R)$. For the two bounds on
    $\theta_{\mathrm{eff}}$: $1+\gamma R\ge1$ gives
    $\theta_{\mathrm{eff}}\le1+\gamma' R$, and
    $\gamma(1+\gamma' R)\le\gamma'(1+\gamma R)$, which reduces to
    $\gamma\le\gamma'$, gives $\theta_{\mathrm{eff}}\le\theta$. Only
    $t\le R$ enters, so no attainment of the supremum defining $R$ is
    needed.
    For~\eqref{eq:transfer-continuity}, bound the five terms
    of~\eqref{eq:aform} separately. Cauchy--Schwarz gives, for the volume
    term,
    $\|\nabla v_h\|_{0,\Omega}\|\nabla w_h\|_{0,\Omega}\le
        \|v_h\|_{\gamma'}\|w_h\|_{\gamma'}$, for the boundary-penalty
    term,
    $\tau_h^{1/2}\|v_h\|_{0,\Gamma}\,\tau_h^{1/2}\|w_h\|_{0,\Gamma}\le
        \|v_h\|_{\gamma'}\|w_h\|_{\gamma'}$, and for the ghost term,
    $\gamma\,g_h(v_h,v_h)^{1/2}g_h(w_h,w_h)^{1/2}\le\gamma'\,g_h(v_h,v_h)^{1/2}
        g_h(w_h,w_h)^{1/2}\le\|v_h\|_{\gamma'}\|w_h\|_{\gamma'}$, using
    $\gamma\le\gamma'$. For the consistency term
    $(\partial_nv_h,w_h)_\Gamma$ use the cut-cell trace
    inequality~\eqref{eq:cuttrace} and then~\eqref{eq:extension} \emph{at the
        larger weight},
    $h^{1/2}\|\partial_nv_h\|_{0,\Gamma}\le
        C_{\mathrm{tr}}(p)\|\nabla v_h\|_{0,\Omega_h}\le
        C_{\mathrm{tr}}(p)(C_E(p)/\gamma')^{1/2}\|v_h\|_{\gamma'}$,
    together with
    $h^{-1/2}\|w_h\|_{0,\Gamma}\le(\tau_0p(p+1))^{-1/2}\|w_h\|_{\gamma'}$
    from~\eqref{eq:tau}. The same two estimates with $v_h$ and $w_h$
    interchanged bound the consistency term $(v_h,\partial_nw_h)_\Gamma$,
    and the resulting factor two is absorbed into $C$.
\end{proof}

Two readings of the lemma matter later. First, \eqref{eq:transfer-continuity} says that \emph{the smaller weight no
    longer enters the continuity constant}: in Theorem~\ref{thm:continuity} the weight branch of $C_a$ grew like
$\gamma^{-1/2}$ because the graph norm was taken at $\gamma$ itself, whereas against the norm at a fixed larger weight
the estimate is $\gamma$-free, and all dependence on $\gamma$ is concentrated in $\theta_{\mathrm{eff}}$. Second, that
factor has the two regimes of $\min\{\theta,1+\gamma'R\}$. For a \emph{fixed} cut configuration and reference weight
$\gamma'$, it approaches $1+\gamma'R$ as $\gamma\to0$; for $R>0$ it does not become constant at any positive weight. At
fixed $0<\gamma<\gamma'$, one has $\theta_{\mathrm{eff}}<\theta$ for every finite $R$, and
$\theta_{\mathrm{eff}}\to\theta$ as $R\to\infty$. Thus $\theta$ is the supremum over any family of cut configurations
on which $R$ is unbounded, as in the slab model.

\begin{remark}[on the amplification quotient]\label{rem:amplification}
    This paper never invokes a criterion $\gamma R\le1$, and nothing here
    requires the ghost energy to be small compared with the physical
    energy. That route is not cut-uniform, since
    $R\gtrsim\kappa^{-(2p-1)}$ in the slab model, and a sufficiency
    criterion built on it would exclude exactly the configurations the
    method handles well. The quotient appears instead as a
    \emph{geometric cap}, rather than an admissibility criterion: in
    Lemma~\ref{lem:transfer}, $\theta_{\mathrm{eff}}\le1+\gamma'R$ is an
    asymptotic cap on weight transfer. In~\eqref{eq:shedfactor}, $1+R$
    is the second branch of $\Xi_\gamma$. On a fixed mesh and cut,
    $\Xi_\gamma$ and $M_\gamma$ become constant once
    $\gamma\le(1+R)^{-1}$; explicit factors such as $\gamma^{1/2}$ in the
    ghost residual and the exact transfer factor still depend on the weight.
    This should be quoted carefully: the level at which the estimates
    saturate is $1+R\gtrsim\kappa_{\min}^{-(2p-1)}$, which at $p=7$ is
    still far above anything reached in practice, so the saturation
    explains at most the \emph{shape} of the measured behavior, flat in
    $\gamma$ below a geometric threshold, and not its magnitude.
\end{remark}

\begin{lemma}[stability along the weight path]\label{lem:gamma-path}
    Let $\tau_0>0$ and $0<\gamma\le\gamma'\le1$, and let
    $u_h^{\gamma},u_h^{\gamma'}\in V_h^p$
    solve~\eqref{eq:discrete-problem} with the respective weights. Then
    $\delta:=u_h^{\gamma}-u_h^{\gamma'}$ satisfies
    \begin{gather}\label{eq:path-identity}
        a_\gamma(\delta,v_h)
        =(\gamma'-\gamma)\,g_h(u_h^{\gamma'},v_h),
        \qquad v_h\in V_h^p,
    \end{gather}
    and
    \begin{gather}\label{eq:path-lipschitz}
        \|\delta\|_\gamma
        \ \le\
        (\gamma'-\gamma)\,\min\{\gamma^{-1/2},R^{1/2}\}\,
        g_h(u_h^{\gamma'},u_h^{\gamma'})^{1/2}.
    \end{gather}
    Under the assumptions of Theorem~\ref{thm:energy}, moreover,
    \begin{gather}\label{eq:path-data}
        g_h(u_h^{\gamma'},u_h^{\gamma'})^{1/2}
        \ \le\
        C(p,\tau_0)\,(1+M_{\gamma'})\,\Xi_{\gamma'}^{1/2}\,
        h^{p}\,\|u\|_{H^{p+1}(\Omega)}.
    \end{gather}
\end{lemma}

\begin{proof}
    Both solutions satisfy $a_{\gamma}(u_h^\gamma,v_h)=\ell(v_h)$ and
    $a_{\gamma'}(u_h^{\gamma'},v_h)=\ell(v_h)$, and
    $a_\gamma=a_{\gamma'}-(\gamma'-\gamma)g_h$;
    subtraction gives~\eqref{eq:path-identity}. Testing with $v_h=\delta$
    and using the coercivity identity~\eqref{eq:coercivity-identity} and
    Cauchy--Schwarz for $g_h$,
    \begin{gather*}
        \|\delta\|_\gamma^2
        =(\gamma'-\gamma)\,g_h(u_h^{\gamma'},\delta)
        \le
        (\gamma'-\gamma)\,
        g_h(u_h^{\gamma'},u_h^{\gamma'})^{1/2}\,
        g_h(\delta,\delta)^{1/2}.
    \end{gather*}
    The unweighted ghost energy of $\delta$ is paid for in the two ways
    of~\eqref{eq:twobranch}: out of the weight,
    $g_h(\delta,\delta)\le\gamma^{-1}\|\delta\|_\gamma^2$, or out of the
    geometry, $g_h(\delta,\delta)\le R\,n(\delta,\delta)
        \le R\,\|\delta\|_\gamma^2$. Dividing by $\|\delta\|_\gamma$ (the case
    $\delta=0$ being trivial) gives~\eqref{eq:path-lipschitz}.
    For~\eqref{eq:path-data}, write
    $u_h^{\gamma'}=I_hu^e-\xi'$ with $\xi':=I_hu^e-u_h^{\gamma'}$;
    then~\eqref{eq:twobranch} and~\eqref{eq:xi-bound} give
    $g_h(\xi',\xi')^{1/2}\le\Xi_{\gamma'}^{1/2}\|\xi'\|_{\gamma'}
        \le C(p,\tau_0)(1+M_{\gamma'})\Xi_{\gamma'}^{1/2}h^p
        \|u\|_{H^{p+1}(\Omega)}$,
    and~\eqref{eq:interp-ghost} bounds $g_h(I_hu^e,I_hu^e)^{1/2}$ by
    $C(p)h^p\|u\|_{H^{p+1}(\Omega)}$; the triangle inequality for
    $g_h^{1/2}$ and $\Xi_{\gamma'}\ge1$ conclude.
\end{proof}

\begin{corollary}[monotone error bound below the saturation threshold]
    \label{cor:monotone}
    Fix the mesh and the cut, let $u\in H^{p+1}(\Omega)$, and let
    $0<\gamma\le\gamma'\le(1+R)^{-1}$. Write
    $M_R:=\min\{M_{\mathrm{geo}},\,M_{\mathrm{wt}}((1+R)^{-1})\}$, the saturated value
    of~\eqref{eq:Mgamma}. Then
    \begin{gather}\label{eq:monotone-bound}
        \|u^e-u_h^\gamma\|_\gamma
        \ \le\
        C(p,\tau_0)\,(1+M_R)\,h^{p}\,\|u\|_{H^{p+1}(\Omega)}
        \ +\
        2\,\gamma^{1/2}\,g_h(\eta,\eta)^{1/2},
    \end{gather}
    with both constants free of $\gamma$, so the right side is
    nondecreasing in $\gamma$. Writing $\mathcal B_\gamma$ for this right
    side, the physical errors satisfy
    \begin{gather*}
        \|\nabla(u-u_h^\gamma)\|_{0,\Omega}\le\mathcal B_\gamma,
        \qquad
        \|u-u_h^\gamma\|_{0,\Gamma}
        \le\tau_h^{-1/2}\mathcal B_\gamma.
    \end{gather*}
    Moreover,
    $g_h(\eta,\eta)^{1/2}\le C(p)h^p\|u\|_{H^{p+1}(\Omega)}$
    by~\eqref{eq:interp-ghost}. Moreover
    \begin{gather}\label{eq:path-frozen}
        \|u_h^{\gamma}-u_h^{\gamma'}\|_\gamma
        \ \le\
        C(p,\tau_0)\,(\gamma'-\gamma)(1+R)\,(1+M_R)\,
        h^{p}\,\|u\|_{H^{p+1}(\Omega)}
        \ \le\
        C(p,\tau_0)\,(1+M_R)\,h^{p}\,\|u\|_{H^{p+1}(\Omega)}:
    \end{gather}
    across the whole sub-threshold range the discrete solutions differ from
    one another by no more than the error floor itself, linearly in the
    weight decrement.
\end{corollary}

\begin{proof}
    For $\gamma\le(1+R)^{-1}$ one has $\Xi_\gamma=1+R$ and
    $M_\gamma=M_R$. In the proof of Theorem~\ref{thm:energy} keep the
    ghost contribution separate: the first four terms give
    $\|\xi\|_\gamma\le C(p,\tau_0)(1+M_R)h^p\|u\|_{H^{p+1}(\Omega)}
        +\gamma^{1/2}g_h(\eta,\eta)^{1/2}$, and splitting the graph norm of
    the interpolation error as
    $\|\eta\|_\gamma\le
        \bigl(\|\nabla\eta\|^2_{0,\Omega}
        +\tau_h\|\eta\|^2_{0,\Gamma}\bigr)^{1/2}
        +\gamma^{1/2}g_h(\eta,\eta)^{1/2}$
    and using~\eqref{eq:interp} yields~\eqref{eq:monotone-bound}.
    For~\eqref{eq:path-frozen}, $\gamma^{-1}\ge1+R>R$ makes the geometric
    branch of~\eqref{eq:path-lipschitz} the smaller one, and
    $\Xi_{\gamma'}=1+R$ in~\eqref{eq:path-data}, so
    $\|\delta\|_\gamma\le(\gamma'-\gamma)R^{1/2}(1+R)^{1/2}
        \,C(p,\tau_0)(1+M_R)h^p\|u\|_{H^{p+1}(\Omega)}$; bound
    $R^{1/2}(1+R)^{1/2}\le1+R$ and, for the second inequality,
    $(\gamma'-\gamma)(1+R)\le1$.
\end{proof}

\begin{remark}[a smaller weight improves the estimate until it saturates]
    \label{rem:improvement}
    Read at face value, Theorem~\ref{thm:energy} says the error constant \emph{grows} as the weight decreases, through
    $1+M_\gamma\le C(p)\gamma^{-1/2}$. That branch is worst-case in the cut; in computations one commonly observes the
    opposite, a mild \emph{improvement} of the error when the weight is moderately reduced. The two statements are
    compatible, and Corollary~\ref{cor:monotone} locates the mechanism. On a fixed configuration the $\gamma^{-1/2}$ branch
    is inactive below the threshold $(1+R)^{-1}$, and the surviving weight dependence of the bound has the opposite sign:
    it is the consistency term $\gamma^{1/2}g_h(\eta,\eta)^{1/2}$, which shrinks with the weight, so the bound decreases
    monotonically as the weight decreases until it reaches the weight-free floor of~\eqref{eq:monotone-bound}. The effect
    changes the coefficient in this bound at full regularity. Below full regularity, Remark~\ref{rem:lowreg} bounds the
    corresponding discrete ghost residual by $C(p,s)\gamma^{1/2}h^{s-1}\|u\|_{H^s(\Omega)}$; no comparison of its actual
    magnitude across regularity classes is asserted. A monotone upper bound does not make the error itself monotone;
    what~\eqref{eq:path-frozen} adds is that below the threshold the error can respond to the weight at most linearly, with
    slope $(1+R)(1+M_R)h^p$, so its variation over the whole sub-threshold range is confined to the scale of the error
    bound.

    Comparisons across weights must be made in weight-free quantities: the graph norm itself weakens as $\gamma$ falls, by
    at most the factor $\sqrt{\theta_{\mathrm{eff}}}$ of Lemma~\ref{lem:transfer}, so an error decay measured in
    $\|\cdot\|_\gamma$ is partly an artifact of the norm, which is why Corollary~\ref{cor:monotone} is stated for the
    physical errors~\eqref{eq:physical-rates} as well.
\end{remark}

The bound of Corollary~\ref{cor:monotone} is monotone; the error itself obeys exact identities. Write
\begin{gather}\label{eq:skew}
    k(v,w):=-(\partial_nv,w)_\Gamma+(v,\partial_nw)_\Gamma,
    \qquad
    a_\gamma=s_\gamma+k,
\end{gather}
for the skew Nitsche pair, $s_\gamma$ being the symmetric part
of~\eqref{eq:graphnorm} read as a bilinear form, so that $k(v,v)=0$ for
every $v$ with the required traces. On a fixed mesh and cut the discrete
solution $u_h^\gamma$ is a rational function of $\gamma$, the solution of a
linear system whose matrix depends affinely on $\gamma$ and is nonsingular
for every $\gamma>0$ by Theorem~\ref{thm:coercivity}, and we write
$\dot u_h:=\frac{d}{d\gamma}u_h^\gamma\in V_h^p$ for its weight derivative.

\begin{lemma}[weight derivatives]\label{lem:weight-derivative}
    Let $\tau_0>0$, $\gamma>0$, $u\in H^{p+1}(\Omega)$, and $e:=u^e-u_h^\gamma$.
    Then, for all $v_h\in V_h^p$,
    \begin{gather}\label{eq:udot}
        a_\gamma(\dot u_h,v_h)=-g_h(u_h^\gamma,v_h),
    \end{gather}
    and
    \begin{gather}\label{eq:ghost-derivative}
        \frac{d}{d\gamma}\,g_h(u_h^\gamma,u_h^\gamma)
        =-2\,\|\dot u_h\|_\gamma^2\ \le\ 0,
        \\
        \label{eq:physical-derivative}
        \frac{d}{d\gamma}\,n(e,e)
        =2\gamma\,\|\dot u_h\|_\gamma^2+2\,k(e,\dot u_h),
        \qquad
        k(e,\dot u_h)=\tfrac12\,g_h(u_h^\gamma,q_h),
    \end{gather}
    where $q_h\in V_h^p$ is the adjoint Ritz projection of the error, the
    unique discrete function with $a_\gamma(w_h,q_h)=a_\gamma(w_h,e)$ for
    all $w_h\in V_h^p$.
\end{lemma}

\begin{proof}
    Differentiating $a_\gamma(u_h^\gamma,v_h)=\ell(v_h)$ with respect to
    $\gamma$, with $a_\gamma=a_0+\gamma g_h$ and $\ell$ independent of the
    weight, gives~\eqref{eq:udot}. Then
    $\frac{d}{d\gamma}g_h(u_h^\gamma,u_h^\gamma)=2g_h(u_h^\gamma,\dot u_h)
        =-2a_\gamma(\dot u_h,\dot u_h)=-2\|\dot u_h\|^2_\gamma$
    by~\eqref{eq:udot} with $v_h=\dot u_h$ and the coercivity
    identity~\eqref{eq:coercivity-identity}. For the physical energy,
    $\frac{d}{d\gamma}n(e,e)=-2n(e,\dot u_h)$, and Galerkin
    orthogonality~\eqref{eq:orthogonality} tested with $\dot u_h$ reads
    $n(e,\dot u_h)=-k(e,\dot u_h)-\gamma g_h(e,\dot u_h)$, where
    $g_h(e,\dot u_h)=-g_h(u_h^\gamma,\dot u_h)=\|\dot u_h\|^2_\gamma$ by
    ghost consistency~\eqref{eq:ghost-consistency} and~\eqref{eq:udot}.
    This is the first identity in~\eqref{eq:physical-derivative}. For the
    second, $q_h$ exists and is unique because $a_\gamma$ is nondegenerate
    on $V_h^p$ by Theorem~\ref{thm:coercivity}. Since $s_\gamma$ is
    symmetric and $k$ skew,
    $a_\gamma(\dot u_h,e)-a_\gamma(e,\dot u_h)=-2k(e,\dot u_h)$, and
    $a_\gamma(e,\dot u_h)=0$, so
    $k(e,\dot u_h)=-\tfrac12a_\gamma(\dot u_h,e)
        =-\tfrac12a_\gamma(\dot u_h,q_h)=\tfrac12g_h(u_h^\gamma,q_h)$
    by the definition of $q_h$ and~\eqref{eq:udot}.
\end{proof}

\begin{corollary}[exact monotonicity in one dimension]\label{cor:onedim}
    In the one-dimensional analogue of the setting, with $\Omega$ an interval, let the active mesh be a chain of cells of
    length $h$ and $\Gamma$ either both endpoints of $\Omega$ or one of them, with a natural boundary condition at the
    other. Boundary pairings are sums over the endpoints in $\Gamma$, and the ghost penalty is
    \begin{gather}\label{eq:ghost-1d}
        g_h(v,w)=\sum_{F\in\mathcal F_G}\sum_{k=1}^p
        \frac{h^{2k-1}}{(k!)^2}
        \llbracket v^{(k)}\rrbracket_F\llbracket w^{(k)}\rrbracket_F,
    \end{gather}
    where $F$ denotes an interior mesh node with at least one cut
    neighbor. Under the hypotheses of Lemma~\ref{lem:weight-derivative},
    $k(e,\dot u_h)=0$ for every $\gamma>0$. Hence, on a fixed mesh and cut, the physical error energy
    \begin{gather}\label{eq:onedim-monotone}
        \mathcal E(\gamma):=n(e,e)
        =\|\nabla(u-u_h^\gamma)\|^2_{0,\Omega}
        +\tau_h\|u-u_h^\gamma\|^2_{0,\Gamma}
    \end{gather}
    is nondecreasing on $(0,\infty)$, and the unweighted ghost energy $g_h(u_h^\gamma,u_h^\gamma)$ of the discrete solution
    is nonincreasing. Decreasing the weight therefore never increases the error and never
    decreases the ghost energy.
\end{corollary}

\begin{proof}
    By Lemma~\ref{lem:weight-derivative} it suffices to show that the
    adjoint Ritz projection $q_h$ of the error is affine, since affine
    functions have no derivative jumps and then $g_h(u_h^\gamma,q_h)=0$.
    Let $x_0<x_1$ be the endpoints of $\Omega$, with outward normals
    $n=\mp1$, and let $q(x)=\alpha+\beta x$ be the affine function with
    $q=2e$ at both endpoints, so $\beta=2(e(x_1)-e(x_0))/(x_1-x_0)$. Affine
    functions belong to $V_h^p$ and carry no ghost jumps, so for
    $w_h\in V_h^p$, using $(w_h',\beta)_\Omega=\beta(w_h,n)_\Gamma$ and
    $\partial_nq=\beta n$,
    \begin{gather*}
        a_\gamma(w_h,q)
        =2\beta(w_h,n)_\Gamma-(\partial_nw_h,q)_\Gamma+\tau_h(w_h,q)_\Gamma .
    \end{gather*}
    Galerkin orthogonality~\eqref{eq:orthogonality} tested with the two
    functions $1$ and $x$, using $(e',1)_\Omega=(e,n)_\Gamma$ for the
    continuous function $e$, gives
    \begin{gather*}
        (\partial_ne-\tau_he,\,1)_\Gamma=0,
        \qquad
        (\partial_ne-\tau_he,\,x)_\Gamma=2(e,n)_\Gamma ,
    \end{gather*}
    two linear equations for the two values of $r:=\partial_ne-\tau_he$ on
    $\Gamma$, whose solution is $r=\beta n$. Hence
    $2k(w_h,e)=-2(\partial_nw_h,e)_\Gamma+2(w_h,\tau_he+\beta n)_\Gamma
        =a_\gamma(w_h,q)$. On the other hand
    $a_\gamma(w_h,e)=a_\gamma(e,w_h)+2k(w_h,e)=2k(w_h,e)$ by Galerkin
    orthogonality. Thus $a_\gamma(w_h,e-q)=0$ for all $w_h\in V_h^p$, and
    $q_h=q$.
    If $\Gamma=\{x_1\}$, Galerkin orthogonality tested with the constant $1$ alone gives $r=0$ at $x_1$, and the same
    computation with the constant $q=2e(x_1)$, that is $\beta=0$, gives $a_\gamma(w_h,e-q)=0$.
\end{proof}

\begin{remark}[the ghost penalty as a regularizer]\label{rem:tikhonov}
    The first term of~\eqref{eq:physical-derivative} is nonnegative and describes the contribution of the symmetric part to
    the change in physical error energy. The second term is the adjoint-inconsistency correction, a ghost pairing of the
    discrete solution with the adjoint Ritz projection of its error. In one dimension this correction vanishes: decreasing
    the weight reduces the combined physical error energy while increasing the unweighted ghost energy. In general
    dimension the ghost energy still increases as the weight decreases by~\eqref{eq:ghost-derivative}, but the correction
    has no controlled sign, so the identity does not imply monotone physical accuracy.

    For the symmetric Nitsche sign, the ghost penalty has a Tikhonov regularization interpretation. Write
    $a_0^{\mathrm{sym}}$ for the physical symmetric form. Whenever the total stabilized form $a_0^{\mathrm{sym}}+\gamma
        g_h$ is positive definite on $V_h^p$, consistency and ghost consistency imply that the symmetric discrete solution
    minimizes $a_0^{\mathrm{sym}}(u^e-v,u^e-v)+\gamma g_h(v,v)$ over $V_h^p$. Comparing the minimizers at two such weights
    shows that the physical symmetric quadratic form $a_0^{\mathrm{sym}}(e,e)$ is nondecreasing in $\gamma$, while the
    ghost energy is nonincreasing. The former quadratic form need not be nonnegative and is not asserted to be an error
    norm. Below the coercivity floor of the total form, the minimization characterization is lost; the dashed curves of
    Figure~\ref{fig:weight-path} illustrate a near-singular spike at this transition.

    Although $u^e$ has zero ghost energy, its interpolant can have nonzero derivative jumps.
    Estimate~\eqref{eq:interp-ghost} gives $g_h(I_hu^e,I_hu^e)\le C(p)h^{2p}\|u\|_{H^{p+1}(\Omega)}^2$. Thus the ghost
    contribution can enter the squared error bound at the same order as the approximation contribution. This upper bound
    supplies neither a lower bound on the ghost energy nor a comparison between the ghost energies of the interpolant and
    the discrete solution. At each fixed positive weight, the estimates retain the stated convergence order with a
    weight-dependent constant; they do not determine whether special cases have higher actual rates.

    In the model of Figure~\ref{fig:weight-path}, reducing the weight improves the cut-cell error by about three orders of
    magnitude before it saturates, while the interior-cell error is observed to remain essentially unchanged. This
    concentration of the improvement near the boundary is an observation for that model. Although the ghost penalty acts on
    faces in the boundary zone, the equation for $\dot u_h$ is a global discrete problem, so the identities do not confine
    changes in the solution to that zone.
\end{remark}

Figure~\ref{fig:weight-path} traces the error along the weight path in a solvable model: two cells of size one, the
first fully inside, the second cut at fraction $\kappa=10^{-4}$, one ghost face between them carrying the all-order
penalty~\eqref{eq:ghost-1d}, boundary penalty $5p(p+1)$ at the cut end, and the manufactured solution $u(x)=\sin(x+1)$,
smooth but not a polynomial, so every degree carries a nonzero approximation error. The uncut end carries a natural
boundary condition: the full cell stands for its connection to the rest of the mesh, and a Dirichlet condition there
would stabilize the model artificially. The error is measured in the $H^1$-seminorm of the physical part of the cut
cell, where the weight acts. As the weight decreases, the nonsymmetric error falls by about three orders of magnitude,
somewhat less for the higher degrees, and then levels off. Corollary~\ref{cor:onedim} proves that the physical
energy~\eqref{eq:onedim-monotone} of this model is monotone in the weight; the plotted cut-cell seminorm is not: for
every degree it passes through a shallow dip, about a factor two below the level at which it settles. At small weights
the symmetric error is nearly the same, but on the way it passes through the loss of definiteness of the symmetric
form, a near-singular spike at its ghost-weight floor; the nonsymmetric path has no such obstacle at any positive
weight, by Theorem~\ref{thm:coercivity}.

Figure~\ref{fig:weight-path}(b) shows the nonsymmetric solutions on both cells, including the part of the cut cell
outside $\Omega$, where the discrete solution is the polynomial of the physical part continued across the cell.
Reducing the weight to $\gamma_p$ changes this extension little: outside $\Omega$ it departs from the extension of $u$
by the extrapolation error of the polynomial, which decreases with the degree. Far below $\gamma_p$ the extension stays
bounded and moves, if anything, closer to that of $u$.

\begin{figure}[htbp]
    \centering
    \pgfplotsset{weight path/.style={paperplot, width=\linewidth, height=0.85\linewidth, mark size=1.5pt}}
    \begin{subfigure}[t]{0.49\textwidth}
        \centering
        \begin{tikzpicture}
            \begin{loglogaxis}[
                    weight path,
                    x dir=reverse,
                    xlabel={ghost weight $\gamma$},
                    ylabel={$H^1$-seminorm of the error on $\cell\cap\Omega$},
                    legend to name=weightpathlegend,
                    legend columns=-1,
                    legend style={font=\scriptsize, /tikz/every even column/.append style={column sep=6pt}},
                    mark repeat=40,
                    xmin=1e-7, xmax=1,
                    ymin=1e-13, ymax=1e0,
                    enlarge x limits=false,
                    enlarge y limits=false,
                    xtick={1e0,1e-2,1e-4,1e-6},
                ]
                \addplot[p2] table[x=gamma, y=p2, col sep=comma]
                    {numerical/results/error_floor_compare_nonsym.csv};
                \addlegendentry{$\Q_2$}
                \addplot[p3] table[x=gamma, y=p3, col sep=comma]
                    {numerical/results/error_floor_compare_nonsym.csv};
                \addlegendentry{$\Q_3$}
                \addplot[p4] table[x=gamma, y=p4, col sep=comma]
                    {numerical/results/error_floor_compare_nonsym.csv};
                \addlegendentry{$\Q_4$}
                \addplot[p5] table[x=gamma, y=p5, col sep=comma]
                    {numerical/results/error_floor_compare_nonsym.csv};
                \addlegendentry{$\Q_5$}
                \addplot[p7] table[x=gamma, y=p7, col sep=comma]
                    {numerical/results/error_floor_compare_nonsym.csv};
                \addlegendentry{$\Q_7$}
                \addplot[p2, densely dashed, mark options={solid, fill=white}]
                table[x=gamma, y=p2, col sep=comma]
                    {numerical/results/error_floor_compare_sym.csv};
                \addplot[p3, densely dashed, mark options={solid, fill=white}]
                table[x=gamma, y=p3, col sep=comma]
                    {numerical/results/error_floor_compare_sym.csv};
                \addplot[p4, densely dashed, mark options={solid, fill=white}]
                table[x=gamma, y=p4, col sep=comma]
                    {numerical/results/error_floor_compare_sym.csv};
                \addplot[p5, densely dashed, mark options={solid, fill=white}]
                table[x=gamma, y=p5, col sep=comma]
                    {numerical/results/error_floor_compare_sym.csv};
                \addplot[p7, densely dashed, mark options={solid, fill=white}]
                table[x=gamma, y=p7, col sep=comma]
                    {numerical/results/error_floor_compare_sym.csv};
                \addplot[teal,   solid, semithick, opacity=0.5, no markers] coordinates
                    {(3.83e-2,1e-13) (3.83e-2,1e0)};
                \addplot[p2, only marks, mark options={solid, fill=white}] coordinates {(3.83e-2,1e0)};
                \addplot[degQthree, solid, semithick, opacity=0.5, no markers] coordinates
                    {(1.95e-2,1e-13) (1.95e-2,1e0)};
                \addplot[p3, only marks, mark options={solid, fill=white}] coordinates {(1.95e-2,0.3)};
                \addplot[red,    solid, semithick, opacity=0.5, no markers] coordinates
                    {(1.18e-2,1e-13) (1.18e-2,1e0)};
                \addplot[p4, only marks, mark options={solid, fill=white}] coordinates {(1.18e-2,1e0)};
                \addplot[violet, solid, semithick, opacity=0.5, no markers] coordinates
                    {(7.86e-3,1e-13) (7.86e-3,1e0)};
                \addplot[p5, only marks, mark options={solid, fill=white}] coordinates {(7.86e-3,0.3)};
                \addplot[degQseven, solid, semithick, opacity=0.5, no markers] coordinates
                    {(4.19e-3,1e-13) (4.19e-3,1e0)};
                \addplot[p7, only marks, mark options={solid, fill=white}] coordinates {(4.19e-3,1e0)};
            \end{loglogaxis}
        \end{tikzpicture}
        \caption{Cut-cell error against the weight.}
    \end{subfigure}\hfill
    \begin{subfigure}[t]{0.49\textwidth}
        \centering
        \begin{tikzpicture}
            \begin{axis}[
                    weight path,
                    xlabel={$x$},
                    ylabel={$u_h$},
                    xmin=-1, xmax=1, ymin=-0.05, ymax=1.3,
                    enlarge x limits=false,
                    mark repeat=40,
                ]
                \fill[gray!12] (axis cs:1e-4,-0.05) rectangle (axis cs:1,1.3);
                \addplot[p2] table[x=x, y=p2_gp, col sep=comma] {numerical/results/onedim_profile.csv};
                \addplot[p3] table[x=x, y=p3_gp, col sep=comma] {numerical/results/onedim_profile.csv};
                \addplot[p4] table[x=x, y=p4_gp, col sep=comma] {numerical/results/onedim_profile.csv};
                \addplot[p5] table[x=x, y=p5_gp, col sep=comma] {numerical/results/onedim_profile.csv};
                \addplot[p7] table[x=x, y=p7_gp, col sep=comma] {numerical/results/onedim_profile.csv};
                \addplot[p2, densely dashed, mark phase=21, mark options={solid, fill=white}]
                table[x=x, y=p2_g1, col sep=comma] {numerical/results/onedim_profile.csv};
                \addplot[p3, densely dashed, mark phase=21, mark options={solid, fill=white}]
                table[x=x, y=p3_g1, col sep=comma] {numerical/results/onedim_profile.csv};
                \addplot[p4, densely dashed, mark phase=21, mark options={solid, fill=white}]
                table[x=x, y=p4_g1, col sep=comma] {numerical/results/onedim_profile.csv};
                \addplot[p5, densely dashed, mark phase=21, mark options={solid, fill=white}]
                table[x=x, y=p5_g1, col sep=comma] {numerical/results/onedim_profile.csv};
                \addplot[p7, densely dashed, mark phase=21, mark options={solid, fill=white}]
                table[x=x, y=p7_g1, col sep=comma] {numerical/results/onedim_profile.csv};
                \node[font=\scriptsize, anchor=south east] at (axis cs:0.97,0) {outside $\Omega$};
            \end{axis}
        \end{tikzpicture}
        \caption{Discrete solutions of the nonsymmetric form on the active mesh.}
    \end{subfigure}

    \smallskip
    \pgfplotslegendfromname{weightpathlegend}
    \caption{The two-cell model at cut fraction $\kappa=10^{-4}$, with a natural boundary condition at the uncut end.
        (a) $H^1$-seminorm of the error on the physical part of the cut cell against the ghost weight $\gamma$,
        decreasing to the right. Solid lines: nonsymmetric form~\eqref{eq:aform}. Dashed lines: symmetric form, with
        its floor~\eqref{eq:sym-floor-nonsymm} marked by a vertical line. (b) Discrete solutions of the nonsymmetric form,
        including the part of the cut cell outside $\Omega$ (shaded), at $\gamma_p$ of~\eqref{eq:reduced-weight}
        (solid) and at $\gamma=1$ (dashed).}
    \label{fig:weight-path}
\end{figure}

\section{Numerical experiments}
\label{sec:numerics}
The experiments in this section concern the discretization; the solver is the subject of Section~\ref{sec:multigrid}.
They address three questions raised by the analysis: whether the classical rates persist at every ghost weight held
fixed during refinement, as Theorems~\ref{thm:energy} and~\ref{thm:ltwo} predict; how the error responds when the
weight is reduced below the coercivity floor of the symmetric Nitsche form; and which $L^2$ rate is observed, the
optimal $h^{p+1}$ or the proved $h^{p+1/2}$ (Remark~\ref{rem:halforder}). All errors are measured in the weight-free
quantities $|u-u_h|_{1,\Omega}$ and $\|u-u_h\|_{0,\Omega}$. The graph norm $\|\cdot\|_\gamma$ changes with the weight,
and comparing errors across weights in that norm would confound the error with the norm.

\subsection{Test problem and meshes}
\label{sec:num-setting}

The physical domain $\Omega$ is the unit disc centered at the origin, embedded in the background square
$[-1.21,1.21]^2$ and described by the level-set function $|x|-1$. The background mesh on level $\ell$ has $2^\ell$
cells in each coordinate direction, so $h=2.42\cdot2^{-\ell}$, and the circle cuts each level at a different set of
volume fractions. The manufactured solution is $u=\sin(\pi x)\cos(\pi y)$, with $f=-\Delta u$ and $g=u|_\Gamma$; it is
smooth but not polynomial, so the approximation error is nonzero at every degree. The discrete space is the
Gauss--Lobatto nodal $\Q_p$ space on the active cells, and the boundary penalty is~\eqref{eq:tau} with $\tau_0=5$. On
each level the level-set function is interpolated in the same $\Q_p$ space, and the computed domain is the zero set of
the interpolant. Cut cells and their boundary parts are integrated by the algorithm of Saye~\cite{Saye2015} applied to
the interpolated level set and built on the $(p+3)$-point Gauss rule, and the uncut cells use the tensor-product
$(p+3)$-point Gauss rule; the ghost-penalty faces use the $(p+1)$-point Gauss rule, which is exact for the polynomial
integrands there. The analysis assumes exact geometry and quadrature and does not cover this perturbation. Unless
stated otherwise, the linear systems are solved by GMRES preconditioned with the multigrid method of
Section~\ref{sec:multigrid}, to a residual reduction of $10^{-11}$. Every iterative solve reported in this section
reaches this reduction; the largest counts, for $\Q_5$ at $\gamma=1$, slightly exceed $100$ iterations, in agreement
with Table~\ref{tab:mg-weights}. At high degree the discretization error approaches machine precision, and with a
common tolerance such as $10^{-9}$ the algebraic error would exceed it on most meshes; where the discretization error
falls to the level of the algebraic error, the computed error stalls. The weight sweep of Section~\ref{sec:num-floor}
and the orders of Section~\ref{sec:num-halforder} measure the discretization alone, down to errors below this level,
and are computed with a sparse direct solver.

All computations use the finite element library deal.II~\cite{dealII97}. The implementation is not optimized: the
matrix of every local problem~\eqref{eq:local-smoother} is inverted exactly and the inverse is stored for each patch,
so memory, not computing time, limits the largest problems in Section~\ref{sec:num-potato}.

A cell is \emph{cut} if the quadrature generated by Saye's algorithm assigns positive measure to both its part inside
and its part outside $\Omega$. A cell with an empty outside part is integrated as an uncut cell, and a cell with an
empty inside part is not active. Coefficients of both signs in the Bernstein expansion of the interpolated level set on
a cell do not imply that its zero set enters the cell, so the sign test alone marks as cut some cells that lie on one
side of the boundary; Section~\ref{sec:num-potato} shows where this matters.

The two additional points are needed although the cut quadrature is generated by Saye's algorithm. On a cut cell the
physical part is bounded by the curved zero set of the interpolated level set, and a rule built on the $q$-point Gauss
rule is then not exact for polynomial integrands of degree $2q-1$, as it is on a full cell. With $q=p+1$ the resulting
violation of the discrete Green identity is a consistency error which, in the patch test below, decays at the same
orders as the discretization error, $h^{p}$ in the $H^1$-seminorm and $h^{p+1}$ in the $L^2$ norm; it changes the error
constant without changing the rate. In that patch test, for $p=2$, $\gamma=1$ and a quadratic $u$, which the method
must reproduce exactly, the $(p+1)$-point rule leaves an $L^2$ error far above the algebraic error on every level,
while the $(p+3)$-point rule reproduces $u$ up to the algebraic error. Using a linear instead of a quadratic
interpolant of the level set gives the same picture, so the effect comes from the quadrature and not from the degree of
the geometry approximation. For all degrees $p=1,\dots,5$, with $u=x^p+x^{p-1}y+y^p$, the $(p+3)$-point rule reproduces
$u$ up to the algebraic error from level $5$ on, at $\gamma=1$ as well as at $\gamma=10^{-12}$. On the two coarsest
levels a very small weight exposes a remaining quadrature error, which the $(p+5)$-point rule removes.

Besides fixed weights, the runs use the \emph{reduced weight}
\begin{gather}\label{eq:reduced-weight}
    \gamma_p:=\gamma_0\,e^{-cp},
    \qquad
    \gamma_0=\tfrac29\,e^{4}\approx12.1,\quad c=4,
\end{gather}
which decreases exponentially with the degree, does not depend on the mesh or the cut, and is held fixed during
refinement, so the estimates of Section~\ref{sec:error} apply with the constants of Remark~\ref{rem:weight-once}.
The constant $\gamma_0$ places $\gamma_1=2/9$ at twice the floor $\gamma_{\mathrm{lim}}(1)=1/9$
of~\eqref{eq:sym-floor-nonsymm}, below which the symmetric form loses cut-uniform stability, and each further degree
divides the weight by $e^{4}\approx55$. Table~\ref{tab:reduced-weight} lists $\gamma_p$ together with
$\gamma_{\mathrm{lim}}(p)$; from $\Q_2$ on, $\gamma_p$ lies below that floor, where only the nonsymmetric form is
covered by a stability estimate.

\begin{table}[htbp]
    \centering
    \caption{The reduced weight $\gamma_p$ of~\eqref{eq:reduced-weight} and the floor $\gamma_{\mathrm{lim}}(p)$
        of~\eqref{eq:sym-floor-nonsymm} of the symmetric form, at $\tau_0=5$.}
    \label{tab:reduced-weight}
    \begin{tabular}{lccccccc}
        \toprule
                                   & $\Q_1$ & $\Q_2$        & $\Q_3$        & $\Q_4$         & $\Q_5$         & $\Q_6$         & $\Q_7$         \\
        \midrule
        $\gamma_p$                 & $2/9$  & $2/9\,e^{-4}$ & $2/9\,e^{-8}$ & $2/9\,e^{-12}$ & $2/9\,e^{-16}$ & $2/9\,e^{-20}$ & $2/9\,e^{-24}$ \\
        $\gamma_{\mathrm{lim}}(p)$ & $1/9$  & $1/26$        & $1/51$        & $1/84$         & $1/125$        & $1/174$        & $1/231$        \\
        \bottomrule
    \end{tabular}
\end{table}

Table~\ref{tab:num-meshes} lists, for each level, the number of cut cells and the smallest physical volume fraction
$\kappa_{\min}=\min_{\cell}|\cell\cap\Omega|/|\cell|$ over the active cells. The amplification quotient $R$
of~\eqref{eq:shedfactor} grows like $\kappa_{\min}^{-(2p-1)}$ in the slab model. We compute it on levels $3$--$7$ as
the largest eigenvalue of $g_h$ relative to $n$, by the Lanczos method. For $p\ge2$, on every mesh on which $R$ is
computable, $(1+R)^{-1}$ lies below every weight tested, so these runs lie on the weight branch of $\Xi_\gamma$
(Remark~\ref{rem:amplification}). For $p=1$ the threshold $(1+R)^{-1}$ is moderate, and the $\Q_1$ runs at the smallest
weights lie on the saturated branch. For high degrees on the finer levels, $R$ is too large to be computed reliably in
double precision.

\begin{table}[htbp]
    \centering
    \caption{Test meshes for the unit disc in $[-1.21,1.21]^2$: level $\ell$, mesh size $h$, number of active and of
        cut cells, and smallest physical volume fraction $\kappa_{\min}$ of an active cell, for the level set interpolated
        in $\Q_p$, $p\ge3$. For $p=1$ and $p=2$ the interpolated circle lowers $\kappa_{\min}$ on levels $3$--$7$ by up
        to $22\%$ and $1\%$, respectively.}
    \label{tab:num-meshes}
    \begin{tabular}{ccccc}
        \toprule
        $\ell$ & $h$                & active cells & cut cells & $\kappa_{\min}$    \\
        \midrule
        $3$    & $3.03\cdot10^{-1}$ & $52$         & $28$      & $3.11\cdot10^{-2}$ \\
        $4$    & $1.51\cdot10^{-1}$ & $164$        & $52$      & $4.40\cdot10^{-2}$ \\
        $5$    & $7.56\cdot10^{-2}$ & $608$        & $108$     & $1.54\cdot10^{-2}$ \\
        $6$    & $3.78\cdot10^{-2}$ & $2316$       & $212$     & $8.25\cdot10^{-4}$ \\
        $7$    & $1.89\cdot10^{-2}$ & $8992$       & $420$     & $6.79\cdot10^{-5}$ \\
        $8$    & $9.45\cdot10^{-3}$ & $35580$      & $844$     & $2.72\cdot10^{-4}$ \\
        \bottomrule
    \end{tabular}
\end{table}
 
\subsection{Convergence at fixed weights}
\label{sec:num-rates}

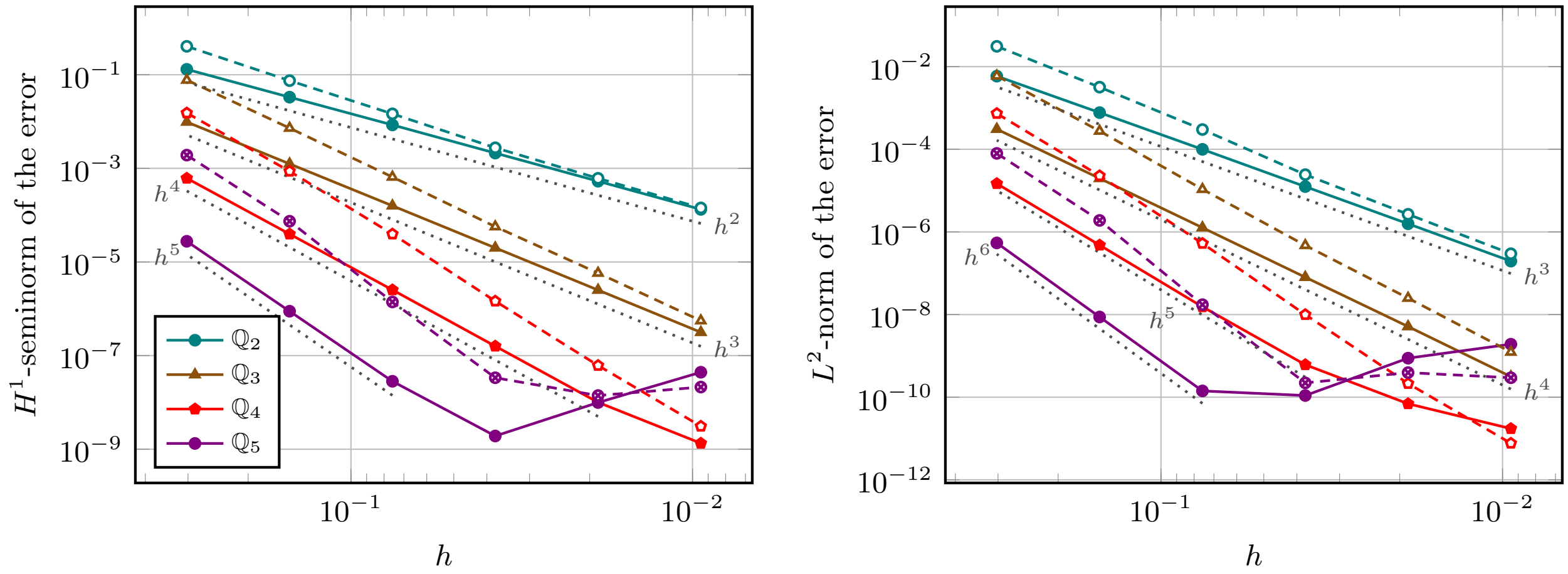
\begin{figure}[htbp]
    \centering
    \def\ratesdata{numerical/results/nonsymm_paper_plots/rates.csv}
    \pgfplotsset{num rates/.style={paperplot, xlabel={$h$}, x dir=reverse, legend pos=south west,
                legend cell align=left, legend style={font=\scriptsize}, mark size=1.5pt},
        g1/.style={densely dashed, mark options={solid, fill=white}, forget plot},
        slope/.style={black!70, dotted, thick, samples=2, forget plot}}
    \begin{tikzpicture}
        \begin{loglogaxis}[num rates, ylabel={$H^1$-seminorm of the error}]
            \addplot[p2] table[x=h, y=H1_p2_gp, col sep=comma] {\ratesdata}; \addlegendentry{$\Q_2$}
            \addplot[p3] table[x=h, y=H1_p3_gp, col sep=comma] {\ratesdata}; \addlegendentry{$\Q_3$}
            \addplot[p4] table[x=h, y=H1_p4_gp, col sep=comma] {\ratesdata}; \addlegendentry{$\Q_4$}
            \addplot[p5] table[x=h, y=H1_p5_gp, col sep=comma] {\ratesdata}; \addlegendentry{$\Q_5$}
            \addplot[p2, g1] table[x=h, y=H1_p2_g1, col sep=comma] {\ratesdata};
            \addplot[p3, g1] table[x=h, y=H1_p3_g1, col sep=comma] {\ratesdata};
            \addplot[p4, g1] table[x=h, y=H1_p4_g1, col sep=comma] {\ratesdata};
            \addplot[p5, g1] table[x=h, y=H1_p5_g1, col sep=comma] {\ratesdata};
            \addplot[slope, domain=0.009453125:0.3025] {0.7439*x^2}
            node[pos=0, anchor=west, font=\scriptsize] {$h^2$};
            \addplot[slope, domain=0.009453125:0.3025] {0.1862*x^3}
            node[pos=0, anchor=west, font=\scriptsize] {$h^3$};
            \addplot[slope, domain=0.01890625:0.3025] {0.03926*x^4}
            node[pos=1, anchor=east, inner sep=1pt, font=\scriptsize] {$h^4$};
            \addplot[slope, domain=0.075625:0.3025] {0.005712*x^5}
            node[pos=1, anchor=east, inner sep=1pt, font=\scriptsize] {$h^5$};
        \end{loglogaxis}
    \end{tikzpicture}\hfill
    \begin{tikzpicture}
        \begin{loglogaxis}[num rates, ylabel={$L^2$-norm of the error}]
            \addplot[p2] table[x=h, y=L2_p2_gp, col sep=comma] {\ratesdata};
            \addplot[p3] table[x=h, y=L2_p3_gp, col sep=comma] {\ratesdata};
            \addplot[p4] table[x=h, y=L2_p4_gp, col sep=comma] {\ratesdata};
            \addplot[p5] table[x=h, y=L2_p5_gp, col sep=comma] {\ratesdata};
            \addplot[p2, g1] table[x=h, y=L2_p2_g1, col sep=comma] {\ratesdata};
            \addplot[p3, g1] table[x=h, y=L2_p3_g1, col sep=comma] {\ratesdata};
            \addplot[p4, g1] table[x=h, y=L2_p4_g1, col sep=comma] {\ratesdata};
            \addplot[p5, g1] table[x=h, y=L2_p5_g1, col sep=comma] {\ratesdata};
            \addplot[slope, domain=0.009453125:0.3025] {0.1163*x^3}
            node[pos=0, anchor=west, font=\scriptsize] {$h^3$};
            \addplot[slope, domain=0.009453125:0.3025] {0.01969*x^4}
            node[pos=0, anchor=west, font=\scriptsize] {$h^4$};
            \addplot[slope, domain=0.0378125:0.3025] {0.003955*x^5}
            node[pos=0.4, anchor=north east, inner sep=1pt, font=\scriptsize] {$h^5$};
            \addplot[slope, domain=0.075625:0.3025] {0.0003792*x^6}
            node[pos=1, anchor=east, inner sep=1pt, font=\scriptsize] {$h^6$};
        \end{loglogaxis}
    \end{tikzpicture}
    \caption{Convergence of the nonsymmetric discretization on the unit disc at two ghost weights, each held fixed
        during refinement: solid lines with filled markers the reduced weight $\gamma_p$ of~\eqref{eq:reduced-weight},
        dashed lines with open markers $\gamma=1$. Left: $H^1$-seminorm of the error, with reference slopes $h^p$;
        right: $L^2$ norm of the error, with reference slopes $h^{p+1}$ (dotted). Where the error reaches the
        algebraic error left by the solver, the curves level off.}
    \label{fig:num-rates}
\end{figure}

For each degree $p=1,\dots,5$ the mesh is refined at the weights $\gamma\in\{1,10^{-2},10^{-4}\}$ and at the reduced
weight $\gamma_p$ of~\eqref{eq:reduced-weight}, each held fixed during refinement. Figure~\ref{fig:num-rates} shows the
errors at $\gamma=1$ and at $\gamma_p$, and Table~\ref{tab:num-eoc} lists the observed orders between levels $4$ and
$5$ for all weights. On finer levels the error of the high degrees approaches the algebraic error left by the solver
and stalls, which Figure~\ref{fig:num-rates} shows for $\Q_5$.

Theorem~\ref{thm:energy} gives the order $p$ in the energy error at every fixed positive weight. The weight enters the
constant twice, through the flux-control constant $M_\gamma$, which may grow as the weight decreases, and through the
ghost consistency term, which shrinks (Remark~\ref{rem:weight-once}). Reducing the weight may therefore move the curves
vertically, in either direction, but should not change their slope. For the smooth solution used here, the expected
$L^2$ order is $p+1$, above the proved $p+1/2$ (Remark~\ref{rem:halforder}).

At the reduced weight the curves follow the expected slopes from the coarsest levels on, and the observed orders in
Table~\ref{tab:num-eoc} match $p$ and $p+1$. Reducing the weight moves the curves down: on every level before the stall
the error at $\gamma_p$ lies below the error at $\gamma=1$, by a factor that grows with the degree. The gap narrows
under refinement, since the large weight raises the error mostly on coarse meshes. While it closes, the observed orders
at $\gamma=1$ exceed $p$ and $p+1$; at the intermediate weights they do so only for the higher degrees.

\begin{table}[htbp]
    \centering
    \caption{Observed orders of convergence of the nonsymmetric discretization between levels $\ell=4$ and
        $\ell=5$, $H^1$-seminorm ($H^1$) and $L^2$ norm ($L^2$) of the error, at three
        fixed ghost weights and at the reduced weight $\gamma_p$ of~\eqref{eq:reduced-weight}, each held fixed during
        refinement. Last row: the orders expected on smooth data.}
    \label{tab:num-eoc}
    \setlength{\tabcolsep}{4pt}
    \begin{tabular}{ccccccccc}
        \toprule
        $\gamma$ & \multicolumn{2}{c}{$1$}        & \multicolumn{2}{c}{$10^{-2}$} & \multicolumn{2}{c}{$10^{-4}$}
                 & \multicolumn{2}{c}{$\gamma_p$}                                                                                                              \\
        \cmidrule(lr){2-3}\cmidrule(lr){4-5}\cmidrule(lr){6-7}\cmidrule(lr){8-9}
                 & $H^1$                          & $L^2$                         & $H^1$                         & $L^2$  & $H^1$  & $L^2$  & $H^1$  & $L^2$  \\
        \midrule
        $\Q_1$   & $1.04$                         & $2.07$                        & $0.98$                        & $1.99$ & $0.98$ & $1.98$ & $0.99$ & $2.03$ \\
        $\Q_2$   & $2.36$                         & $3.41$                        & $1.98$                        & $2.97$ & $1.97$ & $2.95$ & $1.97$ & $2.95$ \\
        $\Q_3$   & $3.47$                         & $4.68$                        & $3.09$                        & $4.10$ & $2.98$ & $3.96$ & $2.97$ & $3.96$ \\
        $\Q_4$   & $4.46$                         & $5.44$                        & $4.48$                        & $5.47$ & $3.99$ & $5.00$ & $3.96$ & $4.96$ \\
        $\Q_5$   & $5.73$                         & $6.77$                        & $5.29$                        & $6.27$ & $5.22$ & $6.10$ & $4.97$ & $5.94$ \\
        \midrule
        expected & $p$                            & $p+1$                         & $p$                           & $p+1$  & $p$    & $p+1$  & $p$    & $p+1$  \\
        \bottomrule
    \end{tabular}
\end{table}
 
\subsection{A potato-shaped domain}
\label{sec:num-potato}

\begin{figure}[htbp]
    \centering
    \def\potatodata{numerical/results/nonsymm_paper_plots/potato.csv}
    \begin{subfigure}[b]{0.36\textwidth}
        \centering
        \includegraphics[width=\linewidth]{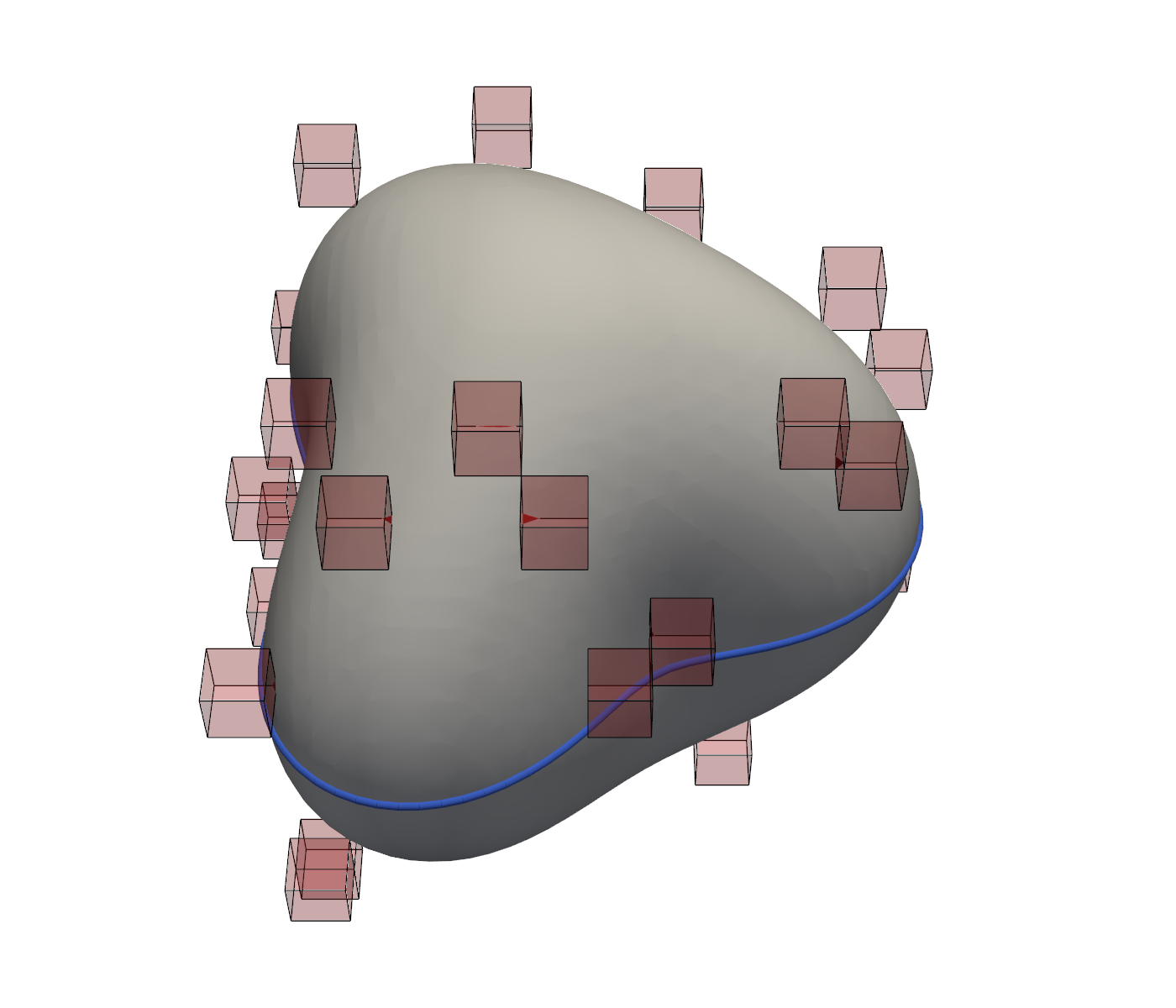}
        \caption{The potato domain~\eqref{eq:potato} and the cells it cuts to less than $10^{-3}$ of their volume.}
    \end{subfigure}\hfill
    \begin{subfigure}[b]{0.6\textwidth}
        \centering
        \pgfplotsset{potato/.style={width=0.85\textwidth, height=0.62\textwidth, xmin=2.7, xmax=6.3,
                    xtick={3,...,6}, every axis/.append style={font=\small}, thick, unbounded coords=jump},
            sym/.style={densely dashed, mark options={solid, fill=white}},
            slope/.style={black!70, dotted, thick, samples=2, forget plot}}
        \begin{tikzpicture}
            \begin{semilogyaxis}[potato, axis y line*=right, axis x line=none,
                    ylabel={smallest cut fraction $\kappa_{\min}$}, ymin=1e-13, ymax=1e-2]
                \addplot[gray, line width=1.5pt, opacity=0.35, mark=*, mark size=0.6pt, mark options={draw=none, fill=gray!70!black, opacity=1}]
                table[x=level, y=kappa_q2, col sep=comma] {\potatodata};
            \end{semilogyaxis}
            \begin{semilogyaxis}[potato, grid=major, xlabel={level $\ell$},
                    ylabel={$H^1$-seminorm of the error}, ymin=1e-5, ymax=1e0,
                    legend pos=north east, legend cell align=left,
                    legend style={font=\scriptsize}, mark size=1.5pt]
                \addplot[p2] table[x=level, y=H1_p2, col sep=comma] {\potatodata}; \addlegendentry{$\Q_2$}
                \addplot[p3] table[x=level, y=H1_p3, col sep=comma] {\potatodata}; \addlegendentry{$\Q_3$}
                \addplot[p4] table[x=level, y=H1_p4, col sep=comma] {\potatodata}; \addlegendentry{$\Q_4$}
                \addlegendimage{gray, line width=1.5pt, opacity=0.35, mark=*, mark size=0.6pt, mark options={draw=none, fill=gray!70!black, opacity=1}}\addlegendentry{$\kappa_{\min}$}
                \addplot[p2, sym, forget plot] table[x=level, y=H1s_p2, col sep=comma] {\potatodata};
                \addplot[p3, sym, forget plot] table[x=level, y=H1s_p3, col sep=comma] {\potatodata};
                \addplot[p4, sym, forget plot] table[x=level, y=H1s_p4, col sep=comma] {\potatodata};
                \addplot[slope, domain=3:6] {0.5274*(2.42*2^(-x))^2} node[pos=1, anchor=west, font=\scriptsize] {$h^2$};
                \addplot[slope, domain=3:5] {0.1347*(2.42*2^(-x))^3} node[pos=1, anchor=west, font=\scriptsize] {$h^3$};
                \addplot[slope, domain=3:4] {0.03162*(2.42*2^(-x))^4} node[pos=1, anchor=west, font=\scriptsize] {$h^4$};
            \end{semilogyaxis}
        \end{tikzpicture}
        \caption{$H^1$-seminorm of the error (left axis), with reference slopes $h^p$ (dotted), and
            $\kappa_{\min}$ (right axis, thick grey line).}
    \end{subfigure}
    \caption{The three-dimensional potato domain. (a) The computed domain for $\Q_3$ on level $4$ and its section at
        $x_3=0.02$ (blue). (b) $H^1$-seminorm of the error against the level, for the nonsymmetric form at the reduced
        weight $\gamma_p$ of~\eqref{eq:reduced-weight} (solid lines with filled markers) and for the symmetric
        form $a_\gamma^{\mathrm{sym}}$ at $\gamma=10^{-1}$ (dashed lines with open markers), and the smallest cut fraction
        $\kappa_{\min}$ of each level, for the level set interpolated in $\Q_2$.}
    \label{fig:num-potato}
\end{figure}
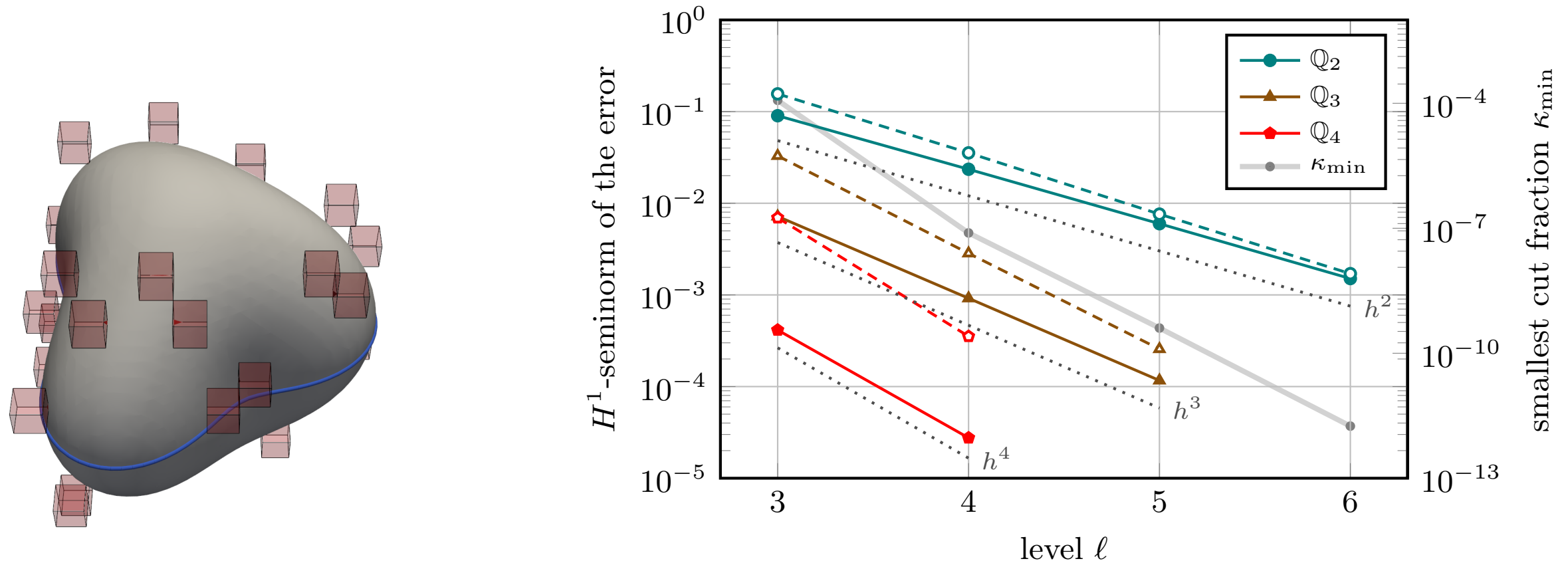

To test the discretization in three dimensions and on a nonconvex boundary of varying curvature, we repeat the
convergence run on a \emph{potato domain}, a deformed ball illustrated in Figure~\ref{fig:num-potato}(a) and given by
the level-set function
\begin{equation}\label{eq:potato}
    \begin{gathered}
        \Omega=\{x\in\mathbb R^3:\ \psi(x)<0\},
        \\
        \psi(x)=|x-c|-r\Bigl(\frac{x-c}{|x-c|}\Bigr),
        \\
        r(\xi)=0.8\bigl(1+0.25\,(\xi_1^3-3\xi_1\xi_2^2)+0.15\,\xi_2\xi_3+0.1\,\xi_3^3\bigr),
    \end{gathered}
\end{equation}
with center $c=(0.03,0.05,0.02)$. The domain is embedded in the background cube $[-1.21,1.21]^3$, and the
discretization is that of Section~\ref{sec:num-setting}, with the manufactured solution
$u=\sin(\pi x_1)\cos(\pi x_2)\cos(\pi x_3)$. The boundary combines convex and concave parts of varying
curvature. The linear systems are solved by GMRES preconditioned with the multigrid method of
Section~\ref{sec:multigrid}, to a residual reduction of $10^{-11}$.

In three dimensions the cut quadrature of Section~\ref{sec:num-setting} needed three changes; without them the error at
high degree was far above the expected one. First, $\psi$ is not smooth at the center $c$, and the cells around $c$
were classified as cut although the zero set does not enter them; cells on which the cut quadrature finds no interface
are therefore integrated as uncut cells. Second, Saye's algorithm~\cite{Saye2015} subdivides a cell when no coordinate
direction represents the zero set as a graph, and after four subdivisions it falls back to a low-order rule; we allow
six. Third, the cut quadrature is built on the $(p+5)$-point Gauss rule instead of the $(p+3)$-point rule.

Figure~\ref{fig:num-potato} shows the potato domain and the $H^1$-seminorm of the error against the level for the
nonsymmetric form at the reduced weight $\gamma_p$ of~\eqref{eq:reduced-weight} and, for comparison, for the symmetric
form $a_\gamma^{\mathrm{sym}}$ at $\gamma=10^{-1}$, above its floor~\eqref{eq:sym-floor-nonsymm} for these degrees,
together with the smallest cut fraction $\kappa_{\min}$ of each level. Both forms are run on $\Q_2$ on levels $3$--$6$,
$\Q_3$ on levels $3$--$5$ and $\Q_4$ on levels $3$ and $4$; for $\Q_4$, level $5$ exceeds the memory of the test
machine. The smallest cut fraction drops by several orders of magnitude from level to level, and
Theorem~\ref{thm:energy} predicts that the rate does not depend on it: at a fixed weight, the cut enters only the
constant, through $C_{\mathrm{Nik}}(p,\mathcal G)$ and the saturation level $1+R$.

The behavior observed in one and two dimensions carries over to three dimensions. At the reduced weight the error
follows the slope $h^p$ from the coarsest level on, although the smallest cut fraction falls by about eight orders of
magnitude over the levels. The symmetric form at $\gamma=10^{-1}$ has a larger error on coarse meshes, by a factor that
grows with the degree, and its higher observed orders reflect this gap closing under refinement.

\subsection{Weights below the symmetric floor}
\label{sec:num-floor}

For comparison we use the symmetric form $a_\gamma^{\mathrm{sym}}$ of~\eqref{eq:aform}, with the load functional
changed accordingly.

The symmetric form is stable uniformly in the cut only above a floor. In the two-cell model of~\mgpaper{}, which is the
model of Figure~\ref{fig:weight-path}, at vanishing cut fraction it is positive semidefinite exactly when
\begin{gather}\label{eq:sym-floor-nonsymm}
    \gamma\ \ge\ \gamma_{\mathrm{lim}}(p) := \frac{1}{h\tau_h-p^2},
    \qquad h\tau_h>p^2,
\end{gather}
and for every $\gamma<\gamma_{\mathrm{lim}}(p)$ it is indefinite at all sufficiently small cut fractions; with the
scaling~\eqref{eq:tau}, $\gamma_{\mathrm{lim}}(p)=1/(\tau_0p(p+1)-p^2)$, which equals $1/(4p^2+5p)$ at $\tau_0=5$. On a given mesh the cut fractions are not
vanishing, and the assembled symmetric matrix loses definiteness at some mesh-dependent weight below
$\gamma_{\mathrm{lim}}(p)$; near a weight at which it becomes singular, the discrete solution is in general unbounded. Whether and
where this happens depends on the cut fractions that happen to occur, and refinement eventually produces cut fractions
at which it does.

\begin{figure}[htbp]
    \centering
    \def\sweepnonsym{numerical/results/nonsymm_paper_plots/sweep_nonsymm.csv}
    \def\sweepsym{numerical/results/nonsymm_paper_plots/sweep_sym.csv}
    \pgfplotsset{num sweep/.style={paperplot, xlabel={ghost weight $\gamma$}, x dir=reverse,
                xmin=1e-12, xmax=1, enlarge x limits=false, xtick={1e0,1e-4,1e-8,1e-12},
                legend pos=north east, legend columns=2, legend cell align=left, legend style={font=\scriptsize},
                width=0.47\columnwidth, mark size=1.5pt},
        nsym/.style={mark repeat=32},
        sym/.style={densely dashed, mark repeat=32, mark phase=17, mark options={solid, fill=white}, forget plot},
        floor/.style={solid, semithick, opacity=0.5, forget plot, no markers},
        floor mark/.style={only marks, mark options={solid, fill=white}, forget plot}}
    \begin{tikzpicture}
        \begin{loglogaxis}[num sweep, ylabel={$H^1$-seminorm of the error}, ymin=1e-10, ymax=1]
            \addplot[p2, nsym] table[x=gamma, y=H1_p2, col sep=comma] {\sweepnonsym}; \addlegendentry{$\Q_2$}
            \addplot[p3, nsym] table[x=gamma, y=H1_p3, col sep=comma] {\sweepnonsym}; \addlegendentry{$\Q_3$}
            \addplot[p4, nsym] table[x=gamma, y=H1_p4, col sep=comma] {\sweepnonsym}; \addlegendentry{$\Q_4$}
            \addplot[p5, nsym] table[x=gamma, y=H1_p5, col sep=comma] {\sweepnonsym}; \addlegendentry{$\Q_5$}
            \addplot[p2, sym] table[x=gamma, y=H1_p2, col sep=comma] {\sweepsym};
            \addplot[p3, sym] table[x=gamma, y=H1_p3, col sep=comma] {\sweepsym};
            \addplot[p4, sym] table[x=gamma, y=H1_p4, col sep=comma] {\sweepsym};
            \addplot[p5, sym] table[x=gamma, y=H1_p5, col sep=comma] {\sweepsym};
            \addplot[teal, floor] coordinates {(3.85e-2,1e-10) (3.85e-2,1)};
            \addplot[p2, floor mark] coordinates {(3.85e-2,1)};
            \addplot[degQthree, floor] coordinates {(1.96e-2,1e-10) (1.96e-2,1)};
            \addplot[p3, floor mark] coordinates {(1.96e-2,0.35)};
            \addplot[red, floor] coordinates {(1.19e-2,1e-10) (1.19e-2,1)};
            \addplot[p4, floor mark] coordinates {(1.19e-2,1)};
            \addplot[violet, floor] coordinates {(8e-3,1e-10) (8e-3,1)};
            \addplot[p5, floor mark] coordinates {(8e-3,0.35)};
        \end{loglogaxis}
    \end{tikzpicture}\hfill
    \begin{tikzpicture}
        \begin{loglogaxis}[num sweep, ylabel={$L^2$-norm of the error}, ymin=1e-13, ymax=1e-3]
            \addplot[p2, nsym] table[x=gamma, y=L2_p2, col sep=comma] {\sweepnonsym};
            \addplot[p3, nsym] table[x=gamma, y=L2_p3, col sep=comma] {\sweepnonsym};
            \addplot[p4, nsym] table[x=gamma, y=L2_p4, col sep=comma] {\sweepnonsym};
            \addplot[p5, nsym] table[x=gamma, y=L2_p5, col sep=comma] {\sweepnonsym};
            \addplot[p2, sym] table[x=gamma, y=L2_p2, col sep=comma] {\sweepsym};
            \addplot[p3, sym] table[x=gamma, y=L2_p3, col sep=comma] {\sweepsym};
            \addplot[p4, sym] table[x=gamma, y=L2_p4, col sep=comma] {\sweepsym};
            \addplot[p5, sym] table[x=gamma, y=L2_p5, col sep=comma] {\sweepsym};
            \addplot[teal, floor] coordinates {(3.85e-2,1e-13) (3.85e-2,1e-3)};
            \addplot[p2, floor mark] coordinates {(3.85e-2,1e-3)};
            \addplot[degQthree, floor] coordinates {(1.96e-2,1e-13) (1.96e-2,1e-3)};
            \addplot[p3, floor mark] coordinates {(1.96e-2,3.5e-4)};
            \addplot[red, floor] coordinates {(1.19e-2,1e-13) (1.19e-2,1e-3)};
            \addplot[p4, floor mark] coordinates {(1.19e-2,1e-3)};
            \addplot[violet, floor] coordinates {(8e-3,1e-13) (8e-3,1e-3)};
            \addplot[p5, floor mark] coordinates {(8e-3,3.5e-4)};
        \end{loglogaxis}
    \end{tikzpicture}
    \caption{Error against the ghost weight on the unit disc at level $\ell=6$, the weight decreasing to the right.
        Solid lines with filled markers: nonsymmetric form~\eqref{eq:aform}; dashed lines with open markers: symmetric
        form $a_\gamma^{\mathrm{sym}}$; the markers are two decades of $\gamma$ apart, those of the symmetric form shifted by one decade. The
        semitransparent vertical line with an open marker of the degree at its top marks $\gamma_{\mathrm{lim}}(p)$
        of~\eqref{eq:sym-floor-nonsymm}. The spikes of the dashed curves occur where an eigenvalue of the symmetric
        matrix passes through zero. Left: $H^1$-seminorm of the error, right: $L^2$ norm of the error.}
    \label{fig:num-sweep}
\end{figure}
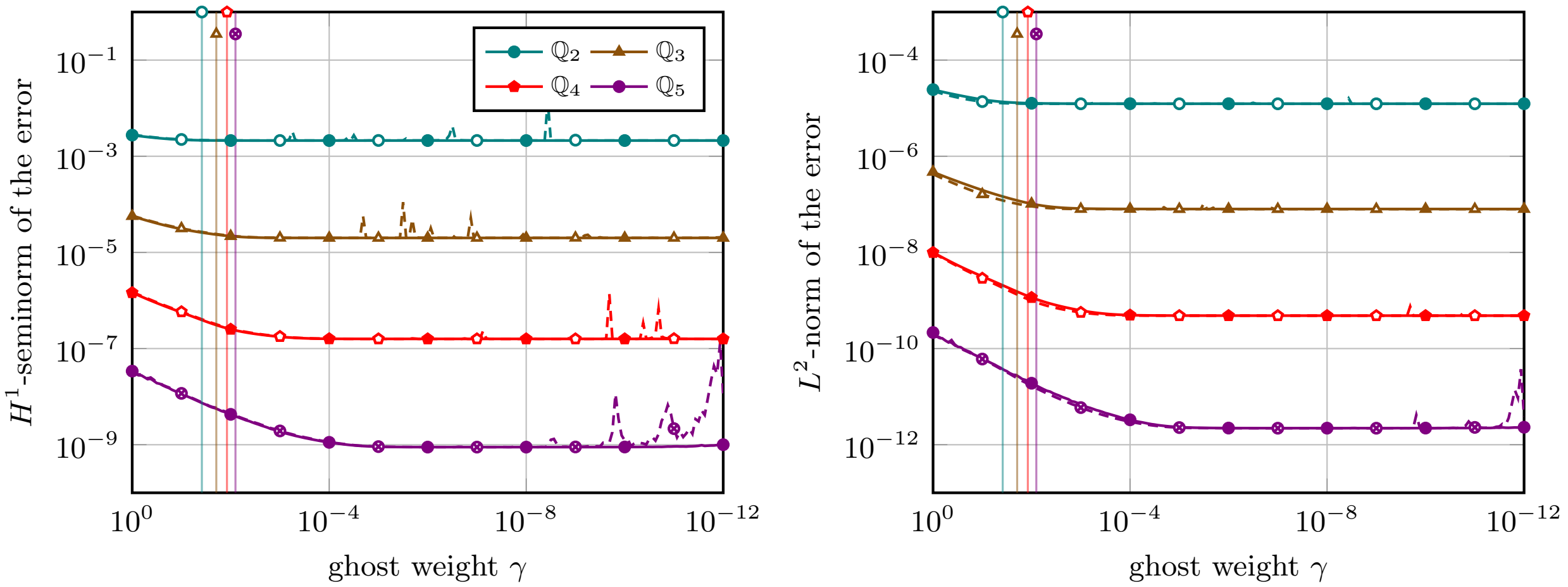

Figure~\ref{fig:num-sweep} fixes the mesh at level $\ell=6$ and sweeps the weight from $1$ to $10^{-12}$ for both
forms, at $193$ weights $\gamma=10^{-k/16}$, $k=0,\dots,192$. The nonsymmetric problem is well posed at every point of
the sweep (Corollary~\ref{cor:wellposed}). The symmetric problem has no stability estimate below
$\gamma_{\mathrm{lim}}(p)$. As the weight decreases, eigenvalues of the symmetric matrix pass through zero; at each
crossing the matrix is singular, and the error in general grows without bound near the crossing weight. In the two-cell
model of Figure~\ref{fig:weight-path} there is a single cut cell and a single crossing at its floor. On the disc the
cut fractions differ from cell to cell, so the crossings occur at several weights below $\gamma_{\mathrm{lim}}(p)$ and
appear as narrow, isolated spikes. Their height depends on how close a sampled weight falls to a crossing. The
eigenvalues of the two-dimensional matrices were not computed; the attribution of the spikes rests on the
one-dimensional model, where the spike coincides with the loss of definiteness. The figure is the two-dimensional
counterpart of Figure~\ref{fig:weight-path}, with the error measured on the whole domain instead of on one cut cell.

The sweep is computed with a direct solver, so that the curves show the discrete solutions themselves. With the
multigrid-preconditioned GMRES the nonsymmetric curves are unchanged up to $\Q_4$, and for $\Q_5$ they are limited by
the algebraic error; the iteration counts are those of Table~\ref{tab:mg-weights}. For the symmetric form below its
floor, GMRES reaches the residual reduction within a few iterations also at nearly singular weights, but the error of
the returned solution can exceed that of the discrete solution by orders of magnitude, and the iteration count gives no
sign of it.
The nonsymmetric error decreases as the weight is reduced and levels off below a weight that depends on the degree;
apart from a slight rise of $\Q_5$ at the smallest weights, it does not grow again anywhere in the sweep, which mirrors
the monotonicity proved in one dimension (Corollary~\ref{cor:onedim}). Above its floor the symmetric form gives nearly
the same error, while below the floor the spikes described above appear for every degree.

\subsection{Smooth data and the half order}
\label{sec:num-halforder}

Theorem~\ref{thm:ltwo} proves the $L^2$ order $p+1/2$, and no lower bound is known. On fitted meshes the loss appears
only at the $H^{p+1}$ endpoint of the regularity scale~\eqref{eq:clz}, while smooth data give $p+1$. To test whether
the same holds on cut meshes, we compare the smooth solution $u$ of Section~\ref{sec:num-setting} with the
boundary-layer solution
\begin{gather}\label{eq:num-boundary-layer}
    u_\sigma=u+\bigl(1-|x|^2\bigr)^{\sigma},
    \qquad
    \sigma=p+\tfrac12+\delta,
    \qquad
    \delta=10^{-4},
\end{gather}
the counterpart on the disc of the solutions with which Chen, Liu and Zhang~\cite{ChenLiuZhang25} show numerically
that~\eqref{eq:clz} is sharp. Its derivatives
of order $p+1$ behave like $\operatorname{dist}(x,\Gamma)^{-1/2+\delta}$ along the whole boundary, so
$u_\sigma\in H^{p+1}(\Omega)$, but $u_\sigma\in W^{p+1,q}(\Omega)$ only for $q<(1/2-\delta)^{-1}$, and~\eqref{eq:clz}
predicts the $L^2$ order $\sigma$. The singular behavior is spread along $\Gamma$ because the loss
in~\eqref{eq:clz} arises in the strip of width $h$ along the boundary. The boundary layer vanishes on $\Gamma$, so the
Dirichlet data are those of $u$. Table~\ref{tab:num-halforder} reports the observed $L^2$ orders for $p=1,\dots,4$ at the
reduced weight $\gamma_p$ of~\eqref{eq:reduced-weight}, together with the symmetric form $a_\gamma^{\mathrm{sym}}$ at $\gamma=1$, above its floor,
as an adjoint-consistent control.
For the smooth solution the nonsymmetric form attains the order $p+1$ for every degree. The boundary-layer solution is
the opposite case,
chosen to be as unfavorable as I could make it: it lies in $H^{p+1}(\Omega)$ only just, and its singular behavior is
spread along the whole of $\Gamma$, exactly where the loss arises. Here the half order does appear. For $\Q_1$ the
nonsymmetric form converges with order $1.5$ on every level, while the symmetric control approaches $2$; for $\Q_2$ the
order decreases toward $2.5$ under refinement, while the control stays close to $3$. The estimate of
Theorem~\ref{thm:ltwo} is therefore sharp, at least for the lowest degrees, but the half order appears only for data
built to provoke it. For $\Q_3$ and $\Q_4$ the run is inconclusive, since the symmetric control itself falls short of
the orders $4$ and $5$ on the finer levels, possibly because of the quadrature of the load, which is not smooth at
$\Gamma$.

\begin{table}[htbp]
    \centering
    \caption{Observed $L^2$ orders for the smooth solution $u$ and the boundary-layer
        solution~\eqref{eq:num-boundary-layer} between levels $\ell-1$ and $\ell$: nonsymmetric form at the reduced
        weight $\gamma_p$, and the symmetric
        form $a_\gamma^{\mathrm{sym}}$ at $\gamma=1$ as adjoint-consistent control. Last row: orders predicted by the
        fitted-mesh estimate~\eqref{eq:clz}, and the optimal order for the symmetric form.}
    \label{tab:num-halforder}
    \setlength{\tabcolsep}{3.5pt}
    \begin{tabular}{ccccccccccccc}
        \toprule
                  & \multicolumn{3}{c}{$\Q_1$}   & \multicolumn{3}{c}{$\Q_2$} & \multicolumn{3}{c}{$\Q_3$}
                  & \multicolumn{3}{c}{$\Q_4$}                                                                                                                                                        \\
        \cmidrule(lr){2-4}\cmidrule(lr){5-7}\cmidrule(lr){8-10}\cmidrule(lr){11-13}
                  & \multicolumn{2}{c}{nonsymm.} & symm.                      & \multicolumn{2}{c}{nonsymm.} & symm.
                  & \multicolumn{2}{c}{nonsymm.} & symm.                      & \multicolumn{2}{c}{nonsymm.} & symm.                                                                                  \\
        $\ell$    & $u$                          & $u_\sigma$                 & $u_\sigma$                   & $u$    & $u_\sigma$ & $u_\sigma$
                  & $u$                          & $u_\sigma$                 & $u_\sigma$                   & $u$    & $u_\sigma$ & $u_\sigma$                                                       \\
        \midrule
        $4$       & $1.94$                       & $1.51$                     & $1.77$                       & $2.93$ & $2.83$     & $3.32$     & $3.96$ & $3.98$ & $4.41$ & $4.94$ & $5.02$ & $4.50$ \\
        $5$       & $2.03$                       & $1.50$                     & $1.91$                       & $2.95$ & $2.70$     & $3.33$     & $3.96$ & $3.84$ & $4.99$ & $4.96$ & $4.97$ & $5.92$ \\
        $6$       & $2.02$                       & $1.49$                     & $2.02$                       & $2.99$ & $2.68$     & $3.28$     & $3.98$ & $3.76$ & $3.92$ & $4.99$ & $4.86$ & $5.90$ \\
        $7$       & $1.99$                       & $1.44$                     & $1.90$                       & $2.99$ & $2.50$     & $2.83$     & $3.99$ & $3.51$ & $3.52$ & $4.99$ & $4.66$ & $4.11$ \\
        $8$       & $2.00$                       & $1.51$                     & $2.00$                       & $3.00$ & $2.60$     & $2.96$     & $4.00$ & $3.66$ & $3.85$ & $4.85$ & $4.59$ & $4.61$ \\
        \midrule
        predicted & $2$                          & $1.5$                      & $2$                          & $3$    & $2.5$      & $3$        & $4$    & $3.5$  & $4$    & $5$    & $4.5$  & $5$    \\
        \bottomrule
    \end{tabular}
\end{table}
 \FloatBarrier

\section{Multigrid with full-residual patches}
\label{sec:multigrid}
This section applies the vertex-patch multigrid method of~\mgpaper{} to the nonsymmetric
system~\eqref{eq:discrete-problem} and tests it by experiment. For the symmetric formulation, \mgpaper{} proves
two-level convergence uniformly in the mesh level and the cut, with a degree dependence that enters through the ghost
weight; that analysis needs coercivity and therefore applies only above the floor. We do not extend it to the
nonsymmetric form. The method itself carries over without modification, and the experiments below measure how its
iteration counts respond to a decreasing weight once the floor is removed.

\subsection{Level hierarchy and patch spaces}

The multigrid method uses a hierarchy $\mesh_0\sqsubset\dots\sqsubset\mesh_L=\mesh_h$ of Cartesian background meshes,
each obtained from its predecessor by uniform dyadic refinement, so that $h_\ell=h_{\ell-1}/2$. All level-$\ell$
objects carry that index: $\mesh_{\ell,\Omega}$, $\Omega_\ell$, $V_\ell$ and $g_\ell$ are the active mesh, the active
domain, the discrete space and the ghost penalty of~\eqref{eq:active}--\eqref{eq:ghost} formed on $\mesh_\ell$, and
$A_\ell$ is the operator of the form~\eqref{eq:aform} on $V_\ell$, with the same weight $\gamma$ on every level. Since
$\Omega_\ell\subseteq\Omega_{\ell-1}$, the natural injection $V_{\ell-1}\to V_\ell$ is well defined and serves as
prolongation; restriction is its adjoint.

The patch centers are every vertex in $\Omega$ together with every exterior vertex having at least one cut adjacent
cell; the exterior vertices must generate their own patches, since otherwise their vertex degrees of freedom are never
smoothed. We enumerate these \emph{admissible} centers $\vertex_1,\dots,\vertex_J$, $J=J(\ell)$. For each of them,
$\widehat\omega_j$ is the full background-grid vertex patch of $\vertex_j$, made of the $2^d$ cells sharing that
vertex, and $\omega_j$ is the union of the active cells in $\widehat\omega_j$. Near $\Gamma$ the patch $\omega_j$ may
consist of fewer than $2^d$ cells; we call such a patch \emph{partial}.

\begin{figure}[!ht]
    \centering
    \tikzset{fr patch/.style={scale=0.95,
                dofD/.style={circle, fill=blue!80!black, inner sep=0pt, minimum size=4pt},
                dofFR/.style={circle, fill=orange, draw=black, inner sep=0pt, minimum size=5pt},
                dofOff/.style={circle, draw=black!60, fill=white, inner sep=0pt, minimum size=3.5pt}}}
    \begin{subfigure}[b]{0.32\textwidth}
        \centering
        \begin{tikzpicture}[fr patch]
            \foreach \x/\y in {0/0,0/1,0/2,1/0,1/1,2/0}
                {\fill[green, opacity=0.3] (\x,\y) rectangle (\x+1,\y+1);}
            \foreach \x/\y in {0/3,1/2,1/3,2/1,2/2,3/0,3/1}
                {\fill[blue, opacity=0.3] (\x,\y) rectangle (\x+1,\y+1);}
            \draw[gray, very thin] (0,0) grid (4,4);
            \draw[line width=1.5pt, red, dashed] (0.6,4) -- (4,0.6);
            \node[red] at (3.7,1.4) {$\Gamma$};
            \node at (0.5,0.5) {$\Omega$};
            \draw[black, line width=1.4pt] (1,1) rectangle (3,3);
            \foreach \i/\j in {1.5/1.5,1.5/2,1.5/2.5,2/1.5,2/2,2/2.5,2.5/1.5,2.5/2,2.5/2.5}
                {\node[dofD] at (\i,\j) {};}
            \foreach \i/\j in {2.5/3,3/2.5,3/3}
                {\node[dofFR] at (\i,\j) {};}
            \foreach \i/\j in {1/1,1/1.5,1/2,1/2.5,1/3,1.5/1,1.5/3,2/1,2/3,2.5/1,3/1,3/1.5,3/2}
                {\node[dofOff] at (\i,\j) {};}
            \draw[red!80!black, line width=1.1pt] (2,2) circle (0.12);
            \node[font=\scriptsize] at (1.75,2.25) {$\vertex_j$};
        \end{tikzpicture}
        \caption{Interior vertex, all patch cells active.}
    \end{subfigure}%
    \hfill
    \begin{subfigure}[b]{0.32\textwidth}
        \centering
        \begin{tikzpicture}[fr patch]
            \foreach \x/\y in {0/0,0/1,1/0,2/0}
                {\fill[green, opacity=0.3] (\x,\y) rectangle (\x+1,\y+1);}
            \foreach \x/\y in {0/2,0/3,1/1,1/2,2/1,3/0,3/1}
                {\fill[blue, opacity=0.3] (\x,\y) rectangle (\x+1,\y+1);}
            \draw[gray, very thin] (0,0) grid (4,4);
            \draw[line width=1.5pt, red, dashed] (0,3.2) -- (4,0.4);
            \node[red] at (3.75,0.95) {$\Gamma$};
            \node at (0.5,0.5) {$\Omega$};
            \draw[black, line width=1.4pt] (1,1) rectangle (3,3);
            \foreach \x/\y in {2/2}
                {\draw[pattern=north east lines, pattern color=gray] (\x,\y) rectangle (\x+1,\y+1);}
            \foreach \i/\j in {1.5/1.5,1.5/2,1.5/2.5,2/1.5,2/2,2/2.5,2.5/1.5,2.5/2}
                {\node[dofD] at (\i,\j) {};}
            \foreach \i/\j in {1.5/3,2/3}
                {\node[dofFR] at (\i,\j) {};}
            \foreach \i/\j in {1/1,1/1.5,1/2,1/2.5,1/3,1.5/1,2/1,2.5/1,3/1,3/1.5,3/2}
                {\node[dofOff] at (\i,\j) {};}
            \draw[red!80!black, line width=1.1pt] (2,2) circle (0.12);
            \node[font=\scriptsize] at (1.75,2.25) {$\vertex_j$};
        \end{tikzpicture}
        \caption{Exterior vertex, partial patch.}
    \end{subfigure}%
    \hfill
    \begin{subfigure}[b]{0.32\textwidth}
        \centering
        \begin{tikzpicture}[fr patch]
            \foreach \x/\y in {0/0}
                {\fill[green, opacity=0.3] (\x,\y) rectangle (\x+1,\y+1);}
            \foreach \x/\y in {0/1,0/2,1/0,1/1,2/0}
                {\fill[blue, opacity=0.3] (\x,\y) rectangle (\x+1,\y+1);}
            \draw[gray, very thin] (0,0) grid (4,4);
            \draw[line width=1.5pt, red, dashed] (0,2.7) -- (2.7,0);
            \node[red] at (0.6,2.55) {$\Gamma$};
            \node at (0.5,0.5) {$\Omega$};
            \draw[black, line width=1.4pt] (1,1) rectangle (3,3);
            \foreach \x/\y in {1/2,2/1,2/2}
                {\draw[pattern=north east lines, pattern color=gray] (\x,\y) rectangle (\x+1,\y+1);}
            \foreach \i/\j in {1.5/1.5,1.5/2,2/1.5,2/2}
                {\node[dofD] at (\i,\j) {};}
            \foreach \i/\j in {1/1,1/1.5,1/2,1/2.5,1/3,1.5/1,2/1,2.5/1,3/1}
                {\node[dofOff] at (\i,\j) {};}
            \draw[red!80!black, line width=1.1pt] (2,2) circle (0.12);
            \node[font=\scriptsize] at (1.75,2.25) {$\vertex_j$};
        \end{tikzpicture}
        \caption{Exterior vertex, one active patch cell.}
    \end{subfigure}
    \caption{Vertex patches for $\Q_2$ around a vertex $\vertex_j$ (red circle) in three cut configurations. Green cells
        are interior, blue cells are cut, and the thick square is the background patch $\widehat\omega_j$; hatched
        cells in it are inactive. Blue dots: degrees of freedom of the Dirichlet space~\eqref{eq:dirichlet-patch-space}.
        Orange dots: the additional degrees of freedom of the full-residual space~\eqref{eq:full-residual-patch-space},
        whose active support lies inside $\omega_j$. White dots: degrees of freedom on $\partial\widehat\omega_j$ whose
        support reaches active cells outside the patch; neither space updates them. In~(c) the two spaces coincide, and each
        of their degrees of freedom also belongs to the full-residual space of a neighboring patch that contains the
        active cell; the \emph{shy} patches of~\cite{wichrowski2026SBMmg} omit such a patch, while here every patch is
        kept.}
    \label{fig:fr-patch}
\end{figure}
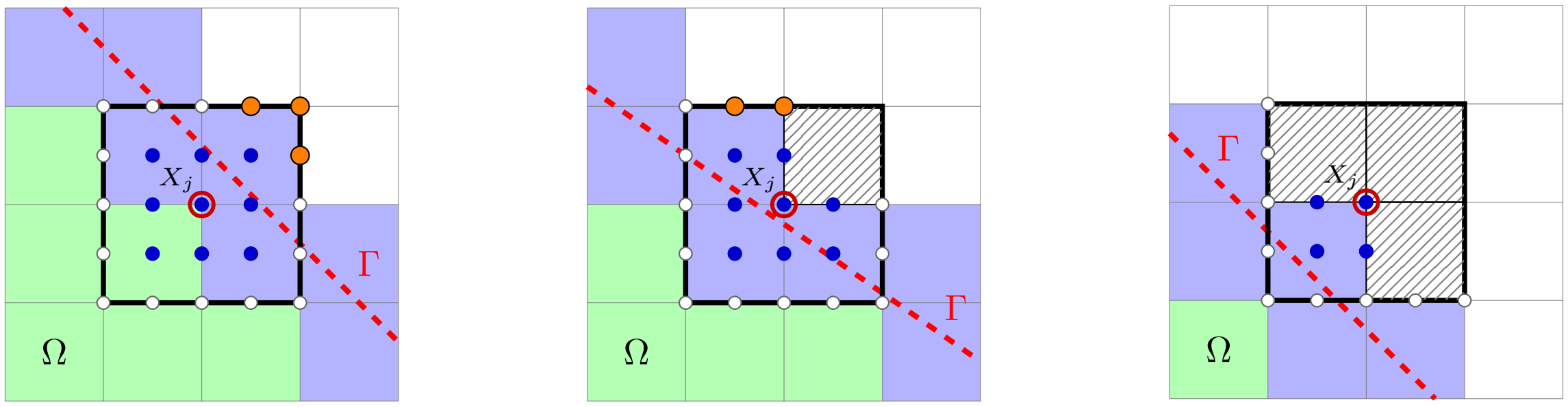

Two local spaces are compared on the same patches. The \emph{Dirichlet} patch space imposes a homogeneous Dirichlet
condition on the whole boundary of the background patch,
\begin{gather}\label{eq:dirichlet-patch-space}
    V_{\ell,j}^{\mathrm D}
    := \operatorname{span}\bigl\{\phi_{\ell,i}:x_{\ell,i}\in
    \operatorname{int}(\widehat\omega_j)\bigr\},
\end{gather}
where $x_{\ell,i}$ is the node of the nodal basis function $\phi_{\ell,i}$. The full-residual
construction~\cite{wichrowski2026SBMmg} instead uses
\begin{gather}\label{eq:full-residual-patch-space}
    V_{\ell,j}^{\mathrm{FR}}
    := \operatorname{span}\bigl\{\phi_{\ell,i}:\omega_i\subseteq\omega_j\bigr\},
    \qquad
    \omega_i:=\bigcup\{\cell\in\mesh_{\ell,\Omega}:\ x_{\ell,i}\in\overline\cell\},
\end{gather}
which updates every degree of freedom whose complete active support is confined to the patch. Since a node interior to
$\widehat\omega_j$ has all its cells in $\widehat\omega_j$, the Dirichlet space is contained in the full-residual space
on every patch. The two differ only where the active domain ends on the boundary of the background patch: a node on
$\partial\widehat\omega_j$ whose neighboring cells outside the patch are all inactive has no active support outside
$\omega_j$. The Dirichlet space excludes such nodes through the artificial Dirichlet condition, while the full-residual
space includes them. Figure~\ref{fig:fr-patch} shows three configurations. The difference arises on partial patches but
also on patches all of whose cells are active, and it disappears when every node on $\partial\widehat\omega_j$ retains
an active cell outside the patch.

\subsection{Smoother and multigrid solver}

Let $V_{\ell,j}$ denote either $V_{\ell,j}^{\mathrm D}$ or $V_{\ell,j}^{\mathrm{FR}}$. Given an iterate $u_\ell$, the
patch correction $s_j\in V_{\ell,j}$ solves
\begin{gather}\label{eq:local-smoother}
    a_\gamma(s_j,v_j)
    = \ell(v_j)-a_\gamma(u_\ell,v_j)
    \qquad\forall v_j\in V_{\ell,j},
\end{gather}
with the global load functional $\ell$ of~\eqref{eq:discrete-problem}, after which $u_\ell\leftarrow u_\ell+s_j$ and
the residual is recomputed before the next patch. Both variants restrict the same operator and the same global residual;
the term \emph{full residual} refers to retaining all residual components of the confined degrees of freedom, not to a
different global operator. The local problem~\eqref{eq:local-smoother} is uniquely solvable for every $\gamma>0$, every
patch and every cut: by Theorem~\ref{thm:coercivity}, $a_\gamma(v_j,v_j)=\|v_j\|_\gamma^2>0$ for every nonzero
$v_j\in V_{\ell,j}$, so the local matrix has a positive definite symmetric part. For the symmetric formulation, the
same property of the local matrices is guaranteed only above the coercivity floor.

Unique solvability alone would survive $\gamma=0$ on a fixed mesh; the positive weight is kept because it gives the
local problems a lower bound that does not depend on the cut. By~\eqref{eq:twobranch}, $\|\nabla
    v_j\|_{0,\Omega}^2+g_h(v_j,v_j)\le\gamma^{-1}a_\gamma(v_j,v_j)$ on every patch and for every cut. At $\gamma=0$ only
the branch $1+R$ remains, and $R$ diverges as the cut degenerates, so a patch around a small cut would become nearly
singular.

A smoothing step visits every patch once and repeats the patches whose cells meet the cut strip $n_c\ge1$
times~\cite{CuiKanschat25CutFEM}; the repetitions strengthen the smoothing near $\Gamma$ but leave the local
problem~\eqref{eq:local-smoother} unchanged. All experiments use the \emph{multiplicative} sweep, which visits the
patches sequentially, with $n_c=1$. Two variants are better suited to parallel execution. The
\emph{semi-multiplicative} sweep splits the patches into $2^d$ colors of non-overlapping
patches~\cite{WitteArndtKanschat21}, processed simultaneously within a color and sequentially across colors; same-color
patches still couple through the ghost faces. The \emph{additive} sweep processes all patches
simultaneously~\cite{Pavarino94additive}. Since the smoother is not symmetrized and the operator is nonsymmetric, the
multigrid cycle preconditions GMRES.

\subsection{Experimental setting}

The two-dimensional experiments use the unit disc and the discretization of Section~\ref{sec:num-setting}; the
three-dimensional ones use the unit ball in $[-1.21,1.21]^3$, with the level-set function $|x|-1$, the manufactured
solution $u=\sin(\pi x_1)\cos(\pi x_2)\cos(\pi x_3)$, and the cut quadrature of Section~\ref{sec:num-potato}. The
coarsest background mesh has two cells in each coordinate direction, and every finer level is obtained by uniform
dyadic refinement, so the finest level $L$ of the hierarchy is the level $\ell=L+1$ of Section~\ref{sec:numerics}. The
operators are assembled as sparse matrices, the ghost penalty~\eqref{eq:ghost} from the tensor-product factorization
of~\cite{wichrowski2026TensorGhostPenalty}, which also admits a matrix-free evaluation as
in~\cite{CuiKanschat25CutFEM}. The cycle is a V-cycle with one pre- and one post-smoothing step, unit relaxation and an
exact solve on the coarsest level. Every run starts from the zero vector and stops when the Euclidean residual norm has
been reduced by a factor $10^{-11}$ relative to the right-hand side, with a cap of $100$ iterations. GMRES is
right-preconditioned, so the criterion applies to the true residual, and is not restarted within the cap.

We report iteration counts and no run times; the implementation used here assembles sparse matrices and is not
optimized for speed. The nonsymmetric form differs from the symmetric one in the sign of one boundary term, so the
operator has the same sparsity and the same quadrature. On uncut patches the boundary terms are absent and the local
problems coincide with those of the symmetric form, which the GPU implementation of~\cite{CuiKanschat25CutFEM} solves
by fast diagonalization, with the residual computed locally in the same
sweep~\cite{WitteArndtKanschat21,WichrowskiMunchKronbichlerKanschat25}. On Cartesian meshes this smoother runs on the
GPU in mixed precision at close to the memory bandwidth~\cite{Wichrowski26MixedPrecisionMG}, which speeds up the uncut
patches, where most of the work lies. On cut patches the implementation of~\cite{CuiKanschat25CutFEM} already uses a
direct solver, which does not rely on symmetry, and the cycle already preconditions GMRES. The cost per iteration is
therefore that of~\cite{CuiKanschat25CutFEM}, where the throughput of the whole solver falls as the iteration count
grows.

Besides fixed weights, the runs use the reduced weight $\gamma_p$ of~\eqref{eq:reduced-weight}, with the values of
Table~\ref{tab:reduced-weight}.

\paragraph{Patch spaces}
Table~\ref{tab:mg-patches} compares the two local spaces on the same discrete problem. The variants share the operator,
the residual, the patch centers and the sweep, so the comparison isolates the effect of the local space.
The full-residual space never needs more iterations than the Dirichlet space and, from $\Q_3$ on, saves between a
quarter and a half of them at both weights; we therefore use it in all remaining experiments.

\begin{table}[htbp]
    \centering
    \caption{GMRES iteration counts for the Dirichlet and the full-residual patch spaces, nonsymmetric form, circle, $L=5$,
        multiplicative sweep with $n_c=1$. Upper block: fixed weight $\gamma=10^{-1}$; lower block: reduced weight
        $\gamma_p$ of~\eqref{eq:reduced-weight}.}
    \label{tab:mg-patches}
    \begin{tabular}{llccccccc}
        \toprule
                                          & local patch space                         & $\Q_1$ & $\Q_2$ & $\Q_3$ & $\Q_4$ & $\Q_5$ & $\Q_6$   & $\Q_7$   \\
        \midrule
        \multirow{2}{*}{$\gamma=10^{-1}$} & Dirichlet, $V_{\ell,j}^{\mathrm D}$       & $9$    & $10$   & $19$   & $41$   & $90$   & $>\!100$ & $>\!100$ \\
                                          & full residual, $V_{\ell,j}^{\mathrm{FR}}$ & $8$    & $7$    & $12$   & $24$   & $52$   & $>\!100$ & $>\!100$ \\
        \midrule
        \multirow{2}{*}{$\gamma_p$}       & Dirichlet, $V_{\ell,j}^{\mathrm D}$       & $9$    & $6$    & $7$    & $7$    & $8$    & $8$      & $10$     \\
                                          & full residual, $V_{\ell,j}^{\mathrm{FR}}$ & $9$    & $5$    & $4$    & $4$    & $5$    & $6$      & $7$      \\
        \bottomrule
    \end{tabular}
\end{table}

\paragraph{Reducing the weight}
Table~\ref{tab:mg-weights} reports the iteration counts as the ghost weight decreases, for the nonsymmetric form
$a_\gamma$ and, for comparison, for the symmetric form $a_\gamma^{\mathrm{sym}}$ of~\eqref{eq:aform}. Every column of
the nonsymmetric block is covered by Theorem~\ref{thm:coercivity}.

The symmetric block is restricted to the weights $\gamma\ge\gamma_{\mathrm{lim}}(p)$ of~\eqref{eq:sym-floor-nonsymm},
at which the symmetric form can be stable uniformly in the cut (Section~\ref{sec:num-floor}). Below this floor, GMRES
requires only nonsingularity and may still return a count, but that count measures the solution of a discrete problem
without a stability estimate rather than the solver for a well-posed discretization. We therefore omit these entries
and mark them with ---. The nonsymmetric form has no such restriction, which is the comparison the table is meant to
show.

For both forms the iteration count grows with the degree at a large weight and falls as the weight decreases, until it
levels off at a small count that depends little on the degree. At equal weights the two forms need nearly the same
number of iterations; they differ in the weights they admit. The symmetric form stops at its floor, where the count
still grows with the degree and exceeds the cap for $\Q_7$, whereas the nonsymmetric form reaches the reduced weight,
where the count has leveled off for every degree. For $\Q_7$ the reduced weight lies below the smallest tabulated
weight and the count is still falling at the end of the table; at $\gamma_p$ itself it drops to the level of the lower
degrees (Table~\ref{tab:mg-patches}).

The growth with the degree reflects the spectrum of the ghost penalty. To illustrate this, we take two uncut cells
sharing one ghost face, compute the largest eigenvalue $\lambda_{\max}(p)$ of the ghost penalty~\eqref{eq:ghost-1d}
relative to the $H^1$-seminorm, and collect the results in Table~\ref{tab:ghost-spectrum}. The eigenvalue does not
depend on $h$, but it grows exponentially with the degree. One dimension is enough here: on tensor-product cells the
face penalty factors into the one-dimensional penalty in the normal direction and mass matrices along the
face~\cite{wichrowski2026TensorGhostPenalty}, so the same value bounds the quotient in any dimension, with equality for
functions constant along the face, when $h_F$ is the cell width. Weighting the penalty shows the difference between the
two forms. From $\Q_3$ on, at the symmetric floor the weighted penalty $\gamma_{\mathrm{lim}}(p)\lambda_{\max}(p)$
exceeds the stiffness by a factor that grows rapidly with the degree, whereas at the reduced weight
$\gamma_p\lambda_{\max}(p)$ stays below one for every degree.

\begin{table}[htbp]
    \centering
    \caption{GMRES iteration counts as the global ghost weight $\gamma$ decreases, circle, $L=5$, multiplicative sweep with
        $n_c=1$, full-residual patch spaces. Upper block: symmetric Nitsche form $a_\gamma^{\mathrm{sym}}$; the column
        $\gamma_\ast$ lists its floor $\gamma_{\mathrm{lim}}(p)$ of~\eqref{eq:sym-floor-nonsymm}, and entries with
        $\gamma<\gamma_{\mathrm{lim}}(p)$ lie below the cut-uniform coercivity floor and are omitted (---). Lower block:
        nonsymmetric form~\eqref{eq:aform}, elliptic at every weight; the column $\gamma_\ast$ lists the reduced weight
        $\gamma_p$ of~\eqref{eq:reduced-weight}. Bold entries compare the two forms at the weight each is meant to
        be used at: for the symmetric form the smallest tabulated weight not below its floor, the lowest weight it
        admits; for the nonsymmetric form the tabulated weight closest to the reduced weight $\gamma_p$ on a
        logarithmic scale. For $\Q_7$, $\gamma_p$ lies below the range of
        the table; at $\gamma_p$ itself GMRES needs $7$ iterations (Table~\ref{tab:mg-patches}).}
    \label{tab:mg-weights}
    \setlength{\tabcolsep}{4pt}
    \begin{tabular}{llccccccccccc}
        \toprule
                                  &        &                & \multicolumn{10}{c}{$\gamma$}                                                                                                                                     \\
        \cmidrule(l){4-13}
                                  &        & $\gamma_\ast$  & $1$                           & $10^{-1}$     & $10^{-2}$         & $10^{-3}$ & $10^{-4}$    & $10^{-5}$ & $10^{-6}$    & $10^{-7}$ & $10^{-8}$    & $10^{-9}$    \\
        \midrule
        \multirow{7}{*}{symm.}    & $\Q_1$ & $1/9$          & $\mathbf{12}$                 & ---           & ---               & ---       & ---          & ---       & ---          & ---       & ---          & ---          \\
                                  & $\Q_2$ & $1/26$         & $13$                          & $\mathbf{7}$  & ---               & ---       & ---          & ---       & ---          & ---       & ---          & ---          \\
                                  & $\Q_3$ & $1/51$         & $22$                          & $\mathbf{13}$ & ---               & ---       & ---          & ---       & ---          & ---       & ---          & ---          \\
                                  & $\Q_4$ & $1/84$         & $50$                          & $\mathbf{25}$ & ---               & ---       & ---          & ---       & ---          & ---       & ---          & ---          \\
                                  & $\Q_5$ & $1/125$        & $>\!100$                      & $55$          & $\mathbf{32}$     & ---       & ---          & ---       & ---          & ---       & ---          & ---          \\
                                  & $\Q_6$ & $1/174$        & $>\!100$                      & $>\!100$      & $\mathbf{72}$     & ---       & ---          & ---       & ---          & ---       & ---          & ---          \\
                                  & $\Q_7$ & $1/231$        & $>\!100$                      & $>\!100$      & $\mathbf{>\!100}$ & ---       & ---          & ---       & ---          & ---       & ---          & ---          \\
        \midrule
        \multirow{7}{*}{nonsymm.} & $\Q_1$ & $2/9$          & $11$                          & $\mathbf{8}$  & $9$               & $9$       & $9$          & $9$       & $9$          & $9$       & $9$          & $9$          \\
                                  & $\Q_2$ & $2/9\,e^{-4}$  & $12$                          & $7$           & $\mathbf{5}$      & $5$       & $5$          & $5$       & $5$          & $5$       & $5$          & $5$          \\
                                  & $\Q_3$ & $2/9\,e^{-8}$  & $21$                          & $12$          & $8$               & $6$       & $\mathbf{4}$ & $4$       & $4$          & $4$       & $4$          & $4$          \\
                                  & $\Q_4$ & $2/9\,e^{-12}$ & $47$                          & $24$          & $15$              & $10$      & $8$          & $5$       & $\mathbf{4}$ & $4$       & $4$          & $4$          \\
                                  & $\Q_5$ & $2/9\,e^{-16}$ & $>\!100$                      & $52$          & $30$              & $21$      & $15$         & $11$      & $8$          & $6$       & $\mathbf{4}$ & $4$          \\
                                  & $\Q_6$ & $2/9\,e^{-20}$ & $>\!100$                      & $>\!100$      & $69$              & $45$      & $30$         & $21$      & $16$         & $12$      & $9$          & $\mathbf{6}$ \\
                                  & $\Q_7$ & $2/9\,e^{-24}$ & $>\!100$                      & $>\!100$      & $>\!100$          & $98$      & $66$         & $42$      & $32$         & $24$      & $18$         & $13$         \\
        \bottomrule
    \end{tabular}
\end{table}

\begin{table}[htbp]
    \centering
    \caption{Largest eigenvalue $\lambda_{\max}(p)$ of the ghost penalty~\eqref{eq:ghost-1d} relative to the
        $H^1$-seminorm on two uncut cells sharing one ghost face, and its products with the floor
        $\gamma_{\mathrm{lim}}(p)$ of~\eqref{eq:sym-floor-nonsymm} and with the reduced weight $\gamma_p$
        of~\eqref{eq:reduced-weight}.}
    \label{tab:ghost-spectrum}
    \begin{tabular}{lccccccc}
        \toprule
                                                      & $\Q_1$ & $\Q_2$            & $\Q_3$            & $\Q_4$            & $\Q_5$            & $\Q_6$            & $\Q_7$            \\
        \midrule
        $\lambda_{\max}(p)$                           & $2.0$  & $8.0$             & $96$              & $1.5\cdot10^{3}$  & $2.4\cdot10^{4}$  & $5.0\cdot10^{5}$  & $1.2\cdot10^{7}$  \\
        $\gamma_{\mathrm{lim}}(p)\,\lambda_{\max}(p)$ & $0.22$ & $0.31$            & $1.9$             & $17$              & $2.0\cdot10^{2}$  & $2.9\cdot10^{3}$  & $5.1\cdot10^{4}$  \\
        $\gamma_p\,\lambda_{\max}(p)$                 & $0.44$ & $3.3\cdot10^{-2}$ & $7.2\cdot10^{-3}$ & $2.0\cdot10^{-3}$ & $6.1\cdot10^{-4}$ & $2.3\cdot10^{-4}$ & $9.8\cdot10^{-5}$ \\
        \bottomrule
    \end{tabular}
\end{table}

\paragraph{Mesh refinement}
Table~\ref{tab:mg-levels} repeats the level sweep for the nonsymmetric form at the fixed weight $\gamma=10^{-1}$ and at
the reduced weight $\gamma_p$ of~\eqref{eq:reduced-weight}, in two and three dimensions.

As for the symmetric form in~\mgpaper{}, the mesh level has little influence on the iteration counts at both weights.
The counts vary somewhat from level to level and rise mildly on the finest levels; such variation with the cut
configuration is also seen in~\cite{CuiKanschat25CutFEM}, where the counts fall again on finer meshes. The nonsymmetric
formulation leaves this behavior unchanged. What the reduced weight changes is the dependence on the degree: at
$\gamma=10^{-1}$ the counts grow with $p$ on every level, whereas at $\gamma_p$ they stay close together for all
degrees, in two and in three dimensions.

\begin{table}[htbp]
    \centering
    \caption{GMRES iteration counts for the nonsymmetric form on the $d$-dimensional ball across mesh levels and
        polynomial degrees, multiplicative sweep with $n_c=1$, full-residual patch spaces. Left: fixed weight
        $\gamma=10^{-1}$; right: reduced weight $\gamma_p$ of~\eqref{eq:reduced-weight}, listed in
        Table~\ref{tab:reduced-weight}. Entries marked with --- are beyond the memory of the test
        machine.}
    \label{tab:mg-levels}
    \begin{subtable}[t]{0.49\textwidth}
        \centering
        \caption{$\gamma=10^{-1}$.}
        \setlength{\tabcolsep}{4pt}
        \begin{tabular}{clccccc}
            \toprule
                                  & $L$ & $\Q_1$ & $\Q_2$ & $\Q_3$ & $\Q_4$ & $\Q_5$ \\
            \midrule
            \multirow{5}{*}{$2D$} & $4$ & $8$    & $6$    & $10$   & $23$   & $49$   \\
                                  & $5$ & $8$    & $7$    & $12$   & $24$   & $52$   \\
                                  & $6$ & $8$    & $6$    & $11$   & $21$   & $45$   \\
                                  & $7$ & $9$    & $7$    & $12$   & $25$   & $56$   \\
                                  & $8$ & $8$    & $8$    & $15$   & $33$   & $79$   \\
            \midrule
            \multirow{3}{*}{$3D$} & $2$ & $9$    & $7$    & $13$   & $26$   & $49$   \\
                                  & $3$ & $9$    & $10$   & $18$   & $45$   & ---    \\
                                  & $4$ & $10$   & $10$   & $20$   & ---    & ---    \\
            \bottomrule
        \end{tabular}
    \end{subtable}%
    \hfill
    \begin{subtable}[t]{0.49\textwidth}
        \centering
        \caption{Reduced weight $\gamma_p$.}
        \setlength{\tabcolsep}{4pt}
        \begin{tabular}{clccccc}
            \toprule
                                  & $L$ & $\Q_1$ & $\Q_2$ & $\Q_3$ & $\Q_4$ & $\Q_5$ \\
            \midrule
            \multirow{5}{*}{$2D$} & $4$ & $8$    & $5$    & $4$    & $4$    & $4$    \\
                                  & $5$ & $9$    & $5$    & $4$    & $4$    & $5$    \\
                                  & $6$ & $9$    & $5$    & $4$    & $4$    & $4$    \\
                                  & $7$ & $9$    & $5$    & $5$    & $5$    & $5$    \\
                                  & $8$ & $9$    & $5$    & $5$    & $6$    & $7$    \\
            \midrule
            \multirow{3}{*}{$3D$} & $2$ & $9$    & $6$    & $5$    & $7$    & $8$    \\
                                  & $3$ & $8$    & $7$    & $8$    & $10$   & ---    \\
                                  & $4$ & $9$    & $8$    & $8$    & ---    & ---    \\
            \bottomrule
        \end{tabular}
    \end{subtable}
\end{table}

The runs on the potato domain of Section~\ref{sec:num-potato} show the same behavior on a nonconvex boundary in three
dimensions; Table~\ref{tab:mg-potato} lists their iteration counts. For the symmetric form at $\gamma=10^{-1}$ the
counts grow with the degree and mildly with the level, while for the nonsymmetric form at $\gamma_p$ they grow far less
with the degree and stay nearly constant across levels. On this geometry smaller weights lower the counts further, to
$5$ for $\Q_3$ at $\gamma=10^{-8}$ and $6$ for $\Q_4$ at $\gamma=10^{-12}$. %
\begin{table}[htbp]
    \centering
    \caption{GMRES iteration counts on the potato domain~\eqref{eq:potato} across the levels $\ell$ of
        Figure~\ref{fig:num-potato}, multiplicative sweep with $n_c=1$, full-residual patch spaces. Left: symmetric form
        $a_\gamma^{\mathrm{sym}}$ at $\gamma=10^{-1}$; right: nonsymmetric form~\eqref{eq:aform} at the reduced weight
        $\gamma_p$ of~\eqref{eq:reduced-weight}. Entries marked with --- are beyond the memory of the test machine.}
    \label{tab:mg-potato}
    \setlength{\tabcolsep}{4pt}
    \begin{tabular}{lccccccc}
        \toprule
                 & \multicolumn{3}{c}{symm., $\gamma=10^{-1}$} &        & \multicolumn{3}{c}{nonsymm., $\gamma_p$}                                    \\
        \cmidrule(lr){2-4}\cmidrule(l){6-8}
        $\ell$   & $\Q_2$                                       & $\Q_3$ & $\Q_4$                                    &  & $\Q_2$ & $\Q_3$ & $\Q_4$ \\
        \midrule
        $3$      & $10$                                         & $19$   & $35$                                      &  & $7$    & $9$    & $13$   \\
        $4$      & $12$                                         & $23$   & $49$                                      &  & $8$    & $9$    & $13$   \\
        $5$      & $13$                                         & $26$   & ---                                       &  & $8$    & $10$   & ---    \\
        $6$      & $14$                                         & ---    & ---                                       &  & $9$    & ---    & ---    \\
        \bottomrule
    \end{tabular}
\end{table}
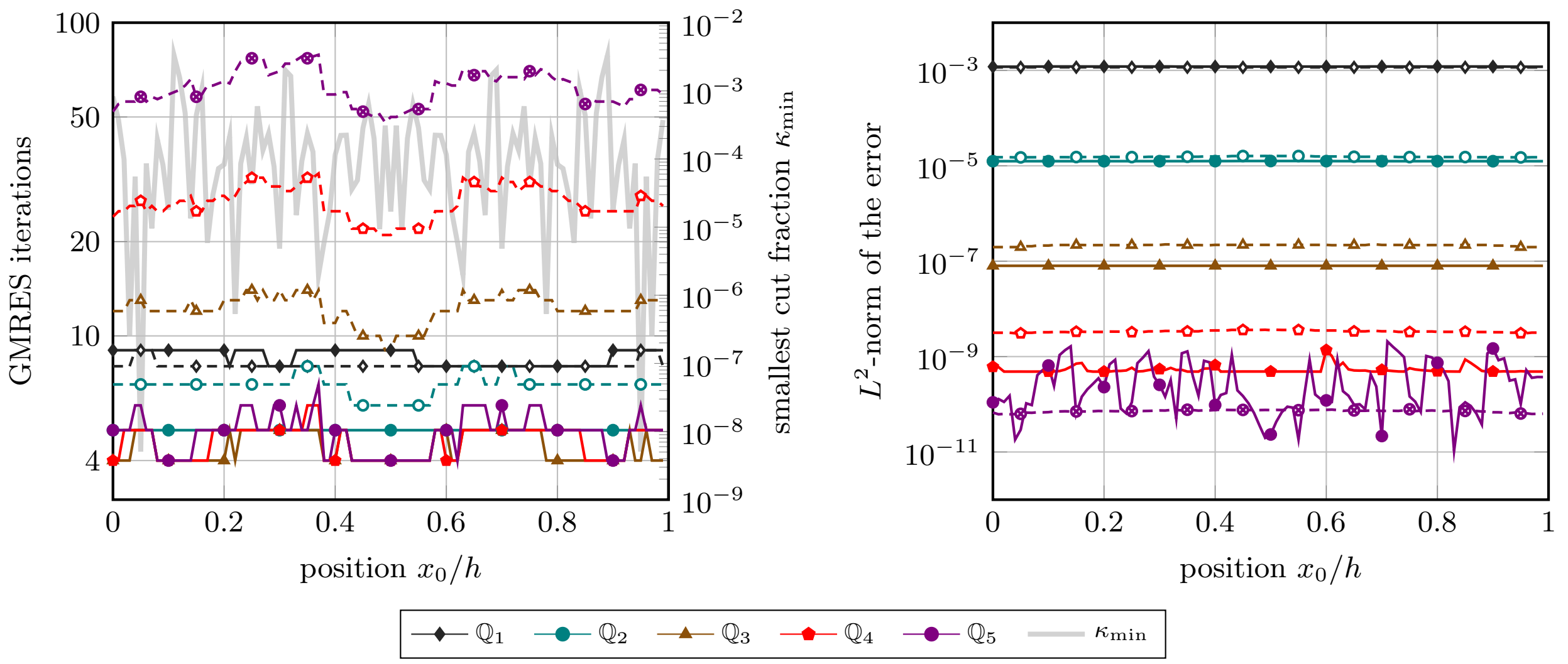
\begin{figure}[htbp]
    \centering
    \def\shiftdata{numerical/results/nonsymm_paper_plots/shift.csv}
    \pgfplotsset{mg shift/.style={paperplot, xlabel={position $x_0/h$}, xmin=0, xmax=1, enlarge x limits=false,
                width=0.45\columnwidth, mark size=1.5pt},
        p1/.style={black!85, mark=diamond*},
        gp/.style={mark repeat=10},
        g01/.style={densely dashed, mark repeat=10, mark phase=6, mark options={solid, fill=white}, forget plot}}
    \begin{tikzpicture}
        \begin{semilogyaxis}[mg shift, axis y line*=right, axis x line=none, grid=none,
                ylabel={smallest cut fraction $\kappa_{\min}$}, ymin=1e-9, ymax=1e-2]
            \addplot[gray, line width=1.5pt, opacity=0.35, no markers] table[x=xh, y=kappa_q3, col sep=comma] {\shiftdata};
        \end{semilogyaxis}
        \begin{semilogyaxis}[mg shift, ylabel={GMRES iterations}, ymin=3, ymax=100,
                ytick={4,10,20,50,100}, yticklabels={$4$,$10$,$20$,$50$,$100$},
                legend to name=mgshiftlegend, legend columns=-1, legend style={font=\scriptsize, /tikz/every even column/.append style={column sep=6pt}}]
            \addplot[p1, gp] table[x=xh, y=it_p1_gp, col sep=comma] {\shiftdata}; \addlegendentry{$\Q_1$}
            \addplot[p2, gp] table[x=xh, y=it_p2_gp, col sep=comma] {\shiftdata}; \addlegendentry{$\Q_2$}
            \addplot[p3, gp] table[x=xh, y=it_p3_gp, col sep=comma] {\shiftdata}; \addlegendentry{$\Q_3$}
            \addplot[p4, gp] table[x=xh, y=it_p4_gp, col sep=comma] {\shiftdata}; \addlegendentry{$\Q_4$}
            \addplot[p5, gp] table[x=xh, y=it_p5_gp, col sep=comma] {\shiftdata}; \addlegendentry{$\Q_5$}
            \addlegendimage{gray, line width=1.5pt, opacity=0.35, no markers}\addlegendentry{$\kappa_{\min}$}
            \addplot[p1, g01] table[x=xh, y=it_p1_g01, col sep=comma] {\shiftdata};
            \addplot[p2, g01] table[x=xh, y=it_p2_g01, col sep=comma] {\shiftdata};
            \addplot[p3, g01] table[x=xh, y=it_p3_g01, col sep=comma] {\shiftdata};
            \addplot[p4, g01] table[x=xh, y=it_p4_g01, col sep=comma] {\shiftdata};
            \addplot[p5, g01] table[x=xh, y=it_p5_g01, col sep=comma] {\shiftdata};
        \end{semilogyaxis}
    \end{tikzpicture}\hfill
    \begin{tikzpicture}
        \begin{semilogyaxis}[mg shift, ylabel={$L^2$-norm of the error}, ymin=1e-12, ymax=1e-2]
            \addplot[p1, gp] table[x=xh, y=L2_p1_gp, col sep=comma] {\shiftdata};
            \addplot[p2, gp] table[x=xh, y=L2_p2_gp, col sep=comma] {\shiftdata};
            \addplot[p3, gp] table[x=xh, y=L2_p3_gp, col sep=comma] {\shiftdata};
            \addplot[p4, gp] table[x=xh, y=L2_p4_gp, col sep=comma] {\shiftdata};
            \addplot[p5, gp] table[x=xh, y=L2_p5_gp, col sep=comma] {\shiftdata};
            \addplot[p1, g01] table[x=xh, y=L2_p1_g01, col sep=comma] {\shiftdata};
            \addplot[p2, g01] table[x=xh, y=L2_p2_g01, col sep=comma] {\shiftdata};
            \addplot[p3, g01] table[x=xh, y=L2_p3_g01, col sep=comma] {\shiftdata};
            \addplot[p4, g01] table[x=xh, y=L2_p4_g01, col sep=comma] {\shiftdata};
            \addplot[p5, g01] table[x=xh, y=L2_p5_g01, col sep=comma] {\shiftdata};
        \end{semilogyaxis}
    \end{tikzpicture}

    \smallskip
    \pgfplotslegendfromname{mgshiftlegend}
    \caption{The circle $|x-(x_0,0)|=1$ moved through one cell of the finest level in steps of $h/100$,
        nonsymmetric form, $L=5$, multiplicative sweep with $n_c=1$, full-residual patch spaces. Solid lines with filled
        markers: reduced weight $\gamma_p$; dashed lines with open markers: $\gamma=10^{-1}$; the markers are ten
        positions apart. Left: GMRES iteration counts, and the smallest cut fraction $\kappa_{\min}$ of the finest level
        (thick grey line, right axis). Right: $L^2$ norm of the error. For $\Q_5$ at $\gamma_p$ the spread of the error is the
        algebraic error left by the solver.}
    \label{fig:mg-shift}
\end{figure}

\paragraph{Small cuts}
The meshes above contain only the cut fractions that the disc happens to produce. To expose the method to arbitrarily
small cuts, we move the circle through one cell of the finest level in steps of $h/100$, replacing the disc by
$|x-(x_0,0)|<1$ with $0\le x_0<h$, on the mesh of Table~\ref{tab:mg-weights}, and solve with the nonsymmetric form at
$\gamma=10^{-1}$ and at $\gamma_p$. Over the $100$ positions the smallest cut fraction of the finest level ranges over
six orders of magnitude. Figure~\ref{fig:mg-shift} shows the iteration counts and the $L^2$ error against the position,
together with the smallest cut fraction. All runs converge. The smallest cut fraction jumps by orders of magnitude
between neighboring positions and falls below $10^{-8}$, where the physical part of the cell is a sliver at the edge of
what the geometry and the quadrature resolve, yet neither the iteration counts nor the errors show any dependence on
it: at the smallest cut both lie within the range met at the other positions. At $\gamma_p$ the counts stay at a few
iterations for every degree. At $\gamma=10^{-1}$ they grow with $p$ as in Table~\ref{tab:mg-weights}, and for the high
degrees they change in steps with the position of the boundary rather than with the cut fraction. The $L^2$ error
depends only weakly on the position of the boundary; for $\Q_5$ at $\gamma_p$ its spread is the algebraic error left by
the solver tolerance, as the discretization error is close to machine precision there.

\FloatBarrier

\section{Conclusions}
\label{sec:conclusions}
Changing the sign of one term of the Nitsche form allows the ghost weight to be reduced as far as one wishes, which
improves both the discretization and the multigrid solver. Most of the degree dependence of multigrid on cut meshes
comes from the ghost penalty. At a reduced weight, it nearly disappears. On cut meshes the Nitsche penalty is no longer
coupled to the size of the intersected cell, so in the symmetric form the ghost penalty must control the boundary flux
on small cuts, which puts a floor under the weight. The nonsymmetric form makes every positive weight admissible: it is
elliptic on every cut configuration and mesh level, whatever the weight. In a norm that does not depend on the weight,
the ellipticity constant shrinks with the weight but never reaches zero, and on a fixed configuration it levels off at
a positive limit once the weight falls below a geometry-determined threshold. Well-posedness and a priori error
estimates hold for every degree-dependent weight held fixed under refinement, with constants that depend on the cut.

This dependence could spoil the error or the solver; neither happens. In the estimates, the cut and the weight act only
through one flux-control constant, which is bounded by $\gamma^{-1/2}$ for every cut and saturates on a fixed
configuration. In the experiments, reducing the weight lowers the error, at the reduced weight by a factor that grows
with the degree, and the error does not grow again as the weight decreases further. Moving the boundary through a cell,
with cut fractions below $10^{-8}$, leaves both the error and the iteration counts within the range met at the other
positions.

The price is adjoint consistency, which costs half an order in the $L^2$ estimate. This half order is real: a
boundary-layer solution concentrated along the boundary shows it for $\Q_1$ and $\Q_2$. Smooth solutions do not notice
it and converge at the optimal order $p+1$. At large weights the observed orders even exceed $p$ and $p+1$, because the
error starts high on coarse meshes and then approaches the error at the reduced weight.

For the multigrid method, the experiments isolate the ghost penalty. The mesh, the degree, the smoother, and the patch
spaces stay fixed; only the weight changes, and most of the growth of the iteration counts with the degree goes with
it, up to $\Q_7$. At a fixed weight this growth follows the spectrum of the ghost penalty, whose largest eigenvalue
relative to the stiffness grows exponentially with the degree. At equal weights the two forms need nearly the same
number of iterations; the sign change matters because it lets the weight go lower. Under mesh refinement the counts
change little, and full-residual patch spaces save a quarter to a half of the iterations from $\Q_3$ on.

Our error estimates are conservative: their constants grow as the weight decreases, at most like $\gamma^{-1/2}$, while
in the experiments the error decreases with the weight and then levels off. Only in one dimension do we prove that the
error decreases with the weight; an estimate that captures this behavior in higher dimensions is missing. The second
question is a convergence theory for multigrid at reduced weights, which would need a stable decomposition below the
symmetric floor and field-of-values bounds for GMRES.

So, how far may the weight be reduced? With the nonsymmetric formulation, apparently quite a lot, and it is worth
considering.

\paragraph{Declaration of assistance} The author declares support of language models (Gemini, Claude, ChatGPT) throughout the preparation of this work.

The manuscript is the result of a long iterative process, in which these models contributed to text drafting,
proof-reading, checks of mathematical consistency, and adversarial review of the statements and their proofs. The paper
was audited by the author, who retains full accountability for all scientific content.

\bibliographystyle{elsarticle-num}

\end{document}